%% file: RMAR_2026_arxive.tex
\documentclass[11pt,a4paper]{article}

\usepackage{header}
\usepackage{abbrev}

\usepackage[normalem]{ulem}

\begin{document}

\title{\algname{}: A Riemannian Multilevel Adaptive Regularization Method with application to fixed-rank problems}

\author{\MakeUppercase{Margherita Porcelli}\thanks{Università degli Studi di Firenze, Dipartimento di Ingegneria Industriale (DIEF), Italy (\email{margherita.porcelli@unifi.it}).  Member of the INdAM Research Group GNCS and of the SIMAI-OPTIMA Activity Group.}
\thanks{ISTI--CNR, Via Moruzzi 1, Pisa, Italia.}\hspace{2mm}\orcidlink{0000-0003-0183-1204}, \MakeUppercase{Elisa Riccietti}\thanks{ENS de Lyon, CNRS, Inria, LIP, UMR 5668, Université Claude Bernard Lyon 1, Lyon, France (\email{elisa.riccietti@ens-lyon.fr}).}\hspace{2mm}\orcidlink{0000-0003-0864-6116}, \MakeUppercase{and Marco Sutti}\thanks{Division of Mathematics, Gran Sasso Science Institute, L'Aquila, Italy (\email{marco.sutti@gssi.it}). Member of the INdAM Research Group GNCS.}\hspace{2mm}\orcidlink{0000-0002-8410-1372}}

\maketitle

\begin{abstract}
We propose the first  Riemannian multilevel adaptive regularization algorithm for solving optimization problems on manifolds endowed with a hierarchical structure. By leveraging regularized first-order models across fine and coarse manifolds, the framework alternates between fine-level tangent space steps and approximated coarse-level model minimizations. We establish a worst-case iteration complexity analysis for general objective functions and retractions, proving convergence bounds to $\varepsilon$-approximate first-order minimizers. Furthermore, we specialize the algorithmic framework to the manifold of fixed-rank matrices. Numerical experiments on large-scale low-rank matrix optimization problems demonstrate that our method outperforms state-of-the-art Riemannian multilevel solvers.



\end{abstract}

\bigskip

\section{Introduction}

In recent years, Riemannian optimization on matrix manifolds, such as the Stiefel, Grassmann, fixed-rank, and Symmetric Positive Definite (SPD) manifolds, has emerged as a powerful paradigm for solving complex constrained optimization problems arising in machine learning, computer vision, and scientific computing \cite{AMS:2004,Vandereycken:2013,Boumal_Absil:2015,agarwal2021adaptive,iannazzo2018riemannian}. By reformulating constrained problems defined in the Euclidean space as unconstrained problems over embedded or quotient submanifolds, Riemannian methods naturally preserve geometric structures and often enjoy superior convergence properties then Euclidean ones; see, e.g., \cite{AMS:2008,boumal2023introduction}.

Just as in the Euclidean setting, large-scale problems often entail substantial computational costs. However, many applications exhibit an intrinsic multiscale or multiresolution structure, either in the underlying optimization domain or in the objective function itself \cite{MariniPorcelliRiccietti2026}. Leveraging this structure offers the potential to significantly reduce computational costs and enhance algorithmic scalability, notably by exploiting the \textit{multilevel paradigm} \cite{Briggs:2000,Trottenberg:2000,Hackbusch:2003,Lewis:2005,WenGoldfarb:2010,eusipco,lauga2024iml,petra_bregman}.

In this work, we propose a multilevel Riemannian optimization framework for problems exhibiting such an intrinsic hierarchical structure. We then specialize this framework to optimization problems defined on \textit{matrix manifolds}, exploiting their geometric and multilevel properties.

We focus on the adaptive regularization paradigm, and we propose \algname{}, a new Riemannian multilevel adaptive regularization gradient method that generalizes the AR1 (Adaptive Regularization with a first-order model) method \cite{book_compl}. Adaptive regularization strategies have been proposed in the Riemannian context by Agarwal et al. \cite{agarwal2021adaptive}, through a Riemannian extension of the Adaptive Regularization with Cubics (ARC) method originally introduced for the Euclidean setting in \cite{book_compl}. Their work demonstrates that the Riemannian formulation preserves optimal iteration complexity bounds, while carefully analyzing the required underlying assumptions and identifying the manifold properties needed to guarantee key theoretical properties. Standard ARC methods employ a second-order Taylor model regularized by a cubic term, categorized as AR2 (Adaptive Regularization with second-order models) by Cartis et al. \cite{book_compl}. In contrast, the present work employs regularized first-order Taylor models, classified as AR1 in \cite{book_compl}, and leverages the hierarchical structure of the underlying problem. \algname{} can thus be interpreted as a Riemannian version of the multilevel AR$q$ method of Calandra et al. \protect{\cite{Calandra:2021}} with $q=1$. 

For the sake of clarity and simplicity of exposition, we first describe our framework within a two-level hierarchy. We assume that the hierarchy involves two distinct manifolds, a \textit{fine manifold} $ \cMh$ and a \textit{coarse manifold} $\cMH$. 
For the matrix manifold case, the manifolds are understood in the sense of classical differential geometry, with a natural representation of elements as matrices \cite{EAS:1998,AMS:2008}; $ \cMh$ is thus defined for $n^h$-by-$n^h$ matrices, while $\cMH$ is defined for $n^H$-by-$n^H$ matrices.

Quantities related to the fine level and to the coarse level are denoted by the superscripts $\cdot^h$ and $\cdot^{H}$, respectively.
The manifolds $\cMh$ and $\cMH$, their tangent bundles $\mathrm{T}\cMh$ and $\mathrm{T}\cMH$, together with two (sufficiently smooth) cost functions $f^{h}:\cMh \mapsto \R$ and
$f^{H} \colon \cMH \mapsto \R$,  define the optimization hierarchy: the fine problem
\begin{equation}\label{eq:pb_fine_manifold}
    \underset{x^{h} \in \cMh}{\min} f^{h}(x^{h}),
\end{equation}
 and its coarser version
\begin{equation}\label{eq:pb_coarse_manifold}
    \underset{x^{H} \in \cMH}{\min} f^{H}(x^{H}).
\end{equation}

Our method alternates between gradient steps for $f^h$ on the fine manifold $\cMh$ and steps that minimize suitable regularized models for $f^H$ on the coarse manifold $\cMH$.  Remarkably, if only fine steps are computed, the method reduces to the Riemannian adaptive first-order regularization method (RAR1), as detailed in \Cref{app:RAR1}.

Inheriting the definition of the coarse model from the class of multilevel methods, e.g., \cite{Trottenberg:2000,Briggs:2000,Hackbusch:2003,Nash:2000,Lewis:2005,WenGoldfarb:2010}, our approach adapts the Euclidean convergence and complexity analysis of multilevel AR algorithms to establish its Riemannian theoretical guarantees; see the multilevel AR$q$ method of Calandra et al. \protect{\cite[Algorithm 3.1]{Calandra:2021}}, with $q=1$. 

To the best of our knowledge, only a few contributions in the literature address general Riemannian multilevel optimization methods, and they all rely on line-search globalization strategies. The main contribution is the paper \cite{sutti2021riemannian}, where Sutti and Vandereycken proposed the Riemannian Multi-Grid Line-Search (RMGLS) method\footnote{The original idea for this algorithm goes back to the PhD thesis of Vandereycken~\cite{Vandereycken_2010}.}, an algorithm combining multilevel optimization of \cite{Nash:2000,WenGoldfarb:2010} with Riemannian optimization on the manifold of fixed-rank matrices \cite{Vandereycken:2013}. Compared with previous approaches based on alternating minimization or rank-truncated iterative methods, RMGLS avoids solving ill-conditioned local systems or using explicit preconditioning \cite{VandereyckenVandewalle:2010,KSV:2016,RakhubaOseledets:2018}, thanks to the intrinsic scalability of multigrid algorithms. Despite its success in computing good low-rank approximations with an almost mesh-independent convergence behavior, the RMGLS method lacked theoretical convergence guarantees. In the recent preprint~\cite{Elshiaty:2026}, Elshiaty et al. introduced a Riemannian extension of multilevel optimization, also based on line search, and on a coarse model that is first-order coherent with the fine-level objective and yields descent directions under mild retraction-convexity assumptions. They formulated two-level and multilevel algorithms and provided a unified Riemannian multilevel optimization framework with a global convergence analysis using a Riemannian Zoutendijk-type argument. While their framework is illustrated through a variety of applications, low-rank matrix optimization problems are not considered.

The primary contributions of our work are summarized as follows:
\begin{itemize}
    \item We introduce \algname{}, the first Riemannian multilevel gradient-type algorithm that combines \textit{adaptive  regularization} techniques with a \textit{hierarchical manifold structure}.

    \item We prove global convergence of \algname{} to stationary points and provide a formal \textit{worst-case iteration complexity} analysis for general minimization problems on Riemannian manifolds. Convergence rate bounds  to $\varepsilon$-approximate first-order minimizers are provided, \Cref{thm worst case complexity \algname{}} being our main theoretical result.
    
    \item We tailor the general theory of \algname{} to the fixed-rank matrix manifold $\mathcal{M}_r$ and specialized the implementation on multiple levels following \cite{sutti2021riemannian}.
    
    \item We conduct extensive numerical experiments illustrating the method and demonstrating that \algname{} consistently outperforms RMGLS, the state-of-the-art Riemannian multilevel solver, on large-scale low-rank matrix optimization problems, with the advantage of being theoretically grounded.

    \item We present a rank-adaptive version of \algname{} for low-rank minimization, useful in practical settings where the rank of the optimization problem at hand is not known a priori.
\end{itemize}

\subsection{Outline of the paper}
This paper is structured as follows. In \Cref{sec:preliminaries}, we provide some preliminaries on Riemannian multilevel optimization. In \Cref{sec:two-level}, we describe our new algorithm, the Riemannian Multilevel Adaptive Regularization First-Order(\algname{}) method, for the solution of problems defined on any two-level hierarchy of manifolds, and we study its theoretical properties in  \Cref{sec:theory}. \Cref{sec:geom_fixed_rank_manif} describes the geometry of fixed-rank matrix manifolds and discusses applying \algname{} in this context. Finally, in~\Cref{sec:numerical_experiments}, we perform extensive numerical experiments using a multilevel implementation of \algname{} for the manifold of fixed-rank matrices and show moderate mesh dependence of \algname{}, comparisons with the one-level version RAR1 and the RMGLS method of~\cite{sutti2021riemannian}, and a rank-adaptive version of \algname{}. Finally, in \Cref{sec:conclusions} we wrap up this work with conclusions and outlook.
\Cref{app:RAR1} presents RAR1, and  \Cref{app:suppl} contains more detailed calculations of the models used in \Cref{sec:two-level} and \Cref{sec:theory} for the case of fixed-rank manifolds.

\section{Preliminaries on Riemannian multilevel optimization}\label{sec:preliminaries}
Let $\cM$ denote either $\cMh$ or $\cMH$, the tangent space to $\mathcal{M}$ in $x$ by $\mathrm{T}_x\mathcal{M}$, the Riemannian inner product on $\mathrm{T}_x\mathcal{M}$ by $\langle\cdot,\cdot\rangle_x$, and its induced norm by $\|\cdot\|_x$. We assume that $ \cMh$ and $ \cMH$ are both endowed with a Riemannian metric. 

Given a smooth function $ f \colon \mathcal{M} \to \mathbb{R} $, we denote its Riemannian gradient by $ \operatorname{grad}f(x) $, characterized by
\[
    \innerproduct{\operatorname{grad}f(x)}{s}_{x} = \operatorname{D} \! f(x)[s], \qquad s \in \mathrm{T}_x\mathcal{M},
\]
where $\operatorname{D} \! f$ is the differential of $f$. 
We use $\nabla$ for gradients of functions defined on Euclidean spaces, including tangent spaces equipped with the fixed inner product at their base point.

Retractions are mappings used to map tangent vectors to the manifold and turn cost functions defined on the manifold into cost functions defined on a tangent space; see, e.g., \cite{AMS:2008,Absil:2012,AbsilOseledets:2015}.
Given a point $ x_k \in \mathcal{M} $ and a smooth retraction $ \operatorname{R} \colon \mathcal{M} \to \R $, we define the \textit{pullback} of $ f $ at $ x_k $ by
\begin{equation}\label{eq:def_pullback}
    \hat{f}_k = f \circ \operatorname{R}_{x_k} \colon \mathcal{U}_k \subseteq \mathrm{T}_{x_k}\mathcal{M} \to \mathbb{R},
\end{equation}
where $ \mathcal{U}_k $ is an open neighborhood of the origin $ 0_{x_k} $, on which the retraction is defined. Thus, the pullback expresses the cost function on a vector space. The retraction properties $ \operatorname{R}_{x_k}(0_{x_k}) = x_k $ and $ \operatorname{D} \! \operatorname{R}_{x_k}(0_{x_k}) = \operatorname{Id}_{\mathrm{T}_{x_k}\mathcal{M}} $ (the latter known as local rigidity condition) imply
\[
    \nabla\hat{f}_k(0_{x_k}) = \operatorname{grad}f(x_k).
\]
We typically retain the subscript $k$ on all pullbacks, as a reminder that they depend on the base point $ x_{k} $. In particular, at level $ \ell \in \{ h, H \} $, we will write $ \hat{f}_{k}^{\ell} = f^{\ell} \circ \operatorname{R}_{x_{k}^{\ell}} $.

\algname{} uses point and tangent transfers to connect the fine and coarse manifolds in~\eqref{eq:pb_fine_manifold} and~\eqref{eq:pb_coarse_manifold}. The point restriction $\mathcal{I}_h^H \colon \cMh \to \cMH$ and prolongation $\mathcal{I}_H^h \colon \cMH \to \cMh$ are assumed to be well-defined at all points considered by the algorithm. At iteration $k$, the coarse base point is computed as
\[
	x_k^H = \mathcal{I}_h^H(x_k^h) = I_m \, x_k^h \, I_n^{\top},
\]
where $ I_m $ and $ I_n $ are the one-dimensional injection matrices.

Tangent vectors are transferred by the linear operators
\[
    \widetilde{\mathcal{I}}_h^H \colon \mathrm{T}_{x_k^h}\cMh \to \mathrm{T}_{x_k^H}\cMH, \qquad \widetilde{\mathcal{I}}_H^h \colon \mathrm{T}_{x_k^H}\cMH \to \mathrm{T}_{x_k^h}\cMh.
\]
For embedded submanifolds of the Euclidean space, a standard construction combines linear transfers $I_h^H$ and $I_H^h$ between the ambient spaces with orthogonal projection onto the target tangent spaces:
\[
    \widetilde{\mathcal{I}}_h^H = \operatorname{P}_{x_k^H} \circ I_h^H\big|_{\mathrm{T}_{x_k^h}\cMh}, \qquad \widetilde{\mathcal{I}}_H^h = \operatorname{P}_{x_k^h} \circ I_H^h\big|_{\mathrm{T}_{x_k^H}\cMH},
\]
where $\operatorname{P}_{x_k^h}$ and $\operatorname{P}_{x_k^H}$ denote the orthogonal projectors onto $\mathrm{T}_{x_k^h}\cMh$ and $\mathrm{T}_{x_k^H}\cMH$, respectively. These projections ensure tangency, which the ambient transfers alone need not preserve. This construction is analogous to projection-based vector transport in retraction-based Riemannian optimization; see~\cite[\S 8.1]{AMS:2008}.

When working with the manifold of fixed-rank matrices, tangent transfers can be implemented efficiently by exploiting the factored representation of tangent vectors. For details on the implementation of point and tangent transfers in the case of the manifold of fixed-rank matrices, we refer the readers to~\cite[\S 3.3]{sutti2021riemannian}.

For theoretical purposes, we assume that the tangent transfer operators satisfy the condition
\[
    \widetilde{\mathcal{I}}_h^H = (\widetilde{\mathcal{I}}_H^h)^{\ast},
\]
where the superscript ${}^{\ast}$ denotes the adjoint with respect to the Riemannian inner products on the fine and coarse tangent spaces; see, e.g.,  \cite[(A.2.8)]{Trottenberg:2000} and \cite[(3.5.3)]{Hackbusch:2003}. We also assume that these operators are uniformly bounded:
\begin{equation}\label{eq:boundop}
    \| \widetilde{\mathcal{I}}_h^H \|_{\mathrm{op}} = \| \widetilde{\mathcal{I}}_H^h \|_{\mathrm{op}} \leq C_{\mathcal{I}},    
\end{equation}
with a constant $ C_{\mathcal{I}} > 0 $ independent of the iteration and the level considered. Here, $ \| \cdot \|_{\mathrm{op}} $ denotes the induced operator norm; in particular,
\[ 
    \| \widetilde{\mathcal{I}}_h^H \|_{\mathrm{op}} \coloneqq \sup_{\substack{s^{h} \in \mathrm{T}_{x_k^h}\cMh, \\ s^{h} \neq 0}} \frac{\| \widetilde{\mathcal{I}}_h^H s^{h} \|_{x_k^H}}{\| s^{h} \|_{x_k^h}},
\]
and analogously for $ \widetilde{\mathcal{I}}_H^h $.

\section{The Riemannian Multilevel Adaptive Regularized Gradient method} \label{sec:two-level}

We describe a multilevel extension of the Riemannian adaptive regularization method with first-order models, cf. \Cref{app:RAR1}. The goal is to minimize a smooth function on a high-dimensional “fine” manifold by exploiting a hierarchy of lower-dimensional, computationally cheaper manifolds that preserve the essential structure of the original problem. At each iteration, the algorithm chooses between two kinds of step. The first is a standard regularized gradient step taken directly on the fine manifold (\textit{fine step}). The second is a \textit{coarse step}, which is used when the current gradient is well represented at the coarse level. In that case, the iterate is transferred to the coarse manifold, where a regularized model of the objective, coherent with the fine problem, is approximately minimized by the same method. When the hierarchy contains more than two levels\footnote{The actual multilevel version is described in \Cref{subsec:multilevel_implementation}.}, this minimization is carried out by applying the same method recursively: at each intermediate level the algorithm again decides whether to take a gradient step or to move further down, until the coarsest level is reached, where only gradient steps are available. The resulting step is then brought back to the fine level. Whichever kind of step is used, it is accepted or rejected by comparing the actual decrease with the one predicted by the model; this comparison also drives the update of the regularization parameter, which becomes smaller after successful steps and larger after unsuccessful ones.

We now describe the details of the main steps of the method considering the case where the hierarchy has only two levels.

\subsection{Step definition}
In this subsection, we describe how to construct the fine step $s^{h}_{k}$ and the coarse step $s^{H}_{k\,\ast}$.
We start in \Cref{subsec:fine_step} with the fine step, which closely follows RAR1, cf. \Cref{app:RAR1}. The coarse step, described in \Cref{subsec:coarse_step}, involves the definition of a coarse model and its regularized version.

\subsubsection{Fine step} \label{subsec:fine_step}
At the fine level, the step $s^{h}_{k}$ is defined following the RAR1 method. Thus, $s^{h}_{k}$ is the minimizer of a linear Taylor model with a quadratic regularization.
At the current iterate $x^{h}_{k}$, the linear model around $x^{h}_{k}$ is $T^{h}_{k} \colon \T_{x^{h}_{k}}\cMh\to\R$ defined as
\begin{align}\label{eq:fine_lin_model}
    T^{h}_{k}(s^{h}) \coloneqq& \hat{f}^h_k(0_{x_k^h})+\innerproduct{\nabla \hat{f}^{h}_k(0_{x_k^h})}{s^{h}}_{x^{h}_{k}}
\end{align}
where $\innerproduct{\cdot}{\cdot}_{x^{h}_{k}}$ is the inner product defined on the tangent space $ \T_{x^{h}_{k}}\cMh $. The regularized model $m^{h}_{k}(s^{h})$ is defined as 
\begin{equation}\label{fine reg model}
    m^{h}_{k}(s^{h}) = T^{h}_{k}(s^{h}) + \frac{1}{2} \, \lambda_k \, \| s^{h} \|^{2}_{x^{h}_{k}},
\end{equation}
where $\lambda_k>0$ is the regularization parameter.
The step $s^{h}_{k}$ is defined as a minimizer of the model~\eqref{fine reg model}; in this case, it can be computed in closed form as
\begin{equation}\label{eq:fine_direction}
    s^{h}_{k} =  - \frac{1}{\lambda_k} \, \grad f^{h}(x^{h}_{k}).
\end{equation}

\subsubsection{Coarse step} \label{subsec:coarse_step}
At iteration $k$, let the coarse base point be $x_k^H=\cI_h^H x_k^h\in\cMH$. Let $\mathcal U_k^H$ be a neighborhood of $0_{x_k^H}$ in $\T_{x_k^H}\cMH$ such that $\operatorname{R}_{x_k^H}$ is a diffeomorphism\footnote{Since $\operatorname{D} \! \operatorname{R}_{x_k^H}(0_{x_k^H}) = \operatorname{Id}_{\mathrm{T}_{x_k^H}\cMH}$, the inverse function theorem on manifolds guarantees the existence of an open neighborhood $\mathcal U_k^H$ of $0_{x_k^H}$ in $\mathrm{T}_{x_k^H}\cMH$ such that $\operatorname{R}_{x_k^H}|_{\mathcal U_k^H}$ is a diffeomorphism onto $\mathcal{V}_k^H \coloneqq \operatorname{R}_{x_k^H}(\mathcal U_k^H)$; see~\cite[Thm.~4.16 and Cor.~4.17]{boumal2023introduction}. Consequently, the local inverse $\operatorname{R}_{x_k^H}^{-1}$ is smooth on $\mathcal{V}_k^H$.} onto $ \mathcal{V}_k^H\coloneqq\operatorname{R}_{x_k^H}(\mathcal U_k^H) \subseteq\cMH $. We emphasize that $x^H_k$ denotes the \emph{fixed} base point, while $x^H$ is the \emph{variable} in $\mathcal V_k^H$. Let the inverse of $\operatorname{R}_{x_k^H}$ on $\mathcal V_k^H$ be $\operatorname{R}_{x_k^H}^{-1}$. The coarse model $\psi_k^H \colon \mathcal V_k^H \subseteq\cMH \to \R $ is defined by
\begin{equation}\label{eq:coarse_model_psiHk}
    \psi^{H}_{k}(x^{H}) =  f^{H}(x^{H}) + \innerproduct{v^{H}_{k}}{\operatorname{R}^{-1}_{x^{H}_{k}}(x^{H})}_{x^{H}_{k}},
\end{equation}
where $v_k^H\in\T_{x_k^H}\cMH$ is known in the multilevel literature as the \textit{first-order coherence term} \cite{Nash:2000,Lewis:2005} since it ensures coherence between levels, cf. \Cref{prop:coherence} below, and it is defined as
\begin{equation}\label{eq:definition_vHk}
    v^{H}_{k} = \widetilde{\cI}_{h}^{H}\grad f^{h} (x^{h}_{k}) - \grad f^{H} ( x^{H}_{k} ).
\end{equation}
For every $s^H\in \mathcal{U}_k^H$, the definition of $\psi_k^H$ gives
\begin{align}\label{eq:coarse_model_psiHk_new}
\psi_k^H(\operatorname{R}_{x_k^H}(s^H)) &= f^H(\operatorname{R}_{x_k^H}(s^H))+\innerproduct{v_k^H}{\operatorname{R}_{x_k^H}^{-1}(\operatorname{R}_{x_k^H}(s^H))}_{x_k^H}\nonumber\\
    &=f^H(\operatorname{R}_{x_k^H}(s^H))+\innerproduct{v_k^H}{s^H}_{x_k^H},
\end{align}
i.e., introducing the pullback of $\psi_k^H$,
\[
    \hat{\psi}^H_k(s^H) = \hat{f}^H(s^H)+\innerproduct{v_k^H}{s^H}_{x_k^H}.
\]
In particular, since $\operatorname{R}_{x_k^H}(0_{x_k^H})=x_k^H$ and $\operatorname{R}_{x_k^H}^{-1}(x_k^H)=0_{x_k^H}$, the definition of $\psi_k^H$ yields 
\begin{equation}\label{eq:uguale}
    \psi_k^H(x_k^H) = f^H(x_k^H) \qquad \text{ and } \qquad \hat{\psi}_k^H(0_{x_k^H})= \hat{f}^H(0_{x_k^H}).
\end{equation}
For $x^H \in \mathcal{V}_k^H$ and $\lambda_k>0$, we define the \textit{regularized coarse model} as the sum of the coarse model $\psi^{H}_{k}$, defined in~\eqref{eq:coarse_model_psiHk}, and a quadratic regularization term, i.e.,
\begin{equation}\label{eq:reg_coarse_model}
    m^{H}_{k}(x^{H}) = \psi^{H}_{k}(x^{H}) + \frac{1}{2} \lambda_k \| \operatorname{R}^{-1}_{x^{H}_{k}}(x^{H}) \|^{2}_{x^{H}_{k}}.
\end{equation}
Let $x^{H}_{k,\ast} \in \mathcal{V}_k^H$ be an approximate minimizer of $m_k^H$ over $\mathcal{V}_k^H$, and define the tangent vector
\[
    s_{k,\ast}^H \coloneqq \operatorname{R}_{x_k^H}^{-1}(x^{H}_{k,\ast})\in\mathcal U_k^H\subseteq\T_{x_k^H}\cMH.
\]
Then, the fine-level step $s_k^h\in\T_{x_k^h}\cMh$ is obtained by prolonging  $s_{k,\ast}^H$:
\begin{equation}\label{eq:coarse_direction}
    s_k^h=\widetilde{\cI}_{H}^{h}s_{k,\ast}^H.
\end{equation}

On $\mathcal{V}_k^H \subseteq \cMH$, starting from $x^H_{(0)}\coloneqq x_k^H$, we iteratively compute $x^{H}_{k,\ast}$ by approximately solving
\begin{equation}\label{eq:coarse_min_problem2}
    \underset{x^H \in \mathcal{V}_k^H}{\min}\;m_k^H(x^H)
\end{equation}
using a Riemannian solver that is globally convergent to a first-order stationary point, such as RAR1 (see \Cref{app:RAR1}) or higher-order methods, such as the Riemannian ARC \cite{agarwal2021adaptive} or the Riemannian trust-region \cite{absil2007trust,boumal2023introduction} algorithms. 

\begin{figure}[htb]
    \centering
    \includegraphics[width=\linewidth]{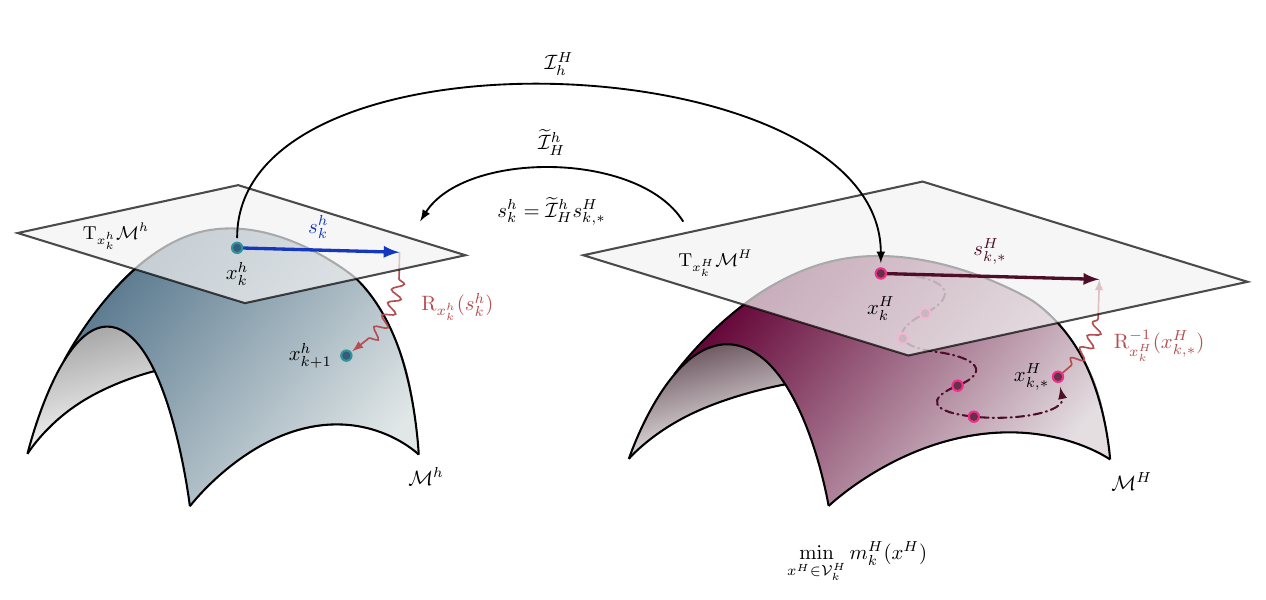}
    \caption{Schematic illustration of a coarse correction in \algname{}.}
    \label{fig:rmar1_scheme}
\end{figure}

\Cref{fig:rmar1_scheme} schematically illustrates a coarse correction step. At iteration $k$, the fine iterate $x_k^{h}$ is restricted to $x_k^{H}=\mathcal{I}_{h}^{H}(x_k^{h})$. The dash-dot curve on $\cMH$ represents the coarse optimization trajectory toward $x_{k,\ast}^{H}$, with intermediate points denoting coarse optimizer iterates. This procedure corresponds to approximately minimizing the model $m_k^{H}$ over the neighborhood $\mathcal{V}_k^{H}$. The tangent vector update $s_{k,\ast}^{H}=\mathrm{R}_{x_k^{H}}^{-1}(x_{k,\ast}^{H})$ is prolonged to the fine tangent space as $s_k^{h}=\widetilde{\mathcal{I}}_{H}^{h}s_{k,\ast}^{H}$. For an accepted correction with unit step length, $x_{k+1}^{h}=\mathrm{R}_{x_k^{h}}(s_k^{h})$. The surfaces and trajectory are schematic; the translucent tangent planes reveal the underlying manifold and optimizer iterates.

As we will see in \Cref{sec:theory}, for the convergence of the multilevel method, an approximate solution of \eqref{eq:coarse_min_problem2} is sufficient, so the coarse iterations $\curly{x^{H}_{(i)}}_{i\in\N}$ are stopped  with $x^{H}_{k,\ast} = x^{H}_{(i+1)}$ satisfying the following conditions:
\begin{equation}\label{eq:coarse_stop_crit_twolevel}
    m^{H}_{k}(x^{H}_{k,\ast})<m^{H}_{k}(x^{H}_{(0)}), \quad \text{and} \quad \norm{\grad m^{H}_{k}(x^{H}_{k,\ast})}_{x_{\ast}^H}\leq \theta^{H}\norm{\operatorname{R}_{x^{H}_{k}}^{-1}(x^{H}_{k,\ast})}_{x^{H}_{k}},
\end{equation}
for some $\theta^{H}>0$, and we then set $ s^{H}_{k,\ast} \coloneqq \operatorname{R}_{x^{H}_{k}}^{-1}(x^{H}_{k,\ast})$. These conditions are the Riemannian equivalent of those used in the Euclidean context in \cite{Calandra:2021,MariniPorcelliRiccietti2026}, and can always be satisfied; see, e.g., \cite[Sec. 3.3]{cartis2011adaptive}.

The coarse step is used only for iterations such that:
\begin{equation}\label{eq:go_down_cond}
    \norm{\widetilde{\cI}_{h}^{H}\grad f^{h}(x^{h}_{k})}_{x^{H}_{k}} \geq \kappa_{H}\norm{\grad f^{h}(x^{h}_{k})}_{x^{h}_{k}},
\end{equation}
where
$ \kappa_{H} \in ( 0, \min\lbrace 1, \| \widetilde{\cI}_{h}^{H}\| \rbrace ) $. This last condition requires the restricted fine-level gradient to retain a prescribed fraction of the norm of the fine gradient and is crucial to ensure that $ s^{h} $ is a descent direction for $f^{h}$ at $ x^{h}_{k} $, if $ s^{H} $ is a descent direction for $ \psi^{H}_{k} $ at $ x^{H}_{k} $.

The next proposition ensures that the behaviors of $f^h$ and $\psi_k^H$ are coherent up to first order in the neighborhood of $x^{H}_{k}$. This allows us to be sure that $f^h$ is decreased, while working with $\psi_k^H$. Indeed, if $s^H$ is a descent direction for $\psi_k^H$, and we define $ s^{h} = \widetilde{\cI}_{H}^{h} s^{H}$, we have that $s^h$ is also a descent direction for the function $f^h$.

\begin{proposition}[First-order coherence]\label{prop:coherence}
For the coarse model $\psi_k^H$ defined in (\ref{eq:coarse_model_psiHk}), we have
$$
 \innerproduct{\grad\psi_k^H(x_k^H)}{s^H}_{x_k^H} =  \innerproduct{\grad f^{h} (x^{h}_{k}) }{s^{h}}_{x^{h}_{k}}
$$
for any $x_k^H \in \mathcal{V}_k^H \coloneqq \operatorname{R}_{x_k^H}(\mathcal U_k^H) \subseteq\cMH $, $ x_k^{h} = \cI_{H}^{h} x_k^{H}$, $s^H\in\T_{x^H}\cM^H$, and $ s^{h} = \widetilde{\cI}_{H}^{h} s^{H} \in\T_{x_k^h}\cM^H$.
\end{proposition}
\begin{proof}
Let $x^H \in \mathcal{V}_k^H$ and $s^H\in\T_{x^H}\cM^H$. We have
\begin{align*}
    \innerproduct{\grad\psi_k^H(x^H)}{s^H}_{x^H} &= \innerproduct{\grad f^H(x^H) + \Big(\operatorname{D} \! \operatorname{R}_{x_k^H}^{-1}(x^H)\Big)^\ast\left[v_k^H\right]}{s^H}_{x^H}\\
    &=\innerproduct{\grad f^H(x^H) + \Big(\operatorname{D} \! \operatorname{R}_{x_k^H}^{-1}(x^H)\Big)^{\ast} \Big[ \widetilde{\cI}_h^H\grad f^h(x_k^h)-\grad f^H(x_k^H) \Big]}{s^H}_{x^H}.
\end{align*}

If $ x^{H} = x^{H}_{k} $ and $s^H\in\T_{x_{k}^H}\cM^H$, then $ \operatorname{D} (\operatorname{R}_{x^{H}_{k}}^{-1})(x^{H}_{k}) = \operatorname{Id}_{\T_{x^{H}_{k}}\cM^H} $, and, using~\eqref{eq:definition_vHk},
\begin{align*}
    \innerproduct{\grad \psi^{H}_{k}(x^{H}_{k})}{s^{H}}_{x^{H}_{k}} &= \innerproduct{\grad f^{H}(x^{H}_{k}) + v^{H}_{k}}{s^{H}}_{x^{H}_{k}} \\
    &= \innerproduct{\grad f^{H}(x^{H}_{k}) + \widetilde{\cI}_{h}^{H}\grad f^{h} (x^{h}_{k}) - \grad f^{H} ( x^{H}_{k} )}{s^{H}}_{x^{H}_{k}} \\
    &= \innerproduct{\widetilde{\cI}_{h}^{H}\grad f^{h} (x^{h}_{k}) }{s^{H}}_{x^{H}_{k}} \\
    &= \innerproduct{\grad f^{h} (x^{h}_{k}) }{s^{h}}_{x^{h}_{k}},
\end{align*} 
where in the last step we used $ \widetilde{\cI}_{h}^{H} = \left({\widetilde{\cI}_{H}^{h}}\right)^{\ast} $ and $ s^{h} = \widetilde{\cI}_{H}^{h} s^{H} $, where the adjoint is taken with respect to the inner products on $\mathrm{T}_{x^{h}_{k}}\cMh$ and $\mathrm{T}_{x^{H}_{k}}\cM^H$.

\end{proof}

In \Cref{app:suppl}, we specialize these calculations to prove the first-order coherence property for the manifold of fixed-rank matrices, when adopting the orthographic retraction \cite[Prop. 24]{Absil:2012}. For the orthographic retraction, see later \Cref{subsec:retraction_choice}.

\subsection{Step acceptance}
Once we obtain the tentative step $ s^{h}_{k} $ either from \eqref{eq:fine_direction} or from
\eqref{eq:coarse_direction}, we evaluate the reduction of the objective function at the trial point $ \operatorname{R}_{x^{h}_{k}}(s^{h}_{k}) $ to decide whether to accept it or not. This is done via the ratio 
\begin{equation}\label{eq:fine_rho}
    \rho^{h}_{k} = \frac{f^{h}(x^{h}_{k})-f^{h}(\operatorname{R}_{x^{h}_{k}}(s^{h}_{k}))}{\text{pred. reduction}}
    =   \frac{\hat{f}^{h}_{k}(0_{x^{h}_{k}}) - \hat{f}^{h}_{k}(s^{h}_{k})}{\text{pred. reduction}}.
\end{equation}
The numerator represents the achieved actual reduction in the cost function. The denominator represents the predicted reduction, calculated as
\begin{equation}\label{eq:fine_pred_reduction}
    \text{pred. reduction} =
    \begin{cases}
        T_k^h(0_{x_k^h})-T_k^h(s_k^h), & \text{if a fine step is used},\\
        \psi_k^H(x_k^H)-\psi_k^H(x^{H}_{k,\ast}), & \text{if a coarse step is used}.
    \end{cases}
\end{equation}
If the ratio $\rho^{h}_{k}$ is larger than or equal to a given threshold $\eta_{1}\in (0,1)$, we accept the step $ s^{h}_{k} $ and the iteration is declared successful and the iterate is updated as $ x^{h}_{k+1} = \operatorname{R}_{x^{h}_{k}}(s^{h}_{k}) $; otherwise, the iteration is declared unsuccessful and the step is rejected, thus we keep the iterate unchanged, i.e., $x^{h}_{k+1} = x^{h}_{k}$.

Next, we have to update the regularization parameter $\lambda_k$ to $\lambda_{k+1}$ according to the following criterion
\begin{equation}\label{eq:fine_reg_update}
\lambda_{k+1} =
    \begin{cases}
        \max \curly{\lambda_{\min}, \gamma_{2} \lambda_k} & \text { if } \rho^{h}_{k} \geq \eta_2 \quad \text{very successful,}\\
        \max \curly{\lambda_{\min}, \gamma_{1} \lambda_k} & \text { if } \rho^{h}_{k} \in\left[\eta_1, \eta_2\right) \quad \text{successful,} \\
        \gamma_{3} \lambda_k & \text { if }\rho^{h}_{k}<\eta_1 \quad  \text{unsuccessful,}
    \end{cases}
\end{equation}
where $\eta_2\in(\eta_1,1)$, $0<\gamma_2<\gamma_1<1<\gamma_3$, and $\lambda_{\min}>0$. On a successful iteration, the regularization parameter is decreased by $\gamma_1$ or remains at its lower bound. A very successful iteration uses the smaller reduction factor $\gamma_2$. On an unsuccessful iteration, the parameter is increased by the factor $\gamma_3$.

\Cref{alg:riemannian_two-level_ar1} summarizes this description by providing a pseudocode for the two-level \algname{} method.

\begin{center}
\begin{algorithm}[htbp]
    \caption{Two-level \algname{}.}
    \label{alg:riemannian_two-level_ar1}
    \SetAlgoLined
    \DontPrintSemicolon
    \KwIn{Initial point $x^{h}_0 \in \cMh$ and initial regularization parameter $\lambda_0 \geq \lambda_{\min}$.}
    \KwOut{An iterate $x^{h}_\ast$ satisfying $\lVert \grad f^{h}(x^{h}_\ast) \rVert_{x^{h}_\ast} \leq \varepsilon^{h}$.}
    \textbf{Parameters:} $\theta^{H} > 0$, $\varepsilon^{h} > 0$, $\lambda_{\min} > 0$, $0 < \eta_1 \leq \eta_2 < 1$, and $0 < \gamma_2 < \gamma_1 < 1 < \gamma_3$.\;
    Set $k \gets 0$ and compute $\grad f^{h}(x^{h}_0)$.\;
    \While{$\lVert \grad f^{h}(x^{h}_k) \rVert_{x^{h}_k} > \varepsilon^{h}$}{
        Set $x^{H}_k \coloneqq \cI_h^{H}x^{h}_k$ and compute $\widetilde{g}^{H}_k \coloneqq \widetilde{\cI}_h^{H} \grad f^{h}(x^{h}_k)$.\;
        \eIf{\eqref{eq:go_down_cond} holds}{
            \textbf{Coarse step.} \;
            Compute $\grad f^{H}(x^{H}_k)$ and construct $\psi^{H}_k$ and $m^{H}_k$ by~\eqref{eq:coarse_model_psiHk} and~\eqref{eq:reg_coarse_model}.\;
            Starting from $x^{H}_k$, approximately minimize $m^{H}_k$ on $\cMH$ to obtain $x^{H}_{k,\ast}$ satisfying~\eqref{eq:coarse_stop_crit_twolevel}.\;
            Set $s^{H}_{k,\ast} \coloneqq \operatorname{R}_{x^{H}_k}^{-1}(x^{H}_{k,\ast})$ and $s^{h}_k = \widetilde{\cI}_H^{h}s^{H}_{k,\ast}$.\;
            Set $\mathrm{pred}^{h}_k = \psi^{H}_k(x^{H}_k) - \psi^{H}_k(x^{H}_{k,\ast})$.\;
        }{
            \textbf{Taylor step.} \;
            Set $s^{h}_k = - \frac{1}{\lambda_k} \grad f^{h}(x^{h}_{k})$.\;
            Set $\mathrm{pred}^{h}_k = \lVert \grad f^{h}(x^{h}_k) \rVert_{x^{h}_k}^{2}/\lambda_k$.\;
        }
        Set $\tilde{x}^{h}_{k} = \operatorname{R}_{x^{h}_k}(s^{h}_k)$ and compute $\rho^{h}_k$ from $\mathrm{pred}^{h}_k$ and the fine objective decrease by~\eqref{eq:fine_rho} and~\eqref{eq:fine_pred_reduction}.\;
        \eIf{$\rho^{h}_k \geq \eta_1$}{
            Set $x^{h}_{k + 1} = \tilde{x}^{h}_{k}$.\;
        }{
            Set $x^{h}_{k + 1} = x^{h}_k$.\;
        }
        Set $\lambda_{k + 1}$ from $\lambda_k$ and $\rho^{h}_k$ by~\eqref{eq:fine_reg_update}.\;
        Set $k \gets k + 1$ and compute $\grad f^{h}(x^{h}_k)$.\;
    }
    Set $x^{h}_\ast \gets x^{h}_k$.\;
    \Return{$x^{h}_\ast$}.\;
\end{algorithm}
\end{center}

\section{Theoretical analysis of \algname{}} \label{sec:theory}

In this section, we provide a theoretical analysis of the proposed two-level \algname{} method\footnote{Note that the analysis conducted in \Cref{sec:theory} extends to the multilevel case due to its recursive nature.} described in \Cref{alg:riemannian_two-level_ar1}.
In \Cref{subsec:global_convg}, we prove global convergence to  critical points, and in \Cref{subsec:worst_case_complexity} we provide a worst-case complexity bound for reaching $\varepsilon$-approximate first-order minimizers. For the analysis, we need the following regularity assumptions.

\begin{assumption}\label{ass:smooth}
The cost functions $f^h$ and $f^H$ are continuously differentiable.
\end{assumption}

\begin{assumption}\label{ass:low}
The cost functions $f^h$ and $f^H$ are bounded below.
\end{assumption}

\begin{assumption}\label[assumption]{ass:pullback_grad_lipschitz}
There exist constants $L_h>0$ and $L_H>0$ such that, for all iterates $x^h_k$, $x^H_k$ and steps $s^h_k$, $s^H_{k,\ast}$ generated by \Cref{alg:riemannian_two-level_ar1}, the gradients of the pullbacks $\hat f^h_k$ and $\hat f^H_k$ satisfy
\begin{align*}
  \| \nabla \hat f^h_k(t s^h_k) - \nabla \hat f^h_k(0_{x^h_k}) \|_{x^h_k} &\le L_h \| t s^h_k \|_{x^h_k}, \quad \forall t\in[0,1],\\
  \| \nabla \hat f^H_k(t s^H_{k,\ast}) - \nabla \hat f^H_k(0_{x^H_k}) \|_{x^H_k} &\le L_H \| t s^H_{k,\ast} \|_{x^H_k}, \quad \forall t\in[0,1].
\end{align*}
\end{assumption}

\begin{assumption}\label[assumption]{ass:jacbound} 
For all $x_k^H$ and $ s^{H}_{k,\ast} $ generated by \Cref{alg:riemannian_two-level_ar1}, the differential of the retraction is uniformly bounded, i.e., there exists a constant $C_{\mathrm{R}}$ such that 
\begin{equation}\label{eq:CR}
     \|\operatorname{D} \! \operatorname R_{x_k^H}(s^{H}_{k,\ast})\|_{\mathrm{op}}
\leq C_{\mathrm{R}}. 
\end{equation}
  \end{assumption}

From~\Cref{ass:pullback_grad_lipschitz} and by applying the Fundamental Theorem of Calculus along the segment $\{t s^h_k \colon t\in[0,1]\}$ in $\T_{x^h_k}\cMh$ (and analogously along $\{t s^H_{k,\ast} \colon t\in[0,1]\}$ in $\T_{x^H_k}\cMH$), we have that, for all iterates $x^h_k$, $x^H_k$ and steps $s^h_k$, $s^H_{k,\ast}$ generated by \Cref{alg:riemannian_two-level_ar1}, the following descent property of the pullbacks holds:
\begin{equation}\label{eq:liph}
    \hat{f}_{k}^h(s^h_k) - \left[ \hat{f}_{k}^h(0_{x_k^h}) + \innerproduct{s^h_k}{\nabla \hat{f}_{k}^h(0_{x^h_k})}_{x^h_k} \right] \leq \frac{L_h}{2} \|s^h_k\|_{x^h_k}^2,
\end{equation}
\begin{equation}\label{eq:lipH}
    \hat{f}_{k}^H(s^H_{k,\ast}) - \left[ \hat{f}_{k}^H(0_{x_k^H}) + \innerproduct{s^H_{k,\ast}}{\nabla \hat{f}^H_k(0_{x_k^H})}_{x^H_k} \right] \leq \frac{L_H}{2} \|s^H_{k,\ast}\|_{x^H_k}^2.
\end{equation}

\subsection{Global convergence} \label{subsec:global_convg}
In this subsection, we establish the global convergence properties of the \algname{} method. Our analysis proceeds in three sequential steps. First, we establish an upper bound for the regularization parameter $\lambda_k$ (\Cref{prop:upper_bound_lambdak}). Second, we derive a fundamental relationship between the norm of the computed step and the norm of the Riemannian gradient (\Cref{prop:steplength_and_gradnorm}). Finally, we combine these two structural ingredients to prove that the sequence of gradient norms converges asymptotically to zero, as stated in our main theoretical result \Cref{thm: Global CV}.

Throughout this subsection, whenever an inverse retraction is used, we assume that the retraction is a diffeomorphism from a neighborhood of the origin containing the relevant tangent steps onto a neighborhood of the base point containing all points at which the inverse is evaluated; see our earlier discussion in \Cref{subsec:coarse_step}.

\begin{proposition}[Upper bound for the regularization parameter $\lambda_k$] \label{prop:upper_bound_lambdak} Let Assumptions \ref{ass:smooth}, \ref{ass:low}, and \ref{ass:pullback_grad_lipschitz} hold. 
If the regularization parameter $\lambda_k$ in \Cref{alg:riemannian_two-level_ar1} is updated according to \eqref{eq:fine_reg_update},  it holds
\begin{equation} \label{eq:lambda_max}
  \lambda_k \le   \lambda_{\max} \coloneqq \max\Big\{\lambda_0,\gamma_3\frac{M}{1-\eta_{1}}\Big\},
\end{equation}
where 
\begin{equation}\label{eq:defM}
    M = \begin{cases} 
            2L_{h} & \text{(fine step)}, \\
            L_H + C_\mathcal{I}^2 L_h & \text{(coarse step)},
        \end{cases}
\end{equation}
$C_\mathcal{I}$ is defined in \eqref{eq:boundop} and $\eta_1,\gamma_3$ and $\lambda_0$ are given in \Cref{alg:riemannian_two-level_ar1}. 
\end{proposition}

\begin{proof}
At any given iteration $k$, the step is obtained either by explicitly minimizing the regularized fine-level Taylor model \eqref{fine reg model}, or by computing an approximate minimizer of the regularized coarse-level model \eqref{eq:reg_coarse_model} such that the stopping conditions \eqref{eq:coarse_stop_crit_twolevel} are satisfied. Consequently, the respective model reductions satisfy the following bounds:
\begin{align}
T_{k}^{h}(0_{x_k^h}) - T_{k}^{h}(s^h_k) &\ge \frac{\lambda_k}{2} \|s^h_k\|^2_{x^h_k}, \label{eq:15a} \\
\psi_{k}^{H}(x^H_k) - \psi_{k}^{H}(x^H_{k,\ast}) & \ge \frac{\lambda_k}{2} \|\operatorname{R}_{x^{H}_{k}}^{-1}(x_{k,\ast}^H)\|^2_{x^H_k}. \label{eq:pred_red_coarse}
\end{align}

We now bound the deviation of the step-acceptance ratio from unity, namely, $|1-\rho_k|$, where $\rho_k$ is defined via the ratio of the actual and predicted reductions in \eqref{eq:fine_rho} and \eqref{eq:fine_pred_reduction}. 

If the fine-level Taylor model is selected at iteration $k$,  the descent property in \eqref{eq:liph}--\eqref{eq:lipH} together with the model reduction bound \eqref{eq:15a} yield:
\begin{equation*}
|1-\rho_k| = \left|1 - \frac{f^{h}(x^h_k) - f^{h}(\operatorname{R}_{x^h_k}(s^h_k))}{T_{k}^{h}(0_{x_k^h}) - T_{k}^{h}(s^h_k)} \right| = \left|\frac{f^{h}(\operatorname{R}_{x^h_k}(s^h_k)) - T_{k}^{h}(s^h_k)}{T_{k}^{h}(0_{x_k^h}) - T_{k}^{h}(s^h_k)} \right| \le \frac{L_{h}}{\lambda_k}.
\end{equation*}

Alternatively, if the coarse-level model is selected, the ratio satisfies:
\begin{equation}\label{eq:coefficient}
    |1-\rho_{k}| = \left|\frac{\psi^{H}_k(x^H_k) - \psi^{H}_k(x^H_{\ast}) - (f^{h}(x^h_k) - f^{h}(\operatorname{R}_{x^h_k}(s^h_k)))}{\psi^{H}_k(x^H_k) - \psi^{H}_k(x^H_{k,\ast})} \right|.
\end{equation}
From the definition of the coarse model $\psi_k^H$ from \eqref{eq:coarse_model_psiHk} and \eqref{eq:uguale}, it holds
\begin{equation}\label{eq:diff}
       \psi_k^H(x^H_k) - \psi_k^H(x^H_{k,\ast}) = f^H(x_k^H)-f^H(x^{H}_{k,\ast})- \innerproduct{v_k^H}{s^{H}_{k,\ast}}. 
\end{equation}
Applying the Fundamental Theorem of Calculus to $f^H$ in $[0, s^{H}_{k,\ast}]$ through the pullback definition yields
\begin{align*}
    f^H(x^H_{k,\ast}) = f^H(\operatorname{R}_{x^H_k}(s^H_{k,\ast})) = f^H(x^H_k) &+ \innerproduct{\grad f^H(x^H_k)}{s^H_{k,\ast}}_{x^H_k} \\ &+ \int_0^1 \innerproduct{\nabla \hat{f}_{k}^H(t s^H_{k,\ast}) - \grad f^H(x^H_k)}{s^H_{k,\ast}}_{x^H_k} \, \mathrm{d}t.
\end{align*}
Substituting this into \eqref{eq:diff}, we can express the predicted reduction as:
\begin{align}\label{eq:denominator}
       \psi_k^H(x^H_k) - \psi_k^H(x^H_{k,\ast}) &= - \innerproduct{\grad f^H(x^H_k)}{s^H_{k,\ast}}_{x^H_k} - \innerproduct{s^H_{k,\ast}}{\itangent{h}{H}\grad f^{h}(x^h_k) - \grad f^{H}(x^H_k)}_{x^H_k} \notag \\
    &\quad - \int_0^1 \innerproduct{\nabla \hat{f}_{k}^H(t s^H_{k,\ast}) - \grad f^H(x^H_k)}{s^H_{k,\ast}}_{x^H_k} \, \mathrm{d}t \notag \\
    &= - \innerproduct{s^H_{k,\ast}}{\itangent{h}{H}\grad f^{h}(x^h_k)}_{x^H_k} - \int_0^1 \innerproduct{\nabla \hat{f}_{k}^H(t s^H_{k,\ast}) - \grad f^H(x^H_k)}{s^H_{k,\ast}}_{x^H_k} \, \mathrm{d}t.
\end{align}
Similarly, applying the integral Taylor expansion to the fine-level objective function gives:
\begin{equation}\label{eq:numerator}
    f^h(\operatorname{R}_{x^h_k}(s^h_k)) - f^h(x^h_k) = \innerproduct{\grad f^h(x^h_k)}{s^h_k}_{x^h_k} + \int_0^1 \innerproduct{\nabla \hat{f}_{k}^h(t s^h_k) - \grad f^h(x^h_k)}{s^h_k}_{x^h_k} \, \mathrm{d}t.
\end{equation}
By the adjoint relation between the tangent transfer operators, namely $(\itangent{H}{h})^\ast = \itangent{h}{H}$, the first-order linear terms match exactly:
\[ 
    \innerproduct{\grad f^{h}(x^h_k)}{s^h_k}_{x^h_k} = \innerproduct{\grad f^{h}(x^h_k)}{\itangent{H}{h}s^H_{k,\ast}}_{x^h_k} = \innerproduct{\itangent{h}{H}\grad f^{h}(x^h_k)}{s^H_{k,\ast}}_{x^H_k}.
\]
and summing \eqref{eq:numerator} and \eqref{eq:denominator} cancels them. Applying the triangle inequality along with the Lipschitz continuity of the pullback gradients (\Cref{ass:pullback_grad_lipschitz}), the absolute value of the numerator in \eqref{eq:coefficient} is bounded by:
\begin{align*}
& \left| \psi_k^H(x_k^H) - \psi_k^H(x^H_{k,\ast}) - \left( f^h(x^h_k) - f^h(\operatorname{R}_{x^h_k}(s^h_k)) \right) \right| \\
&\quad \le \int_0^1 \left| \innerproduct{\nabla \hat{f}_{k}^H(t s^H_{k,\ast}) - \nabla \hat{f}_{k}^H(0)}{s^H_{k,\ast}}_{x^H_k} \right| \mathrm{d}t + \int_0^1 \left| \innerproduct{\nabla \hat{f}_{k}^h(t s^h_k) - \nabla \hat{f}_{k}^h(0)}{s^h_k}_{x^h_k} \right| \mathrm{d}t \\
&\quad \le \frac{L_H}{2}\|s^H_{k,\ast}\|_{x^H_k}^{2} + \frac{L_h}{2}\|s^h_k\|_{x^h_k}^{2}.
\end{align*}

Using the boundedness of the prolongation operator (\ref{eq:boundop}) it holds  $\|s^h_k\|_{x^h_k}=\|\itangent{H}{h}s^{H}_{k,\ast}\|_{x^h_k} \le C_\mathcal{I} \|s^H_{k,\ast}\|_{x^H_k}$, and the bound on the numerator simplifies to
\[
    \frac{L_H + C_\mathcal{I}^2 L_h}{2} \, \|s^H_{k,\ast}\|_{x^H_k}^{2}.
\]
Combining this result with the predicted reduction bound \eqref{eq:pred_red_coarse}, we obtain:
\begin{equation*}
    |1-\rho_{k}| \le \frac{L_H + C_\mathcal{I}^2 L_h}{\lambda_k}.
\end{equation*}

Consequently, for both fine and coarse step selections, there exists a strictly positive constant $M$, defined as in \eqref{eq:defM}, such that
\begin{equation} \label{eq:bound_on_rho_deviation}
    |1-\rho_{k}| \le \frac{M}{\lambda_k}.
\end{equation}
    
From \eqref{eq:bound_on_rho_deviation} and the update mechanism for the regularization parameter, we deduce that $\lambda_k$ must be bounded from above. If an iteration is unsuccessful, $\lambda_k$ is scaled up by a factor $\gamma_3 > 1$. As $\lambda_k$ increases, the right-hand side of \eqref{eq:bound_on_rho_deviation} decreases monotonically until it drops below $1-\eta_{1}$. At this threshold, it is guaranteed that $\rho_{k} \ge \eta_{1}$, forcing a successful step where the regularization parameter is reduced. Hence, $\lambda_k$ is bounded globally 
by $\lambda_{\max} \coloneqq \max\{\lambda_0,\gamma_3\frac{M}{1-\eta_{1}}\}$.
\end{proof}

We now show that the step length cannot become arbitrarily small relative to the norm of the objective function's gradient. 

\begin{proposition}[Relating the step length to the norm of the gradient] \label{prop:steplength_and_gradnorm}
Let \Cref{ass:smooth,ass:low,ass:pullback_grad_lipschitz,ass:jacbound} hold. 
For the steps computed by \Cref{alg:riemannian_two-level_ar1}, it holds:
\begin{align}
    \|s^h_k\|_{x^h_k} &= \frac{1}{\lambda_k} \|\grad f^h(x^h_k)\|_{x^h_k},\label{eq:19}\\
    \|s^{H}_{k,\ast}\|_{x^H_k} & \ge \frac{\kappa_{H}}{c_1} \,  \|\operatorname{grad} f^h(x^h_k)\|_{x^h_k},\label{eq:bound_s_ast_sq_2}
\end{align}
with 
\begin{equation}\label{eq:def_c1}
   c_{1} \coloneqq 2 \, C_{\mathcal{I}}^2L_h  + L_H + \lambda_{\max} + C_{\mathrm{R}} \, \theta^H,
\end{equation}
where $C_\mathcal{I}, \lambda_{\max}$, $\kappa_H$ and $C_{\mathrm{R}}$ are defined in \eqref{eq:boundop}, \eqref{eq:lambda_max}, \eqref{eq:go_down_cond}, and \eqref{eq:CR} respectively, while $\theta^H$ is defined in \Cref{alg:riemannian_two-level_ar1}.   
\end{proposition}

\begin{proof}
If a fine-level Taylor step is chosen, the relation simply follows from the definition of the step. 

When the coarse-level model is chosen instead, we consider the restricted fine-level pullback gradient $ \widetilde{\mathcal I}_h^H\nabla\hat{f}_{k}^h(s_k^h) $. Let $\hat{m}_k^H$ denote the pullback of $m_k^H$ at $x_k^H$, defined by $ \hat{m}_k^H(s^H) = m_k^H \bigl(\operatorname R_{x_k^H}(s^H)\bigr), \; s^H\in \T_{x_k^H}\mathcal M^H$. Adding and subtracting $\nabla\hat{\psi}_k^H(s^{H}_{k,\ast})-\nabla\hat{m}_k^H(s^{H}_{k,\ast})$, we obtain
\[
\begin{aligned}
    \widetilde{\mathcal I}_h^H\nabla\hat{f}_{k}^h(s_k^h) &= \widetilde{\mathcal I}_h^H\nabla\hat{f}_{k}^h(s_k^h) +\bigl[ \nabla\hat{\psi}_k^H(s^{H}_{k,\ast}) -\nabla\hat{m}_k^H(s^{H}_{k,\ast})\bigr] -\bigl[ \nabla\hat{\psi}_k^H(s^{H}_{k,\ast}) -\nabla\hat{m}_k^H(s^{H}_{k,\ast}) \bigr]\\
    &= \bigl[ \widetilde{\mathcal I}_h^H\nabla\hat{f}_{k}^h(s_k^h) -\nabla\hat{\psi}_k^H(s^{H}_{k,\ast}) \bigr] +\bigl[ \nabla\hat{\psi}_k^H(s^{H}_{k,\ast}) -\nabla\hat{m}_k^H(s^{H}_{k,\ast}) \bigr] +\nabla\hat{m}_k^H(s^{H}_{k,\ast}).
\end{aligned}
\]

Taking norms on both sides and applying the triangle inequality gives
\begin{equation}\label{eq:triangle_ineq}
    \bigl\| \widetilde{\mathcal I}_h^H\nabla\hat{f}_{k}^h(s_k^h) \bigr\|_{x^H_k} \leq \bigl\| \widetilde{\mathcal I}_h^H\nabla\hat{f}_{k}^h(s_k^h) -\nabla\hat{\psi}_k^H(s^{H}_{k,\ast}) \bigr\|_{x^H_k} +\bigl\| \nabla\hat{\psi}_k^H(s^{H}_{k,\ast}) -\nabla\hat{m}_k^H(s^{H}_{k,\ast}) \bigr\|_{x^H_k} +\bigl\|\nabla\hat{m}_k^H(s^{H}_{k,\ast})\bigr\|_{x^H_k}.
\end{equation}

Next, we bound the terms on the right-hand side individually. The second term is easy to bound. Indeed, 
from the definition of the regularized model,
\[
    \nabla\hat{m}_k^H(s^H) = \nabla\hat{\psi}_k^H(s^H)+\lambda_k s^H,
\]
which implies
\[
    \bigl\| \nabla\hat{\psi}_k^H(s^{H}_{k,\ast}) -\nabla\hat{m}_k^H(s^{H}_{k,\ast}) \bigr\|_{x^H_k} = \lambda_k\|s^{H}_{k,\ast}\|_{x^H_k}\leq \lambda_{\max}\|s^{H}_{k,\ast}\|_{x^H_k},
\]
where $\lambda_{\max}$ is given in \eqref{eq:lambda_max}.

The third term in \eqref{eq:triangle_ineq} is bounded via the algorithm's internal coarse stopping criterion \eqref{eq:coarse_stop_crit_twolevel}. The chain rule yields
\[
    \nabla\hat{m}_k^H(s^{H}_{k,\ast}) = \Bigl( \operatorname{D} \! \operatorname R_{x_k^H}(s^H_{k,\ast}) \Bigr)^{\ast} \Bigl[
\grad m_k^H(x^{H}_{\ast}) \Bigr],
\]
where the adjoint is taken with respect to the Riemannian metrics.
Consequently, using the stopping criterion
\eqref{eq:coarse_stop_crit_twolevel}, we obtain
\begin{align*}
    \|\nabla\hat{m}_k^H(s^{H}_{k,\ast})\|_{x_k^H}
    &\leq
    \|\operatorname{D} \! \operatorname R_{x_k^H}(s^{H}_{k,\ast})\|_{\mathrm{op}} \, 
    \|\grad m_k^H(x^{H}_{\ast})\|_{x^{H}_{\ast}}\\
    &\leq
    \theta^H
    \|\operatorname{D} \! \operatorname R_{x_k^H}(s^{H}_{k,\ast})\|_{\mathrm{op}} \, 
    \|s^{H}_{k,\ast}\|_{x_k^H}.
\end{align*}
Using Assumption \ref{ass:jacbound}, we can  conclude that
\begin{equation}\label{eq:bound_pull_m}
    \| \nabla \hat{m}_k^H(s^{H}_{k,\ast}) \|_{x_k^H} \leq C_{\mathrm{R}} \, \theta^H \, \|s^{H}_{k,\ast}\|_{x_k^H}.
\end{equation}

To bound the first term in \eqref{eq:triangle_ineq}, we observe that, since \begin{equation}\label{eq:grad_coarse_model_zero}
    \nabla\hat{\psi}_k^H(0_{x_k^H}) - \widetilde{\mathcal I}_h^H\nabla\hat{f}_{k}^h(0_{x_k^h}) = 0, 
\end{equation}
 it holds
\begin{equation}\label{eq:diff-grads}
    \widetilde{\mathcal I}_h^H\nabla\hat{f}_{k}^h(s_k^h) -\nabla\hat{\psi}_k^H(s^{H}_{k,\ast}) = \widetilde{\mathcal I}_h^H \bigl[\nabla\hat{f}_{k}^h(s_k^h)-\nabla\hat{f}_{k}^h(0_{x_k^h})\bigr] - \bigl[\nabla\hat{\psi}_k^H(s^{H}_{k,\ast}) -\nabla\hat{\psi}_k^H(0_{x_k^H})\bigr].
\end{equation}
Since the expression for the gradient of $\hat{\psi}_k^H$ is
\[
	\nabla\hat{\psi}_k^H(s^H) = \nabla\hat{f}_{k}^H(s^H) + v^H_{k},
\]
we get 
\[
    \nabla\hat{\psi}_k^H(s^{H}_{k,\ast}) -\nabla\hat{\psi}_k^H(0_{x_k^H}) = \nabla\hat{f}_{k}^H(s^{H}_{k,\ast})-\nabla\hat{f}_{k}^H(0_{x_k^H})
\]
and, using the Lipschitz continuity of \Cref{ass:pullback_grad_lipschitz} with $t=1$,
\[
\begin{aligned}
    \bigl\| \widetilde{\mathcal I}_h^H\nabla\hat{f}_{k}^h(s_k^h) -\nabla\hat{\psi}_k^H(s^{H}_{k,\ast}) \bigr\|_{x^H_k} &\leq \bigl\| \widetilde{\mathcal I}_h^H [\nabla\hat{f}_{k}^h(s_k^h)-\nabla\hat{f}_{k}^h(0_{x_k^h})] \bigr\|_{x^H_k} +\bigl\| \nabla\hat{\psi}_k^H(s^{H}_{k,\ast}) -\nabla\hat{\psi}_k^H(0_{x_k^H}) \bigr\|_{x^H_k}\\
    &\leq C_\mathcal{I} L_h\|s_k^h\|_{x^h_k}+L_H\|s^{H}_{k,\ast}\|_{x^H_k}.
\end{aligned}
\]
Since $ s_k^h = \widetilde{\mathcal I}_H^h s^{H}_{k,\ast} $, then
\begin{equation}\label{eq:bound_on_shk}
    \|s_k^h\|_{x^h_k} = \|\widetilde{\mathcal I}_H^h s^{H}_{k,\ast}\|_{x^h_k} \leq C_\mathcal{I}\|s^{H}_{k,\ast}\|_{x^H_k},
\end{equation}
and we obtain
\[
    \bigl\| \widetilde{\mathcal I}_h^H\nabla\hat{f}_{k}^h(s_k^h) -\nabla\hat{\psi}_k^H(s^{H}_{k,\ast}) \bigr\|_{x^H_k} \leq \bigl( C_{\mathcal{I}}^2 L_h  +L_H \bigr)\|s^{H}_{k,\ast}\|_{x^H_k}.
\]
Finally, we can conclude\footnote{Note that inequality \eqref{eq:bound_restr_grad_pb_fine} is analogous to \cite[(4.12)]{Calandra:2021}.}
\begin{equation}\label{eq:bound_restr_grad_pb_fine}
    \bigl\| \widetilde{\mathcal I}_h^H\nabla\hat{f}_{k}^h(s_k^h) \bigr\|_{x^H_k} \leq \bigl( C_{\mathcal{I}}^2 L_h +L_H+\lambda_{\max} + C_{\mathrm{R}} \, \theta^H \bigr)\|s^{H}_{k,\ast}\|_{x^H_k}.
\end{equation}

Overall, from~\Cref{ass:pullback_grad_lipschitz}  with $t=1$, we can write:
\begin{align*}
    \|\widetilde{\mathcal I}_{h}^{H} \operatorname{grad} f^h(x^h_k)\|_{x^H_k} &\le \|\widetilde{\mathcal I}_{h}^{H} \operatorname{grad}f^h(x^h_k) - \widetilde{\mathcal I}_{h}^{H}  \nabla \hat{f}_{k}^h(s^h_k)\|_{x^H_k} + \|\widetilde{\mathcal I}_{h}^{H}  \nabla \hat{f}_{k}^h(s^h_k)\|_{x^H_k}\notag  \\
    &= \|\widetilde{\mathcal I}_{h}^{H} \nabla \hat{f}_{k}^h(0_{x^h_k})   - \widetilde{\mathcal I}_{h}^{H}  \nabla \hat{f}_{k}^h(s^h_k)\|_{x^H_k} + \|\widetilde{\mathcal I}_{h}^{H}  \nabla \hat{f}_{k}^h(s^h_k)\|_{x^H_k} \notag \\
    &\le C_{\mathcal{I}} L_h \|s^h_k\|_{x^h_k} + \|\widetilde{\mathcal I}_{h}^{H}\nabla \hat{f}_{k}^h(s^h_k)\|_{x^H_k}\notag \\
    & \le ( C_{\mathcal{I}}^2L_h + C_{\mathcal{I}}^2 L_h + L_H + \lambda_{\max} + C_{\mathrm{R}} \, \theta^H)\|s^{H}_{k,\ast}\|_{x^H_k} \coloneqq c_{1} \, \|s^{H}_{k,\ast}\|_{x^H_k},
\end{align*}
where we used \eqref{eq:bound_on_shk} and \eqref{eq:bound_restr_grad_pb_fine} in the last inequality.
Then by (\ref{eq:go_down_cond}) we get \eqref{eq:bound_s_ast_sq_2}.
\end{proof}

We are now ready to state our convergence result in the following theorem.

\begin{theorem}[Global Convergence]\label{thm: Global CV}
    Let \Cref{ass:smooth,ass:low,ass:pullback_grad_lipschitz,ass:jacbound} hold. Let \(\{x_k^h\}\) be the sequence of fine-level iterates generated by \Cref{alg:riemannian_two-level_ar1}. Then, \(\norm{\operatorname{grad} f^h(x_k^h)}_{x_k^h}\) converges to zero on the subsequence of successful iterations. 
\end{theorem}

\begin{proof}
Consider the sequence of successful iterations and divide them into two groups: $\mathcal{S}^{h}$ when the fine level was used and $\mathcal{S}^{H}$ when the lower level was employed. Let $k_{1}$ be the index of the first successful iteration. From the relations above, it follows that:
\begin{equation}\label{eq:telescoping_bound}
    \begin{aligned}
        f^h(x_{k_1}^h) - f_{\text{low}} &\geq f^h(x_{k_1}^h) - \lim_{k \rightarrow \infty} \inf f^h(x_k^h) \\
        &\geq \sum_{k \in \mathcal{S}^H\cup \mathcal{S}^h} \left( \hat{f}_{k}^h(0_{x_k^h}) - \hat{f}_{k}^h(s_k^h) \right) \\
        &\geq \eta_1 \sum_{k \in \mathcal{S}^H} \left(\psi_k^H(0_{x_k^H})-\psi_k^H(s^{H}_{k,\ast})\right) + \eta_1 \sum_{k \in \mathcal{S}^h} \left(T_k^h(0_{x_k^h})-T_k^h(s_k^h)\right) \\
        &\underset{\eqref{eq:15a},\,\eqref{eq:pred_red_coarse}}{\geq} \frac{\eta_1 }{2} \left( \sum_{k \in \mathcal{S}^H} \lambda_{k}\norm{s^{H}_{k,\ast}}^2_{x_k^H} + \sum_{k \in \mathcal{S}^h} \lambda_k\norm{s_k^h}^2_{x_k^h}\right) \\
        &\underset{\eqref{eq:19},\,\eqref{eq:bound_s_ast_sq_2}}{\geq} \frac{\eta_1 }{2} \left( \sum_{k \in \mathcal{S}^H}\frac{\lambda_{\min}\kappa_{H}^2}{c_1^2} \norm{\grad f^h(x_k^h)}_{x_k^h}^2 + \sum_{k \in \mathcal{S}^h} \frac{1}{\lambda_k}\norm{\grad f^h (x_k^h)}_{x_k^h}^2\right) \\
        &= 
        \frac{\eta_1 }{2}
        \min \Bigg\{
         \frac{\lambda_{\min}\kappa_{H}^2}{c_1^2},\frac{1}{\lambda_{\max}}\Bigg\} \left( \sum_{k \in \mathcal{S}^H} \norm{\grad f^h(x_k^h)}_{x_k^h}^2 + \sum_{k \in \mathcal{S}^h} \norm{\grad f^h (x_k^h)}_{x_k^h}^2\right),
    \end{aligned}
\end{equation}
where $f_{\text{low}}$ is the lower bound of $f^h$ guaranteed by \Cref{ass:low}, and the global upper bound for the regularization parameter is given by \eqref{eq:lambda_max}. 

Hence, we conclude that the series  \(\sum_{\mathcal{S}^h\cup\mathcal{S}^H} \norm{\grad f^h (x_k^h)}_{x_k^h}\) is a bounded series, and therefore has a convergent subsequence. Then, \(\norm{\grad f^h (x_k^h)}_{x_k^h}\)  converges to zero on the subsequence of successful iterations. 

\end{proof}

\subsection{Worst-case complexity} \label{subsec:worst_case_complexity}
We now evaluate the worst-case iteration complexity of our method, that is, the number of iterations needed in the worst-case scenario to obtain $ \norm{\operatorname{grad} f^h(x_k^h)}_{x_k^h} \leq \varepsilon$ for any $\varepsilon>0$.


\begin{theorem}[Worst-case complexity]\label{thm worst case complexity \algname{}}
    Let \Cref{ass:smooth,ass:low,ass:pullback_grad_lipschitz,ass:jacbound} hold and let $k_{1}$ denote the index of the first successful iteration in \Cref{alg:riemannian_two-level_ar1}. Then, given an absolute accuracy level $\varepsilon>0$, \Cref{alg:riemannian_two-level_ar1} needs at most
    \begin{equation}\label{eq:bound_com}
             c \, \frac{f^{h}(x^{h}_{k_1})-f_{\text{low}}}{\varepsilon^2} \left( 1+\frac{| \log \gamma_1|}{\log \gamma_3} \right) +   \frac{1}{\log \gamma_3}\log \! \left(\frac{\lambda_{\max}}{\lambda_0}\right)
    \end{equation}
     iterations in total to produce an iterate $x_{\varepsilon}^h$ such that $\norm{\operatorname{grad} f^h(x_{\varepsilon}^h)}_{x_{\varepsilon}^h}\leq \varepsilon$, where
    \[
        c \coloneqq \frac{\eta_1 }{2}\min \Bigg\{
         \frac{\lambda_{\min}\kappa_{H}^2}{c_1^2},\frac{1}{\lambda_{\max}}\Bigg\},
    \]
    and $\kappa_{H}$, $\lambda_{\max}$, $c_1$ are the constants defined in \eqref{eq:go_down_cond}, \eqref{eq:lambda_max} and \eqref{eq:def_c1} respectively, while $\eta_1$ and $\lambda_{\min}$ are defined in \Cref{alg:riemannian_two-level_ar1}.  
\end{theorem}

\begin{proof}
We first bound the number of successful iterations.
Let $k_{1}$ denote the index of the first successful iteration and let  $k_{\ast}$ be the index of the first successful iterate for which \( \norm{\operatorname{grad} f^h(x_k^h)}_{x_k^h}> \varepsilon\) and \( \norm{\operatorname{grad} f^h(x_{k+1}^h)}_{x_{k+1}^h}< \varepsilon\). Define \(\mathcal{S} \coloneqq \{ 0<k \leq k_{\ast} \colon \rho_{k} \geq \eta_1\} \) the set of successful iterations before $k_{\ast}$. 
It then follows from \eqref{eq:telescoping_bound} that
\begin{equation*}
    \begin{aligned}
        f^h(x_{k_1}^h) - \underset{k \rightarrow\infty}{\lim} f^h(x_k^h) &\geq f^h(x_{k_1}^h)-f^h(x_{{k_{\ast}+1}}^h)= \underset{k \in \mathcal{S}}{\sum}f^h(x_k^h)-\hat{f}_{k}^h(s_k^h)\\ 
        &
        \geq \frac{\eta_1 }{2}\min \Bigg\{
         \frac{\lambda_{\min}\kappa_{H}^2}{c_1^2},\frac{1}{\lambda_{\max}}\Bigg\}|\mathcal{S}| \, \varepsilon^2 \coloneqq c \, |\mathcal{S}| \, \varepsilon^2,
    \end{aligned}
\end{equation*}
which gives a bound on the total number of successful iterations. 

By denoting $\mathcal{U}$ the set of unsuccessful iterations, the updating rule of $\lambda_k$ gives (by induction)
\[
    \lambda_0 \gamma_1^{|\mathcal{S}|}\gamma_3^{|\mathcal{U}|} \leq \lambda_{k_{\ast}} \leq \lambda_{\max},
\]
so,
\[
    |\mathcal{S}|\log \gamma_1 + |\mathcal{U}| \log \gamma_3 \leq \log \! \left(\frac{\lambda_{\max}}{\lambda_0}\right).
\]
Given that $\gamma_1 <1$, 
\[
    |\mathcal{U}| \leq \frac{1}{\log \gamma_3}\log \! \left(\frac{\lambda_{\max}}{\lambda_0}\right) + |\mathcal{S}| \frac{| \log \gamma_1|}{\log \gamma_3},
\]
we obtain that $|\mathcal{S}|+|\mathcal{U}|$ is bounded by \eqref{eq:bound_com}.
\end{proof}

\Cref{thm worst case complexity \algname{}} shows that under suitable assumptions, incorporating lower-level steps does not degrade the overall complexity of the method, preserving the optimal $O(\varepsilon^{-2})$ bound of the RAR1 algorithm (see \Cref{sec:compRAR1}). 

\section{The geometry of fixed-rank matrices} \label{sec:geom_fixed_rank_manif}
So far, our presentation focused on a generic manifold. In this section we specifically consider the manifold of fixed-rank matrices $\mathcal{M}_{r}$, which will be the focus of the numerical experiments in \Cref{sec:numerical_experiments}. In the first subsection, we briefly recall the geometry of this manifold, while in the second subsection, we discuss \Cref{ass:pullback_grad_lipschitz,ass:jacbound} for the manifold of fixed-rank matrices. 

\subsection{Background material} \label{subsec:background}
The manifold of fixed-rank matrices is defined as
\begin{equation}\label{eq:fixed_rank}
    \cMr = \lbrace X \in \mathbb{R}^{m \times n} \colon \rank(X) = r \rbrace.
\end{equation}
Using the singular value decomposition (SVD), one has the equivalent characterization 
\begin{equation*}
    \begin{split}
       \cMr = \lbrace U\Sigma V\tr \colon & U \in \mathrm{St}^{m}_{r}, \ V \in \mathrm{St}^{n}_{r}, \\
              & \Sigma = \operatorname{diag}(\sigma_{1}, \sigma_{2}, \ldots, \sigma_{r} ) \in \R^{r \times r}, \ \sigma_{1} \geq \cdots \geq \sigma_{r} > 0 \rbrace,
    \end{split}
\end{equation*}
where $ \mathrm{St}^{m}_{r} $ is the Stiefel manifold of $ m \times r $ real matrices with orthonormal columns (and analogously $ \mathrm{St}^{n}_{r} $), and  $ \operatorname{diag}(\sigma_{1}, \sigma_{2}, \ldots, \sigma_{r} ) $ is a square matrix with $ \sigma_{1}, \sigma_{2}, \ldots, \sigma_{r} $ on its main diagonal.

Vandereycken~\protect{\cite[Prop.~2.1]{Vandereycken:2013}} showed that $ \cMr $ is a smooth submanifold of dimension $ (m+n-r)r $ embedded in $ \mathbb{R}^{m \times n} $ with a compact representation for its tangent space $ \mathrm{T}_{X}\mathcal{M}_{r} $ at $ X = U\Sigma V\tr \in \cMr $. Indeed, every tangent vector $Z \in \mathrm{T}_{X}\mathcal{M}_{r}$ can be written as
\begin{equation}\label{eq:tan_vec_format_small_param}
Z = UMV\tr + U_{\mathrm{p}} V\tr + UV_{\mathrm{p}}\tr,
\end{equation}
with $M \in \R^{r \times r}$, $U_{\mathrm{p}} \in \R^{m\times r}$, $V_{\mathrm{p}} \in \R^{n \times r}$ such that $U_{\mathrm{p}}\tr U = V_{\mathrm{p}}\tr V = 0$ (known as \emph{gauging conditions} \protect{\cite[\S 9.2.3]{UschmajewV:2020}}). Since $ \cMr \subset \mathbb{R}^{m \times n} $, we represent tangent vectors in~\eqref{eq:tan_vec_format_small_param} as matrices of the same dimensions.

The Riemannian metric on $ \mathcal{M}_{r} $ is the restriction of the Euclidean metric on $ \mathbb{R}^{m \times n} $ to the submanifold $ \cMr $, i.e.,
\[
    g_{X}(W,Z) = \langle W, Z \rangle_{\mathrm{F}} = \trace(W\tr Z), \quad \text{with} \ X \in \cMr \ \text{and} \ W, Z \in \mathrm{T}_{X}\mathcal{M}_{r}.
\]
The norm induced by this metric on each tangent space coincides with the Frobenius norm, denoted by $ \| \cdot \|_{\mathrm{F}} $.

Defining $ P_{U} = UU\tr $ and $ P_{U}^{\perp} = I - P_{U} $ for any $ U \in \mathrm{St}^{m}_{r} $, and similarly for $ V \in \mathrm{St}^{n}_{r} $, the orthogonal projection onto the tangent space at $ X $ is \protect{\cite[(2.5)]{Vandereycken:2013}}
\begin{equation}\label{eq:def_proj}
   \operatorname{P}_{X} \colon \mathbb{R}^{m \times n} \to \mathrm{T}_{X}\mathcal{M}_{r}, \quad Z \mapsto P_{U} Z P_{V} + P_{U}^{\perp} Z P_{V} + P_{U} Z P_{V}^{\perp}.
\end{equation}
For embedded submanifolds, the Riemannian gradient of $f$ at $X$ is given by the orthogonal projection onto the tangent space at $X$ of the Euclidean gradient of $ \bar{f} $, a smooth extension of $f$ to the embedding space $ \mathbb{R}^{m \times n} $; see, e.g., \protect{\cite[(3.37)]{AMS:2008}}:
\begin{equation}\label{eq:Riemannian_grad_as_projection}
    \grad f(X) = \operatorname{P}_{X} \! \big( \nabla \bar{f}(X) \big).
\end{equation}

As mentioned in \Cref{sec:preliminaries}, to map tangent-space updates onto the manifold, we need retractions. A retraction $ \operatorname{R}_{X} $ is a smooth map from the tangent space to the manifold, $ \operatorname{R}_{X}\colon \mathrm{T}_{X}\mathcal{M}_{r} \to \cMr $; it is, essentially, any smooth first-order approximation of the exponential map of the manifold; see, e.g., \protect{\cite{AMS:2008,Absil:2012}}. To establish convergence of Riemannian algorithms, it suffices for the retraction to be defined only locally.
For retractions on the manifolds of fixed-rank matrices, an excellent survey is given in~\cite{AbsilOseledets:2015}; in the experiments in~\Cref{sec:numerical_experiments}, we have chosen the \emph{orthographic retraction}.
In \Cref{app:sec:egrad_psiHk,app:sec:egrad_mHk}, we compute the Euclidean gradients of ambient extensions of the coarse model $\psi_k^H$ and the regularized coarse model $m_k^H$ for this specific choice. The Riemannian gradients are obtained by orthogonally projecting these Euclidean gradients onto the tangent space at the point of evaluation.

\subsection{Discussion on~\texorpdfstring{\Cref{ass:pullback_grad_lipschitz,ass:jacbound}}{Assumption 3 and 4} for the manifold of fixed-rank matrices}\label{subsec:discussion}

In \Cref{sec:theory}, we introduced \Cref{ass:pullback_grad_lipschitz}, which requires Lipschitz bounds on the variation of the pullback gradients from their values at the origin along the computed tangent segments, with constants $L_h$ and $L_H$ independent of $k$, and \Cref{ass:jacbound}, which requires uniform boundedness of the differential of the coarse-level retraction  at the computed coarse steps. We now discuss \emph{sufficient conditions} under which these assumptions hold on the manifold $\cMr$ of fixed-rank matrices defined in \eqref{eq:fixed_rank}.

The main difficulty is that $\cMr$ is not closed: a bounded sequence of rank-$r$ matrices may converge to a matrix of smaller rank. Near such rank-deficient points, the geometry of the manifold can become increasingly ill-conditioned, and the derivatives of retractions need not remain uniformly bounded.

We show that the required assumptions hold when the iterates remain bounded and uniformly separated from the rank-deficient boundary, and the steps remain in a region where the retraction is sufficiently smooth.

\paragraph{Sufficient conditions.}

Let $f$ denote either $f^h$ or $f^H$, and let $\{X_k\}$ and $\{S_k\}$ denote the corresponding sequences of base points and tangent steps generated by \Cref{alg:riemannian_two-level_ar1} with $\cMh$ and $\cMH$ being the fine and coarse approximations of the manifold $\cMr$, respectively. More precisely, these are $(X_k^h,S_k^h)$ at the fine level and $(X_k^H,S_{k,\ast}^H)$ at the coarse level.

We assume the following conditions hold true.

\begin{enumerate}[(i)]
    \item \emph{Boundedness and regularity.} The sequence $\{X_k\}$ is bounded. Moreover, $f$ is the restriction to $\cMr$ of a function, denoted $\bar{f}$, defined on $\mathbb{R}^{m\times n}$ with locally Lipschitz-continuous Euclidean gradient $\nabla \bar{f}$.
    \item \emph{Uniform separation from rank deficiency.} There exists $\sigma_{\min} > 0$ such that
    \[
        \inf_k\sigma_r(X_k) \ge \sigma_{\min} > 0.
    \]
    \item \emph{Control of the step length.} There exists $c \in (0,1)$ such that
    \[
        \|S_k\|_{\mathrm{F}} \le c\,\sigma_r(X_k) \qquad \text{for all }k.
    \]
\end{enumerate}

The boundedness requirement in condition (i) holds for the fine-level iterates if $f^h$ has bounded sublevel sets, since the acceptance mechanism of \algname{} ensures that $f^h(X_{k + 1}^h) \le f^h(X_k^h)$. The corresponding coarse-level base points are bounded whenever the point restriction $\mathcal{I}_h^H$ maps bounded sets to bounded sets, as is the case for the linear restriction operators considered here. The regularity requirement in condition (i) is a separate assumption on the ambient extensions of the objective functions.

Condition (ii) prevents the iterates from approaching the rank-deficient boundary of $\cMr$. This is an assumption on the generated sequences, rather than a property of the objective functions.

Finally, condition (iii) requires the computed steps $S_{k}$ to remain within a sufficiently small neighborhood of their base points $X_{k}$. This is an additional requirement on the steps produced by \algname{}, rather than a consequence of the coarse stopping criterion in \eqref{eq:coarse_stop_crit_twolevel}.

\paragraph{Uniform boundedness of the retraction derivatives.}

Conditions (i) and (ii) ensure that the base points belong to a compact subset of $\cMr$. Moreover, condition (iii) implies that the tangent steps are uniformly bounded, since
\[
    \|S_k\|_{\mathrm{F}} \le c\,\sigma_r(X_k) \le c\,\|X_k\|_{\mathrm{F}}.
\]

Let $\mathcal{B}$ be a compact subset of the ambient space containing the base points and consider
\[
    \mathcal{C} = \big\{(X,S) \in \mathrm T\cMr \colon \;X \in \mathcal{B}, \quad \sigma_r(X) \ge \sigma_{\min}, \quad \|S\|_{\mathrm{F}} \le c\,\sigma_r(X)\big\}.
\]
The set $\mathcal{C}$ is compact: both components are bounded, the singular-value bound prevents loss of rank in limits, and tangency is preserved along such limits. Furthermore, it contains all pairs $(X_k,tS_k)$, for $t \in [0,1]$.

Suppose that the retraction $\operatorname{R}$ is $C^2$ on an open neighborhood of $\mathcal{C}$. By compactness, there exist constants $C_{\mathrm{R}},C_{\mathrm{R}}' > 0$, independent of $k$, such that
\begin{equation}
    \label{eq:uniform_retraction_bounds}
    \|\operatorname{D} \! \operatorname{R}_X(S)\|_{\mathrm{op}} \le C_{\mathrm{R}}, \qquad \text{and} \qquad \|\operatorname{D}^2 \! \operatorname{R}_X(S)\|_{\mathrm{op}} \le C_{\mathrm{R}}',
\end{equation}
for every $(X,S) \in \mathcal{C}$, where both derivatives are taken with respect to the tangent variable and the retraction is viewed as an ambient-matrix-valued map. 

In particular, for the coarse-level step $S_{k,\ast}^H$, the first bound gives
\[
    \|\operatorname{D} \! \operatorname{R}_{X_k^H}(S_{k,\ast}^H)\|_{\mathrm{op}} \le C_{\mathrm{R}},
\]
which establishes \Cref{ass:jacbound}.

We now verify the required smoothness property for the orthographic retraction on $\cMr$.

Let $X = U\Sigma V^{\top}$ and write a tangent vector as $ S = UMV^{\top} + U_{\mathrm{p}}V^{\top} + UV_{\mathrm{p}}^{\top} $, following \eqref{eq:tan_vec_format_small_param}.
The orthographic retraction is
\begin{equation}
    \label{eq:orthographic_retraction_bound}
    \operatorname{R}_X(S) = X + S + U_{\mathrm{p}}(\Sigma + M)^{-1}V_{\mathrm{p}}^{\top},
\end{equation}
whenever $\Sigma + M$ is invertible.
Since
\[
    \sigma_{\min}(\Sigma + M) \ge \sigma_r(X) - \|M\|_2 \ge \sigma_r(X) - \|S\|_{\mathrm{F}},
\]
condition (iii) implies
\[
    \sigma_{\min}(\Sigma + M) \ge (1 - c)\,\sigma_r(X) \ge (1 - c)\,\sigma_{\min} > 0.
\]

Therefore, $\Sigma + M$ is uniformly invertible on the admissible region. Interpreted in smooth local coordinates, the explicit expression \eqref{eq:orthographic_retraction_bound} defines a smooth map wherever $\Sigma + M$ is invertible. The strict lower bound above therefore ensures that the orthographic retraction is smooth on an open neighborhood of $\mathcal{C}$. Consequently, its first and second derivatives admit the uniform bounds in \eqref{eq:uniform_retraction_bounds}.

\paragraph{Lipschitz continuity of the pullback gradients.}

We now establish \Cref{ass:pullback_grad_lipschitz}. Since $\mathcal{C}$ is compact and the retraction is continuous, the set of retracted points $\{ \operatorname{R}_X(S) \colon (X,S) \in \mathcal{C} \}$ is bounded. Hence, there exists a compact convex set $\mathcal{K}$ containing all these retracted points and the corresponding base points.

By condition (i) and the compactness and convexity of $\mathcal{K}$, there exist constants $G,L_{\mathrm{e}} > 0$ such that $\|\nabla \bar{f}(Y)\|_{\mathrm{F}} \le G$ and
\[
    \|\nabla \bar{f}(Y) - \nabla \bar{f}(Z)\|_{\mathrm{F}} \le L_{\mathrm{e}}\|Y - Z\|_{\mathrm{F}},
\]
for all $Y,Z \in \mathcal{K}$.

For each base point $X_k$, consider the pullback $\hat{f}_k(S) = f(\operatorname{R}_{X_k}(S)) = \bar{f}(\operatorname{R}_{X_k}(S))$. The chain rule gives
\[
    \nabla\hat{f}_k(S) = \bigl(\operatorname{D} \! \operatorname{R}_{X_k}(S)\bigr)^{\ast}\big[\nabla \bar{f}(\operatorname{R}_{X_k}(S))\big].
\]
Using $\operatorname{R}_{X_k}(0_{X_{k}}) = X_k$, we obtain
\[
    \nabla\hat{f}_k(S) - \nabla\hat{f}_k(0_{X_{k}}) = \bigl(\operatorname{D} \! \operatorname{R}_{X_k}(S)\bigr)^{\ast}\big[\nabla \bar{f}(\operatorname{R}_{X_k}(S)) - \nabla \bar{f}(X_k)\big] + \big(\operatorname{D} \! \operatorname{R}_{X_k}(S) - \operatorname{D} \! \operatorname{R}_{X_k}(0_{X_{k}})\big)^{\ast} \big[\nabla \bar{f}(X_k)\big].
\]
For every $S$ with $(X_k,S) \in \mathcal{C}$, and in particular for $S = tS_k$, $t \in [0,1]$, the uniform bounds on the retraction derivatives imply
\[
    \|\operatorname{R}_{X_k}(S) - X_k\|_{\mathrm{F}} \le C_{\mathrm{R}}\|S\|_{\mathrm{F}} \quad \text{and} \quad \|\operatorname{D} \! \operatorname{R}_{X_k}(S) - \operatorname{D} \! \operatorname{R}_{X_k}(0_{X_{k}})\|_{\mathrm{op}} \le C_{\mathrm{R}}'\|S\|_{\mathrm{F}}.
\]
Combining these inequalities with the bounds on the ambient gradient, we obtain
\[
    \begin{aligned}
        \|\nabla\hat{f}_k(S) - \nabla\hat{f}_k(0_{X_{k}})\|_{\mathrm{F}}& \le C_{\mathrm{R}}L_{\mathrm{e}}\|\operatorname{R}_{X_k}(S) - X_k\|_{\mathrm{F}} + C_{\mathrm{R}}'G\|S\|_{\mathrm{F}}\\
        & \le (C_{\mathrm{R}}^2L_{\mathrm{e}} + C_{\mathrm{R}}'G)\|S\|_{\mathrm{F}}.
    \end{aligned}
\]
For each fixed $X_k$, the admissible tangent steps form a convex ball. Repeating the same argument with any two steps in this ball in place of $S$ and $0_{X_{k}}$ gives the corresponding two-point Lipschitz estimate. Thus, \Cref{ass:pullback_grad_lipschitz} holds with $L = C_{\mathrm{R}}^2L_{\mathrm{e}} + C_{\mathrm{R}}'G$, which is independent of $k$. Applying this argument separately to $f^h$ and $f^H$ yields the constants $L_h$ and $L_H$ required in \Cref{ass:pullback_grad_lipschitz}.

\paragraph{Geometric interpretation.}
The role of condition (ii) can also be understood through the curvature of the fixed-rank manifold. In general, the normal component of the acceleration of a smooth curve $ c \colon \mathbb{R} \to \mathcal{M} $ on a submanifold is governed by the second fundamental form $ \mathrm{I\!I} $ of the manifold, through the relation\footnote{See, e.g., \cite[end of \S 4]{Olikier:2023} and \cite[(5.44)]{boumal2023introduction}.}
\[
    c''(t) = \ddot{c}(t) + \mathrm{I\!I}(c'(t),c'(t)),
\]
where $c''(t)$ is the acceleration in the ambient space, $\ddot c = \operatorname{P}_{c(t)}(c''(t))$ is the intrinsic Riemannian acceleration\footnote{Here, $\operatorname{P}_{c(t)}$ denotes the orthogonal projector onto the tangent space at $c(t)$.}, and $\mathrm{I\!I}(c'(t),c'(t))$ is the normal component of $c''(t)$. The second fundamental form is in turn related to the \emph{Weingarten map}\footnote{Equivalent definitions appear in the classical references \cite[Def.~18]{ONeill1983}, \cite[p.~128]{doCarmo:1992} and \cite[p.~47, Eq.~(3.24)]{Sakai:1996}, although under different names and notations.} \cite[Def. 1]{AMT:2013}. The Weingarten map of the submanifold $\mathcal{M}$ at $X$ is the operator $\mathfrak{A}_{X} \colon \mathrm{T}_{X}\mathcal{M} \times \mathrm{T}_{X}^{\perp}\mathcal{M} \to \mathrm{T}_{X}\mathcal{M} $, defined by
\[
    \mathfrak{A}_{X}(Z,V) = - \operatorname{P}_{X} \mathrm{D}_{Z} \bar{V},
\]
where $\bar{V}$ is any local extension of $V$ to a normal vector field on $\mathcal{M}$. For all $Z,Y \in \mathrm{T}_X\mathcal{M}$ and $W \in \mathrm{T}_X^\perp\mathcal{M}$, the following relationship between $\mathrm{I\!I}$ and $\mathfrak{A}_{X}$ holds:
\[
    \langle \mathrm{I\!I}_X(Z,Y),W \rangle = \langle (\operatorname{D}_Z\operatorname{P})Y,W \rangle = \langle Y,(\operatorname{D}_Z\operatorname{P})W \rangle = \langle Y,\mathfrak{A}_{X}(Z,W) \rangle.
\]

In particular, for $\cMr$, the last equation in \cite[\S 4.5]{AMT:2013} links the Weingarten map of $\cMr$ at $X$ and the Moore--Penrose pseudo-inverse of $X$, i.e., $ X^{+} = V\Sigma^{-1}U^{\top} $:
\[
    \mathfrak{A}_X(Z,W) = \operatorname{P}_{X} (\operatorname{D}_{Z} \operatorname{P}) W = W Z^{\top} (X^{+})^{\top} + (X^{+})^{\top} Z^{\top} W.
\]
This expression shows how $X^{+}$ governs the possible growth of the Weingarten map (and the second fundamental form). Indeed, since $ X^{+} = V\Sigma^{-1}U^{\top} $, its norm is governed by $1/\sigma_{r}(X)$. More precisely, $\|X^{+}\|_2 = 1/\sigma_r(X)$, so that
\[
    \|\mathfrak{A}_X(Z,W)\|_{\mathrm{F}} \le 2\,\|Z\|_{\mathrm{F}}\|W\|_{\mathrm{F}}/\sigma_r(X).
\]
This curvature information enters the second-order behavior of $f$ on $\cMr$ through the normal component of the Euclidean gradient\footnote{Recall from~\eqref{eq:Riemannian_grad_as_projection} that the Riemannian gradient is $\operatorname{P}_{X}(\nabla \bar{f}(X))$.}. To illustrate this, suppose additionally that $\bar{f}$ is $C^2$. Then, for $Z \in \mathrm{T}_X\cMr$, the Riemannian Hessian of $f$ at $X$ is~\cite[(10)]{AMT:2013}
\[
    \operatorname{Hess} f(X)[Z] = \operatorname{P}_{X} \! \big(\nabla^2 \bar{f}(X)[Z]\big) + \mathfrak{A}_X\big(Z,\operatorname{P}_{X}^{\perp}(\nabla \bar{f}(X))\big),
\]
where $ \nabla^2 \bar{f}(X) $ is the Euclidean Hessian of $\bar{f}$.
If the Euclidean Hessian satisfies $\|\nabla^2 \bar{f}(X)\|_{\mathrm{op}} \le L_{\mathrm{e}}$ on $\mathcal{K}$, then, using $\|\nabla \bar{f}(X)\|_{\mathrm{F}} \le G$, we obtain, for $X \in \mathcal{K}\cap\cMr$,
\[
    \|\operatorname{Hess} f(X)\|_{\mathrm{op}} \le L_{\mathrm{e}} + \frac{2\|\nabla \bar{f}(X)\|_{\mathrm{F}}}{\sigma_r(X)} \le L_{\mathrm{e}} + \frac{2G}{\sigma_r(X)}.
\]
This estimate illustrates the role of $1/\sigma_r(X)$ in second-order bounds and the importance of remaining away from the rank-deficient boundary.

Specifically, if the curve under consideration is a retraction curve, then the normal component of the second derivative of a retraction at the origin is governed by the second fundamental form. Its tangential component also depends on the chosen retraction. Thus, the Weingarten estimate provides a geometric explanation for the possible deterioration of second-order bounds, whereas the uniform retraction bounds and the Lipschitz estimate for the pullback gradients follow from the compactness and smoothness argument above, under the regularity assumption in condition (i).

Conditions (i)--(iii) therefore provide sufficient (but not necessary) conditions for \Cref{ass:pullback_grad_lipschitz,ass:jacbound} to hold on $\cMr$ when using the orthographic retraction.

\section{Numerical experiments}\label{sec:numerical_experiments}

In this section, we illustrate and assess the effectiveness of our \algname{} method across a range of numerical experiments using three test problems, solving their discretizations on the manifold of fixed-rank matrices $ \mathcal{M}_{r} $. We begin by introducing the formulation of the test problems in \Cref{subsec:test_problems}. In \Cref{subsec:retraction_choice}, we recall the orthographic retraction on $ \mathcal{M}_{r} $ and its inverse. Then, in \Cref{subsec:multilevel_implementation}, we describe the actual \textit{multilevel} implementation (extension of \Cref{alg:riemannian_two-level_ar1} to multiple levels) used in the numerical experiments. In \Cref{subsec:setup_and_quantities}, we present the experimental setup and the quantities monitored across the experiments. We then examine convergence within individual cycles (\Cref{subsec:convg_behavior}) and under mesh refinement (\Cref{subsec:mesh_study}). We compare \algname{} against the plain one-level RAR1 (\Cref{subsec:comparisons_rar1}), and RMGLS (\Cref{subsec:comparisons_rmgls}). Finally, in \Cref{subsec:rank_adaptivity}, we present an automatic rank-adaptive version of \algname{}.

\subsection{Description of the test problems} \label{subsec:test_problems}
The first two test problems below (LYAP and NPDE), already considered in~\cite{sutti2021riemannian,Sutti2024}, arise from the discretization of variational problems. Their traditional Euclidean counterparts have been used as benchmarks for nonlinear multilevel methods in the classical literature~\cite{Henson:2003,Gratton:2008,WenGoldfarb:2010}. The third test problem is the stationary Allen--Cahn equation \cite{DaiLiLuong:2020,Bonheure:2022}.
For background on the underlying variational formulations, see~\cite{Brenner:2007,LeDret:2016}.

All test problems are posed on $\Omega=(0,1)^2$, with homogeneous Dirichlet  boundary conditions and the source term is 
\[
    \gamma(x,y) = e^{x-2y} \sum_{q=1}^{5} 2^{q-1} \sin(q\pi x)\sin(q\pi y).
\]
The Lyapunov (LYAP) problem is
\begin{equation}\label{eq:var_pb_1}
    \min_{u\in H_0^1(\Omega)} f(u), \qquad f(u) \coloneqq \int_{\Omega} \left[ \tfrac{1}{2}\|\nabla u\|^2-\gamma u \right]\dx\dy.
\end{equation}
Its variational derivative is
\begin{equation}\label{eq:var_pb_1_grad}
    \frac{\delta f}{\delta u} = -\Delta u-\gamma,
\end{equation}
so its stationary points satisfy the Poisson equation $-\Delta u=\gamma$.

The nonlinear PDE (NPDE) benchmark is based on the functional
\begin{equation}\label{eq:modified_NPDE_variational}
    f(u) \coloneqq \int_{\Omega} \left[ \tfrac{1}{2}\|\nabla u\|^2 + \lambda u^2\left(\tfrac{1}{3}u+\tfrac{1}{2}\right) - \gamma u \right]\dx\dy, \qquad u\in H_0^1(\Omega).
\end{equation}
The parameter $\lambda$ in this functional is fixed at $10$. Its variational derivative is
\[
    \frac{\delta f}{\delta u} = -\Delta u+\lambda u(u+1)-\gamma,
\]
and stationarity therefore yields the semilinear elliptic equation $-\Delta u+\lambda u(u+1)=\gamma$.

The stationary Allen--Cahn equation (for brevity, SACE) is posed on the same square domain $\Omega = (0,1)^2$, with homogeneous Dirichlet boundary conditions, but with zero forcing. The corresponding functional is \cite[(3.1)]{DaiLiLuong:2020}
\[
    f(u) \coloneqq \int_{\Omega} \left[ \tfrac{\epsilon}{2}\|\nabla u\|^2 - \tfrac{1}{2}u^2 + \tfrac{1}{4}u^4 \right]\dx\dy, \qquad u \in H_0^1(\Omega).
\]
The diffusion coefficient $\epsilon$ is fixed at $0.02$. The variational derivative is
\[
    \frac{\delta f}{\delta u} = - \epsilon \Delta u + u^3 - u,
\]
so stationarity yields the equation $- \epsilon \Delta u + u^3 - u = 0$, with homogeneous Dirichlet boundary conditions; see \cite[(4.5)--(4.6)]{DaiLiLuong:2020}.

For all three test problems, discretizing and adding a rank constraint give optimization problems of the form
\[
    \min_{X \in \cMr} f(X),
\]
where $X$ is the matrix that represents the discretization of the solution $u$ on a uniform square grid, and, with a little abuse of notation, we still use $f$ to denote the discrete counterpart of the objective function. Further details of the discretizations for LYAP and NPDE are given in~\cite[Sec.~5]{sutti2021riemannian}.

The choice of $r$ affects both the accuracy attainable under the rank constraint and the cost of each iteration. Since the rank needed to achieve a prescribed accuracy may be unknown a priori, we also consider a rank-adaptive variant of \algname{} in \Cref{subsec:rank_adaptivity}. This variant starts at a low rank and increases it if the ambient gradient norm at the end of a fixed-rank stage exceeds the prescribed tolerance.

\subsection{Choice of retraction on \texorpdfstring{$\mathcal{M}_{r}$}{cMr}} \label{subsec:retraction_choice}
Throughout \Cref{sec:numerical_experiments}, we use the orthographic retraction on $\mathcal{M}_{r}$~\cite[Sec.~3.2]{AbsilOseledets:2015}. Here, we briefly recall its expression, suppressing level superscripts and iteration subscripts for readability. As in \Cref{subsec:background}, $X = U \Sigma V^{\top}\in\mathcal{M}_r$, where $U,V\in\mathbb{R}^{n\times r}$ have orthonormal columns and $\Sigma\in\mathbb{R}^{r\times r}$ is nonsingular. Let $Z\in\mathrm{T}_X\mathcal{M}_r$ have the representation $(M,U_{\mathrm p},V_{\mathrm p})$ introduced in~\eqref{eq:tan_vec_format_small_param}. Provided that $\Sigma+M$ is nonsingular, the orthographic retraction $ \operatorname{R}_{X} \colon \mathrm{T}_X\mathcal{M}_{r} \to \mathcal{M}_{r} $ is defined by
\[
    \mathrm{R}_X(Z) = \bigl(U(\Sigma+M)+U_{\mathrm p}\bigr)(\Sigma+M)^{-1}\bigl((\Sigma+M)V^{\top}+V_{\mathrm p}^{\top}\bigr).
\]
Given $X \in \mathcal{M}_{r}$, we denote by $\operatorname{R}^{-1}_{X} \colon \mathcal{U}_{X} \to \mathrm{T}_{X} \mathcal{M}_{r}$ the \emph{inverse orthographic retraction}\footnote{See, e.g., \cite[\S 3.2.3]{AbsilOseledets:2015} or~\cite[(2.9)]{sutti2021riemannian}.}, where $\mathcal{V}_{X} \subseteq \mathcal{M}_{r}$ is a neighborhood of $X$ on which the orthographic retraction is a diffeomorphism, as in \Cref{subsec:coarse_step}. For $Y \in \mathcal{U}_{X}$, it is given by
\begin{equation}\label{eq:def_inv_orth_retr}
    \operatorname{R}^{-1}_{X}(Y) = \operatorname{P}_{X}(Y-X) = \operatorname{P}_{X}(Y) - X,
\end{equation}
$\operatorname{P}_{X} \colon \mathbb{R}^{m \times n} \to \mathrm{T}_{X} \mathcal{M}_{r}$ being the orthogonal projector onto $\mathrm{T}_{X} \mathcal{M}_{r}$. In this work, we choose the orthographic retraction because it has explicit closed-form formulas for both the retraction and its inverse, which facilitate related calculations and implementations.

\subsection{Multilevel implementation} \label{subsec:multilevel_implementation}
In our theoretical discussion in \Cref{sec:preliminaries,sec:theory}, we presented and analyzed the convergence of a Riemannian two-level AR1 method. In contrast, in the numerical experiments of this section, we use an actual \emph{multilevel version}. We describe it here, using the notation for matrices of the fixed-rank manifold $\mathcal{M}_{r}$.

We consider two-dimensional elliptic PDE problems discretized on a hierarchy of grids, indexed by $\ell = \ell_{\mathrm{c}},\ldots,\ell_{\mathrm{f}}$, where $\ell_{\mathrm{f}}$ and $\ell_{\mathrm{c}}$ denote the finest and coarsest levels, respectively. At level $\ell$, the number of grid points in each coordinate direction is $n_{\ell} = 2^{\ell} + 1$, including boundary points. We write $f^\ell$ and $\mathcal{M}_{r}^{\ell}$ for the discrete objective and fixed-rank manifold at this resolution. The finest-level objective is $f^{\ell_{\mathrm{f}}}$, acting on matrices of size $n_{\ell_{\mathrm{f}}}\times n_{\ell_{\mathrm{f}}}$.

At each level $\ell_{\mathrm{c}} < \ell \leq \ell_{\mathrm{f}}$, the main step of an \algname{} cycle is either a current-level step, called a \emph{fine step}, or a \emph{recursive coarse step} exploiting the lower levels of the discretization hierarchy. We denote by $F^\ell$ the objective supplied to the current call at level $\ell$. At the finest level, $F^{\ell_{\mathrm{f}}} = f^{\ell_{\mathrm{f}}}$. At a recursive level, $F^\ell$ is the regularized model constructed at the immediately finer level and held fixed throughout that call.

For discrete elliptic problems, classical multigrid and multilevel methods combine relaxation, which damps oscillatory error components, with coarse-grid correction, which efficiently treats the remaining smooth components \cite{Nash:2000}. The implementation used in the numerical experiments includes RAR1 pre- and post-smoothing. For clarity, we omit both smoothing stages from the algorithm description below and from the pseudocode.

At the current iterate $X_k^\ell \in \mathcal{M}_r^\ell$ on a level $\ell > \ell_{\mathrm{c}}$, recursion is governed by the multilevel analogue of the two-level condition~\eqref{eq:go_down_cond}:
\begin{equation}\label{eq:go_down_cond_multilevel}
    \begin{aligned}
        \norm{\widetilde{\cI}_{\ell}^{\ell - 1}\grad F^\ell(X_k^\ell)}_{\mathrm{F}} & \geq \kappa_{\ell - 1}\,\norm{\grad F^\ell(X_k^\ell)}_{\mathrm{F}}.
    \end{aligned}
\end{equation}
If this test fails, the algorithm takes a current-level RAR1 trial step on $F^\ell$, as in \Cref{subsec:fine_step}. A rejected trial leaves $X_k^\ell$ unchanged and increases $\lambda_k$ for the next iteration.

If a coarse step is chosen at iteration $k$ on level $\ell$, the algorithm constructs the objective for the immediately coarser level $\ell - 1$. Set $X_{0,k}^{\ell - 1} \coloneqq \mathcal{I}_{\ell}^{\ell - 1}(X_k^\ell)$ and define the first-order coherence correction by
\[
    V_k^{\ell - 1} \coloneqq \widetilde{\mathcal{I}}_{\ell}^{\ell - 1}\bigl(\operatorname{grad} F^\ell(X_k^\ell)\bigr) - \operatorname{grad} f^{\ell - 1}(X_{0,k}^{\ell - 1}).
\]
On the inverse-retraction domain at $X_{0,k}^{\ell - 1}$, define the coarse model and its regularization, in analogy with \eqref{eq:coarse_model_psiHk} and \eqref{eq:reg_coarse_model}, by
\[
    \begin{aligned}
        \psi_k^{\ell - 1}(X^{\ell - 1}) & \coloneqq f^{\ell - 1}(X^{\ell - 1}) + \innerproduct{V_k^{\ell - 1}}{\operatorname{R}^{-1}_{X_{0,k}^{\ell - 1}}(X^{\ell - 1})}_{\mathrm{F}},\\
        m_k^{\ell - 1}(X^{\ell - 1}) & \coloneqq \psi_k^{\ell - 1}(X^{\ell - 1}) + \frac{\lambda_k}{2}\lVert \operatorname{R}^{-1}_{X_{0,k}^{\ell - 1}}(X^{\ell - 1}) \rVert_{\mathrm{F}}^2.
    \end{aligned}
\]
In \Cref{alg:rmar1-cycle}, we use the pullbacks $\hat{\psi}_k^{\ell - 1} \coloneqq \psi_k^{\ell - 1} \circ \operatorname{R}_{X_{0,k}^{\ell - 1}}$ and $\hat{m}_k^{\ell - 1} \coloneqq m_k^{\ell - 1} \circ \operatorname{R}_{X_{0,k}^{\ell - 1}}$, defined on the tangent neighborhood corresponding to the inverse-retraction domain. The models take the anchor $X_{0,k}^{\ell - 1}$ and a tangent step $S^{\ell - 1}$ as arguments. 

Each descent constructs a fresh model from $f^{\ell - 1}$. The coherence and regularization terms in the supplied objective $F^\ell$ contribute to the new coherence correction through $\operatorname{grad} F^\ell(X_k^\ell)$.

The recursive call at level $\ell - 1$ starts at the fixed point $X_{0,k}^{\ell - 1}$ with the supplied objective $F^{\ell - 1} = m_k^{\ell - 1}$. The value of $\lambda_k$ from the calling level is held fixed in this model. The coarser-level solve initializes its own adaptive parameter at $\lambda_0$. The same construction applies whenever a recursive level chooses to descend further. Point restrictions are assumed to yield points on the corresponding coarse manifolds, and retracted trial points are assumed to lie in the domain of the supplied objective for their call.

At every recursive level $\ell_{\mathrm{c}} \leq \ell < \ell_{\mathrm{f}}$, the call stops when the current iterate $X_k^\ell$ satisfies both conditions below (cf.~\eqref{eq:coarse_stop_crit_twolevel} in the two-level version):
\begin{equation}\label{eq:coarse_stop_crit_multilevel}
    F^\ell(X_k^\ell) < F^\ell(X_0^\ell), \qquad \lVert \grad F^\ell(X_k^\ell) \rVert_{\mathrm{F}} \leq \theta^\ell \lVert S^\ell \rVert_{\mathrm{F}},
    \qquad S^\ell \coloneqq \operatorname{R}^{-1}_{X_0^\ell}(X_k^\ell).
\end{equation}
Here, $X_0^\ell$ is the fixed initial point of that call, and $\theta^\ell > 0$. The test follows the adaptive-regularization principle with a first-order model; see \cite[\S 2.4.1]{CartisGouldToint:2022}. On the coarsest level $\ell_{\mathrm{c}}$, we use RAR1 (cf. \Cref{app:RAR1}) on the supplied objective with the same stopping rule. At the finest level, the stopping test is $\lVert \grad f^{\ell_{\mathrm{f}}}(X_k^{\ell_{\mathrm{f}}}) \rVert_{\mathrm{F}} \leq \epsilon^{\ell_{\mathrm{f}}}$.

The recursive call at level $\ell - 1$ returns its total inverse-retraction displacement, denoted by $S_k^{\ell - 1}$ in the calling routine. The calling routine prolongs it as $S_k^\ell = \widetilde{\mathcal{I}}_{\ell - 1}^{\ell}(S_k^{\ell - 1})$. For this coarse trial, the predicted reduction is the decrease in $\psi_k^{\ell - 1}$ from the fixed initial point $X_{0,k}^{\ell - 1}$ to the point obtained by the coarser-level solve. The actual reduction is the decrease in the supplied objective $F^\ell$ from $X_k^\ell$ to $\operatorname{R}_{X_k^\ell}(S_k^\ell)$, as in \Cref{alg:rmar1-cycle}.

\begin{figure}[htb]
    \centering
    \includegraphics[width=\linewidth]{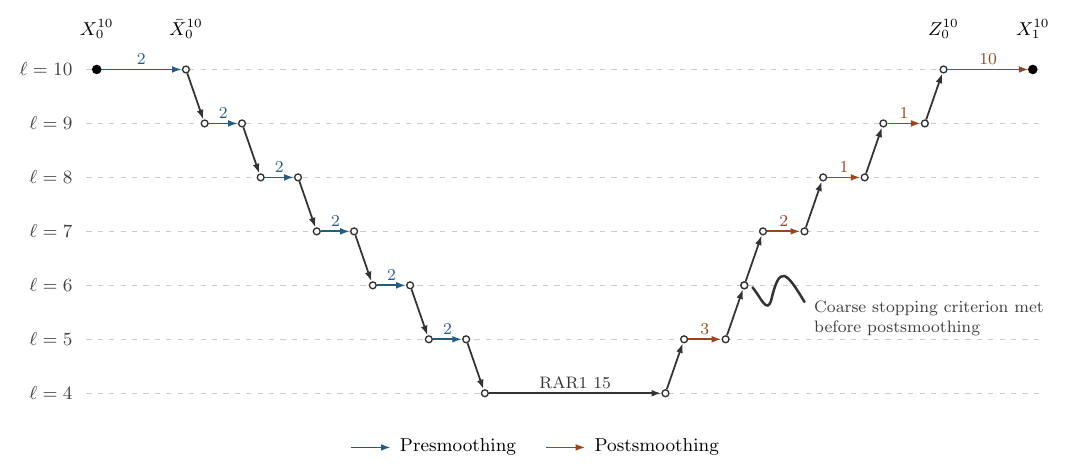}
    \caption{First \algname{} cycle for the LYAP experiment from \Cref{subsec:convg_behavior}.}
    \label{fig:RMAR1_LYAP_first_cycle}
\end{figure}

\Cref{fig:RMAR1_LYAP_first_cycle} illustrates the first \algname{} cycle of the LYAP convergence experiment described in \Cref{subsec:convg_behavior}. The experiment setup is discussed in \Cref{sec:numerical_experiments}, but we anticipate it here to provide an illustration of an actual cycle performed by our \algname{} on a concrete problem. The rank is $r = 5$, and the parameters $\theta^{\ell - 1} = 0.02$ and $\kappa_{\ell - 1} = 0.1$ are used at every coarse level. Blue and orange arrows indicate presmoothing and postsmoothing, respectively, and the horizontal labels count accepted steps. The downward arrows trace the recursive calls from level 10 to level 4, where RAR1 takes 15 accepted steps. On the ascent, each upward arrow indicates the application of the correction returned by the level below. Level 6 meets its coarse stopping criterion immediately after applying the correction and returns before postsmoothing.
Two presmoothing steps from $X_0^{10}$ give the presmoothed iterate $\bar{X}_0^{10}$. Applying the accepted coarse correction to $\bar{X}_0^{10}$ gives $Z_0^{10}$, the starting point for postsmoothing.

A high-level pseudocode for \algname{} is given in \Cref{alg:rmar1-cycle,alg:rmar1-convergenceConditions,alg:rmar1-driver}. \Cref{alg:rmar1-driver} is the main driver of the code, which calls the subroutine $\operatorname{RMAR1\_cycle}$, described in \Cref{alg:rmar1-cycle}. \Cref{alg:rmar1-convergenceConditions} checks the stopping conditions in \eqref{eq:coarse_stop_crit_multilevel} at recursive levels and the absolute gradient tolerance $\epsilon^{\ell_{\mathrm{f}}}$ at the finest level. Please note that the actual MATLAB implementation, publicly available at \url{https://github.com/MarcoSutti/RMAR1}, includes other technicalities omitted from the pseudocode for clarity of exposition.

\clearpage

\begin{adjustbox}{max totalsize={\linewidth}{\textheight},center}
\begin{minipage}{\linewidth}
\begin{algorithm}[H]
    \setlength{\baselineskip}{11pt}
    \setlength{\abovedisplayskip}{3pt}
    \setlength{\belowdisplayskip}{3pt}
    \setlength{\abovedisplayshortskip}{2pt}
    \setlength{\belowdisplayshortskip}{2pt}
    \setlength{\jot}{1pt}
    \DontPrintSemicolon
    \caption{$\operatorname{RMAR1\_cycle}(\ell, F^\ell, X_0^\ell, \{f^j\}_{j = \ell_{\mathrm{c}}}^{\ell})$.}\label{alg:rmar1-cycle}
    \KwIn{The current level $\ell$, the supplied real-valued objective $F^\ell$ on a domain in $\mathcal{M}^\ell$, an initial point $X_0^\ell$ in that domain, and the reference functions $\{f^j\}_{j = \ell_{\mathrm{c}}}^{\ell}$ on their respective manifolds.}
    \KwOut{A point $X_\ast^{\ell_{\mathrm{f}}}$ if $\ell = \ell_{\mathrm{f}}$, or a tangent vector $S^\ell$ at $X_0^\ell$ otherwise; stopping as in \Cref{alg:rmar1-convergenceConditions}.}
    Set $k = 0$ and $\lambda_k = \lambda_0$.\;
    \While{$\operatorname{convergenceConditions}(\ell, \ell_{\mathrm{f}}, F^\ell, X_0^\ell, X_k^\ell, \epsilon^{\ell_{\mathrm{f}}}, \theta^\ell)$}{
        $\bullet$ \textbf{Model choice:} Set $\mathsf{useCoarse} \gets \mathsf{false}$.\;
        \If{$\ell > \ell_{\mathrm{c}}$}{
            Set $X_{0,k}^{\ell - 1} \coloneqq \mathcal{I}_{\ell}^{\ell - 1} (X_k^\ell)$.\;
            Compute $\widetilde{\mathcal{I}}_{\ell}^{\ell - 1}(\operatorname{grad} F^\ell(X_k^\ell))$.\;
            \If{$\lVert \widetilde{\mathcal{I}}_{\ell}^{\ell - 1}(\operatorname{grad} F^\ell(X_k^\ell)) \rVert \geq \kappa_{\ell - 1}\lVert \operatorname{grad} F^\ell(X_k^\ell) \rVert$ \textbf{\textup{and}} $\lVert \widetilde{\mathcal{I}}_{\ell}^{\ell - 1}(\operatorname{grad} F^\ell(X_k^\ell)) \rVert > 0$}{
                $\mathsf{useCoarse} \gets \mathsf{true}$.\;
            }
        }
        \eIf{$\mathsf{useCoarse}$}{
            $\bullet$ \textbf{Coarse-model construction:} Define
            \begin{algomathdisplay}
                \begin{aligned}
                    V_k^{\ell - 1} & \coloneqq \widetilde{\mathcal{I}}_{\ell}^{\ell - 1}(\operatorname{grad} F^{\ell}(X_k^\ell)) - \operatorname{grad} f^{\ell - 1}(X_{0,k}^{\ell - 1}),\\
                    \hat{\psi}_k^{\ell - 1}(S^{\ell - 1}) & \coloneqq f^{\ell - 1}(\operatorname{R}_{X_{0,k}^{\ell - 1}}(S^{\ell - 1})) + \langle V_k^{\ell - 1}, S^{\ell - 1}\rangle,\\
                    \hat{m}_k^{\ell - 1}(S^{\ell - 1}) & \coloneqq \hat{\psi}_k^{\ell - 1}(S^{\ell - 1}) + \frac{\lambda_k}{2}\lVert S^{\ell - 1}\rVert^2.
                \end{aligned}
            \end{algomathdisplay}
            \textbf{Recursive call:} With the defining quantities from level $\ell$ fixed, define $F^{\ell - 1}$ on the inverse-retraction domain at $X_{0,k}^{\ell - 1}$ and compute
            \begin{algomathdisplay}
                \begin{aligned}
                    F^{\ell - 1}\!\left(\operatorname{R}_{X_{0,k}^{\ell - 1}}(S^{\ell - 1})\right) & \coloneqq \hat{m}_k^{\ell - 1}(S^{\ell - 1}),\\
                    S_k^{\ell - 1} & \gets \operatorname{RMAR1\_cycle}\!\left(\ell - 1, F^{\ell - 1}, X_{0,k}^{\ell - 1}, \{f^j\}_{j = \ell_{\mathrm{c}}}^{\ell-1} \right).
                \end{aligned}
            \end{algomathdisplay}
            Set $S_k^\ell = \widetilde{\mathcal{I}}_{\ell - 1}^{\ell} (S_k^{\ell - 1})$.\;
            Set $\mathrm{pred}_k^\ell = \hat{\psi}_k^{\ell - 1}(0) - \hat{\psi}_k^{\ell - 1}(S_k^{\ell - 1})$.\;
        }{
            $\bullet$ \textbf{Taylor step computation:} Set $S_k^\ell = - \operatorname{grad} F^\ell(X_k^\ell)/\lambda_k$.\;
            Set $\mathrm{pred}_k^\ell = \lVert \operatorname{grad} F^\ell(X_k^\ell)\rVert^2/\lambda_k$.\;
        }
        $\bullet$ \textbf{Acceptance of the trial point:} Compute $\displaystyle\rho_k^\ell = \frac{F^\ell(X_k^\ell) - F^\ell(\operatorname{R}_{X_k^\ell}(S_k^\ell))}{\mathrm{pred}_k^\ell}$.\;
        \eIf{$\rho_k^\ell \geq \eta_1$}{
            $X_{k + 1}^{\ell} = \operatorname{R}_{X_k^\ell}(S_k^\ell)$.\;
        }{
            $X_{k + 1}^\ell = X_k^\ell$.\;
        }
        $\bullet$ \textbf{Regularization parameter update:}\;
        \eIf{$\rho_k^\ell \geq \eta_1$}{
            $\displaystyle\lambda_{k + 1} = \begin{cases}
                \max\{\lambda_{\min}, \gamma_2\lambda_k\} & \text{if }\rho_k^\ell \geq \eta_2,\\
                \max\{\lambda_{\min}, \gamma_1\lambda_k\} & \text{if }\rho_k^\ell < \eta_2.
            \end{cases}$\;
        }{
            $\lambda_{k + 1} = \gamma_3\lambda_k$.\;
        }
        $k = k + 1$.\;
    }
    \eIf{$\ell = \ell_{\mathrm{f}}$}{
        \Return{$X_k^\ell$}.\;
    }{
        $S^\ell = \operatorname{R}^{-1}_{X_0^\ell}(X_k^\ell)$.\;
        \Return{$S^\ell$}.\;
    }
\end{algorithm}
\end{minipage}
\end{adjustbox}

\newpage

\begin{algorithm}[htbp]
    \small
    \DontPrintSemicolon
    \LinesNumbered
    \caption{$[\,\mathsf{notConverged}\,] = \operatorname{convergenceConditions}(\ell, \ell_{\mathrm{f}}, F^\ell, X_0^\ell, X_k^\ell, \epsilon^{\ell_{\mathrm{f}}}, \theta^\ell)$.}\label{alg:rmar1-convergenceConditions}
    \KwIn{The finest level $\ell_{\mathrm{f}}$, the current level $\ell$, the supplied objective $F^\ell$, the fixed initial point $X_0^\ell$ and current point $X_k^\ell$, the finest-level tolerance $\epsilon^{\ell_{\mathrm{f}}} > 0$, and $\theta^\ell > 0$.}
    \KwOut{A boolean $\mathsf{notConverged}$, equal to $\mathsf{true}$ if further iterations are required and $\mathsf{false}$ otherwise.}
    $\mathsf{notConverged} \gets \mathsf{false}$.\;
    \eIf{$\ell = \ell_{\mathrm{f}}$}{
        \If{$\lVert\operatorname{grad} F^\ell(X_k^\ell)\rVert > \epsilon^{\ell_{\mathrm{f}}}$}{
            $\mathsf{notConverged} \gets \mathsf{true}$.\;
        }
    }{
        $S^\ell \coloneqq \operatorname{R}^{-1}_{X_0^\ell}(X_k^\ell)$.\;
        \If{$F^\ell(X_k^\ell) \geq F^\ell(X_0^\ell)$ \textbf{\textup{or}} $\lVert\operatorname{grad} F^\ell(X_k^\ell)\rVert > \theta^\ell\lVert S^\ell\rVert$}{
            $\mathsf{notConverged} \gets \mathsf{true}$.\;
        }
    }
    \Return{$\mathsf{notConverged}$}.\;
\end{algorithm}

\begin{algorithm}[htbp]
    \small
    \DontPrintSemicolon
    \LinesNumbered
    \caption{Driver for RMAR1\_cycle.}\label{alg:rmar1-driver}
    \KwIn{The coarsest and finest levels $\ell_{\mathrm{c}}$ and $\ell_{\mathrm{f}}$, a hierarchy $\{f^\ell\colon \mathcal{M}^\ell \to \mathbb{R}\}_{\ell = \ell_{\mathrm{c}}}^{\ell_{\mathrm{f}}}$ with $f^{\ell_{\mathrm{f}}} = f$, point restrictions $\mathcal{I}_{\ell}^{\ell - 1}$ and metric-adjoint tangent transfers $\widetilde{\mathcal{I}}_{\ell}^{\ell - 1}, \widetilde{\mathcal{I}}_{\ell - 1}^{\ell}$ for $\ell = \ell_{\mathrm{c}} + 1,\ldots,\ell_{\mathrm{f}}$, retractions on each level, an initial point $X_0^{\ell_{\mathrm{f}}} \in \mathcal{M}^{\ell_{\mathrm{f}}}$, and a tolerance $\epsilon^{\ell_{\mathrm{f}}} > 0$.}
    \KwOut{An iterate $X_\ast^{\ell_{\mathrm{f}}}$ satisfying $\lVert\operatorname{grad} f(X_\ast^{\ell_{\mathrm{f}}})\rVert \leq \epsilon^{\ell_{\mathrm{f}}}$.}
    Choose $0 < \eta_1 \leq \eta_2 < 1$, $0 < \gamma_2 \leq \gamma_1 < 1 < \gamma_3$, and $\lambda_{\min} > 0$.\;
    Choose $\kappa_{\ell - 1} \in (0,1)$ for $\ell = \ell_{\mathrm{c}} + 1,\ldots,\ell_{\mathrm{f}}$.\;
    Choose $\theta^\ell > 0$ for $\ell = \ell_{\mathrm{c}},\ldots,\ell_{\mathrm{f}}$.\;
    Choose a common initial value $\lambda_0 \geq \lambda_{\min}$ and set $F^{\ell_{\mathrm{f}}} = f^{\ell_{\mathrm{f}}}$.\;
    Make $\ell_{\mathrm{c}}$, $\ell_{\mathrm{f}}$, the transfers, the retractions, the fixed parameters, $\lambda_0$, and $\epsilon^{\ell_{\mathrm{f}}}$ available to every RMAR1\_cycle call.\;
    $X_\ast^{\ell_{\mathrm{f}}} \gets \operatorname{RMAR1\_cycle}\!\left(\ell_{\mathrm{f}}, F^{\ell_{\mathrm{f}}}, X_0^{\ell_{\mathrm{f}}}, \{f^\ell\}_{\ell = \ell_{\mathrm{c}}}^{\ell_{\mathrm{f}}}\right)$.\;
    \Return{$X_\ast^{\ell_{\mathrm{f}}}$}.\;
\end{algorithm}

\subsection{Experimental setup and monitored quantities} \label{subsec:setup_and_quantities}

We perform all the experiments in MATLAB R2025a on an ASUS Vivobook S 16 laptop running Windows 11 Home, equipped with an Intel Core Ultra 5 225H processor, 16 GB of RAM, and an Intel Arc 130T GPU. Computations use the CPU, with MATLAB restricted to a single computational thread.

We used Manopt version 8.0 \cite{manopt} as the computational framework for our MATLAB implementation of \algname{}, which seamlessly integrates with Manopt. We adapted Manopt's embedded fixed-rank manifold factory to incorporate the orthographic retraction and its inverse, and implemented the multilevel scheme and RAR1 algorithm within this framework.

Depending on the experiment, we monitor the following quantities, where $X_k$ denotes the current finest-level iterate:
\begin{itemize}
    \item The relative objective error
    \[
        \text{err-}f \coloneqq \frac{\lvert f(X_k) - f_{\mathrm{ref}}\rvert}{\lvert f_{\mathrm{ref}}\rvert},
    \]
    where $ f_{\mathrm{ref}} $ is the objective value at a highly accurate rank-$r$ reference solution of the same discretized problem, obtained by solving the discretized optimization problem directly on $ \mathcal{M}_{r} $ with the Riemannian ARC solver of \cite{agarwal2021adaptive}, requiring the Riemannian gradient norm to be below $10^{-12}$.
    \item The absolute Riemannian gradient norm
    \[
        \text{R-grad} \coloneqq \|\grad f(X_k)\|_{\mathrm F}.
    \]
\end{itemize}

Unless otherwise stated, we use the following parameter values:
\[
    \lambda_{0} = 10^{-6}, \quad \eta_{1} = 0.25, \quad \eta_{2} = 0.70, \quad \gamma_{1} = 0.85, \quad \gamma_{2} = 0.3, \quad \gamma_{3} = 3, \quad \lambda_{\min} = 10^{-10}.
\]
Moreover, we use $\kappa_{\ell} = 0.1$ in the test~\eqref{eq:go_down_cond_multilevel} that determines whether to attempt a coarse correction and $\theta^{\ell} = 0.02$ in condition \eqref{eq:coarse_stop_crit_multilevel} at each coarse level for all test problems.

In each experiment, we fix the coarsest level at $\ell_{\mathrm{c}} = 4$. To obtain the initial iterate $X_0$ on the finest grid at cycle $k = 0$, we approximately solve the rank-constrained problem on the coarsest grid with Riemannian ARC~\cite{agarwal2021adaptive} and prolong the resulting approximation level by level to the finest grid. This initialization is inspired by the nested-iteration principle underlying full multigrid methods; see, e.g., \cite[\S 2.6]{Trottenberg:2000}.

The \algname{} optimization process is declared converged when R-grad is below $10^{-7}$. The following experimental subsections specify the finest level, rank, smoothing budgets, and initial regularization parameter, and any departures from the above settings.

\subsection{Convergence behavior of \algname{}} \label{subsec:convg_behavior}
\begin{figure}[htb]
  \centering
  \begin{minipage}[t]{0.32\textwidth}
     \centering
     \includegraphics[width=\linewidth]{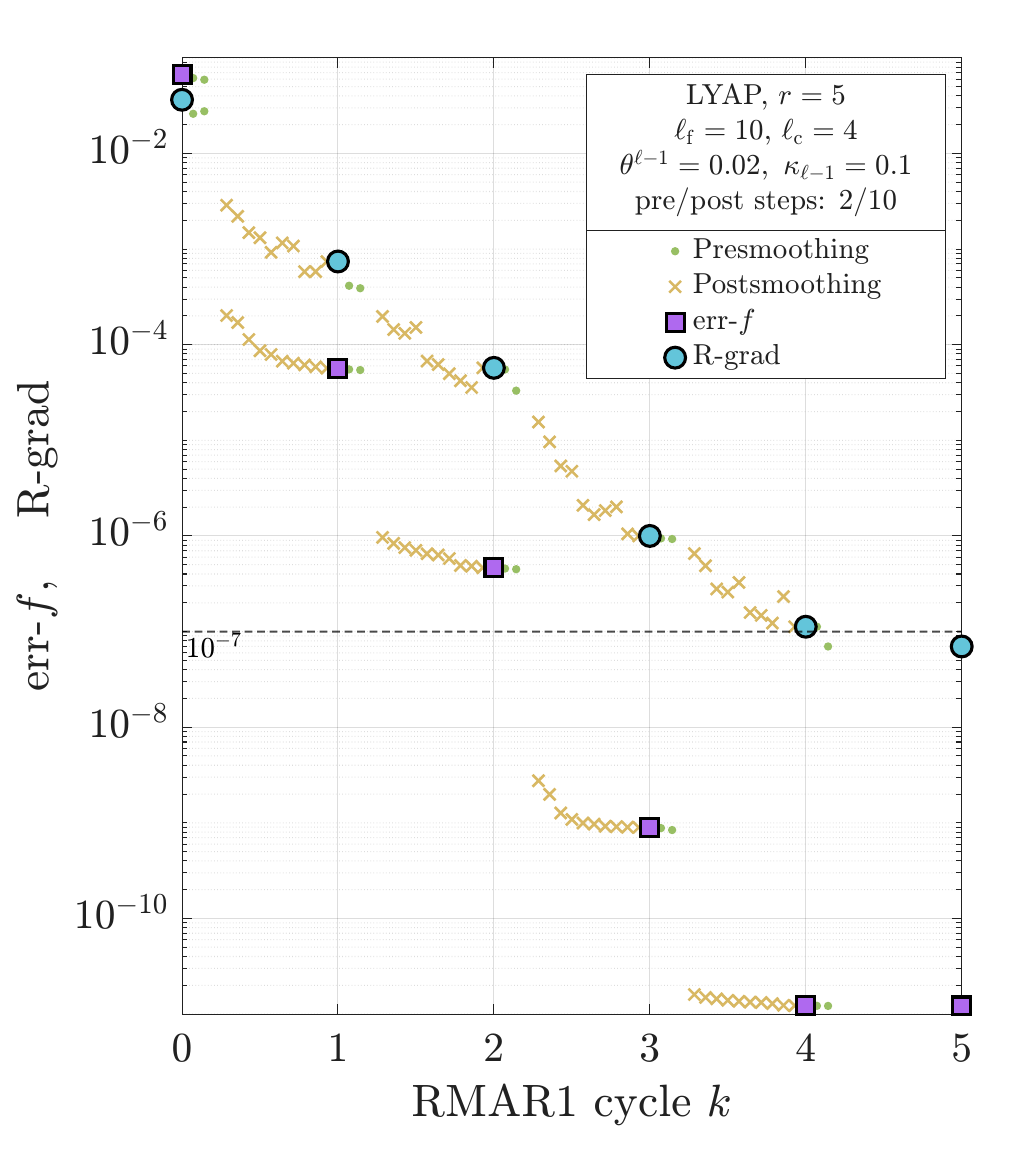}
     \par\smallskip (a)
  \end{minipage}\hfill
  \begin{minipage}[t]{0.32\textwidth}
     \centering
     \includegraphics[width=\linewidth]{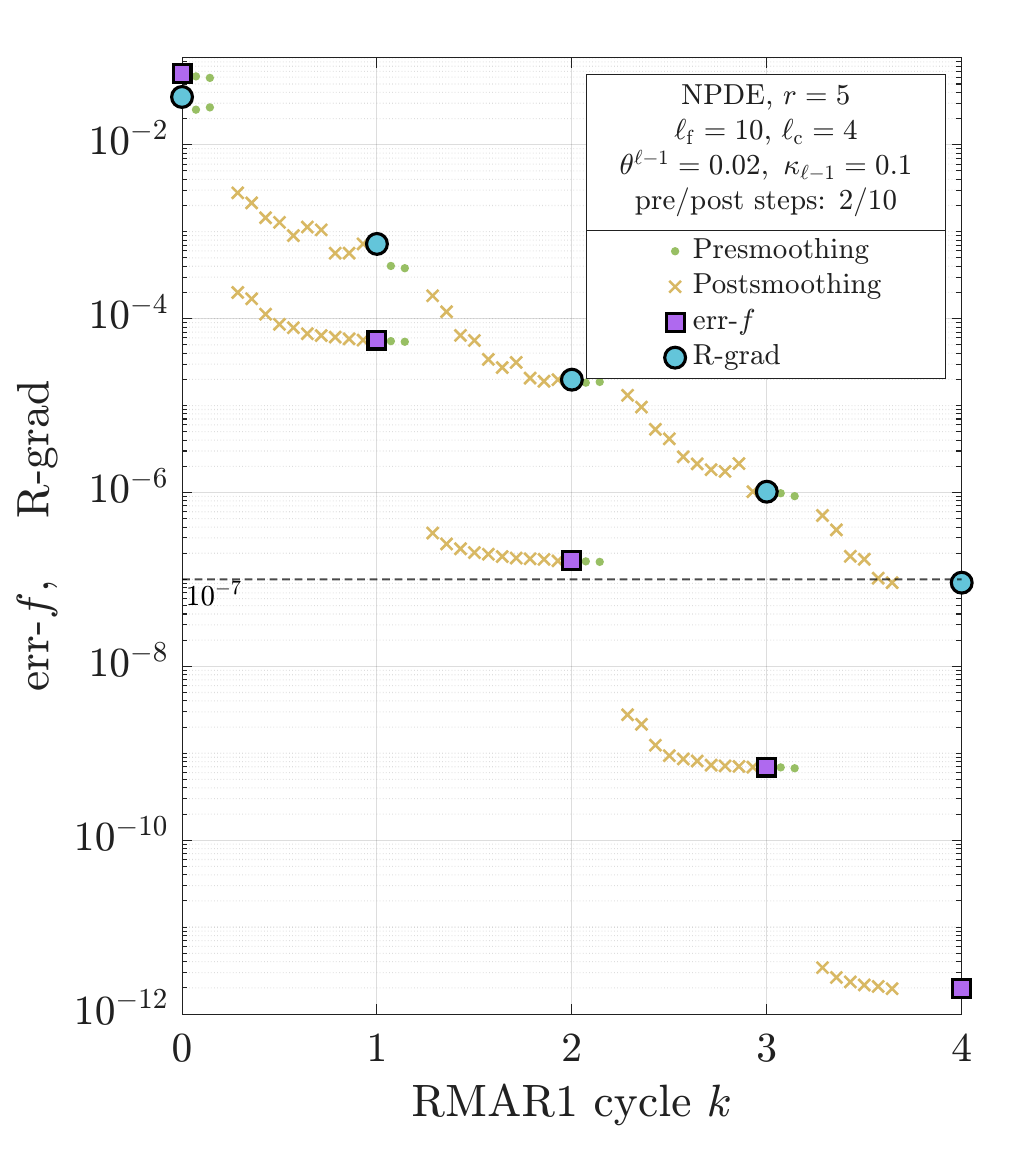}
     \par\smallskip (b)
  \end{minipage}\hfill
  \begin{minipage}[t]{0.32\textwidth}
     \centering
     \includegraphics[width=\linewidth]{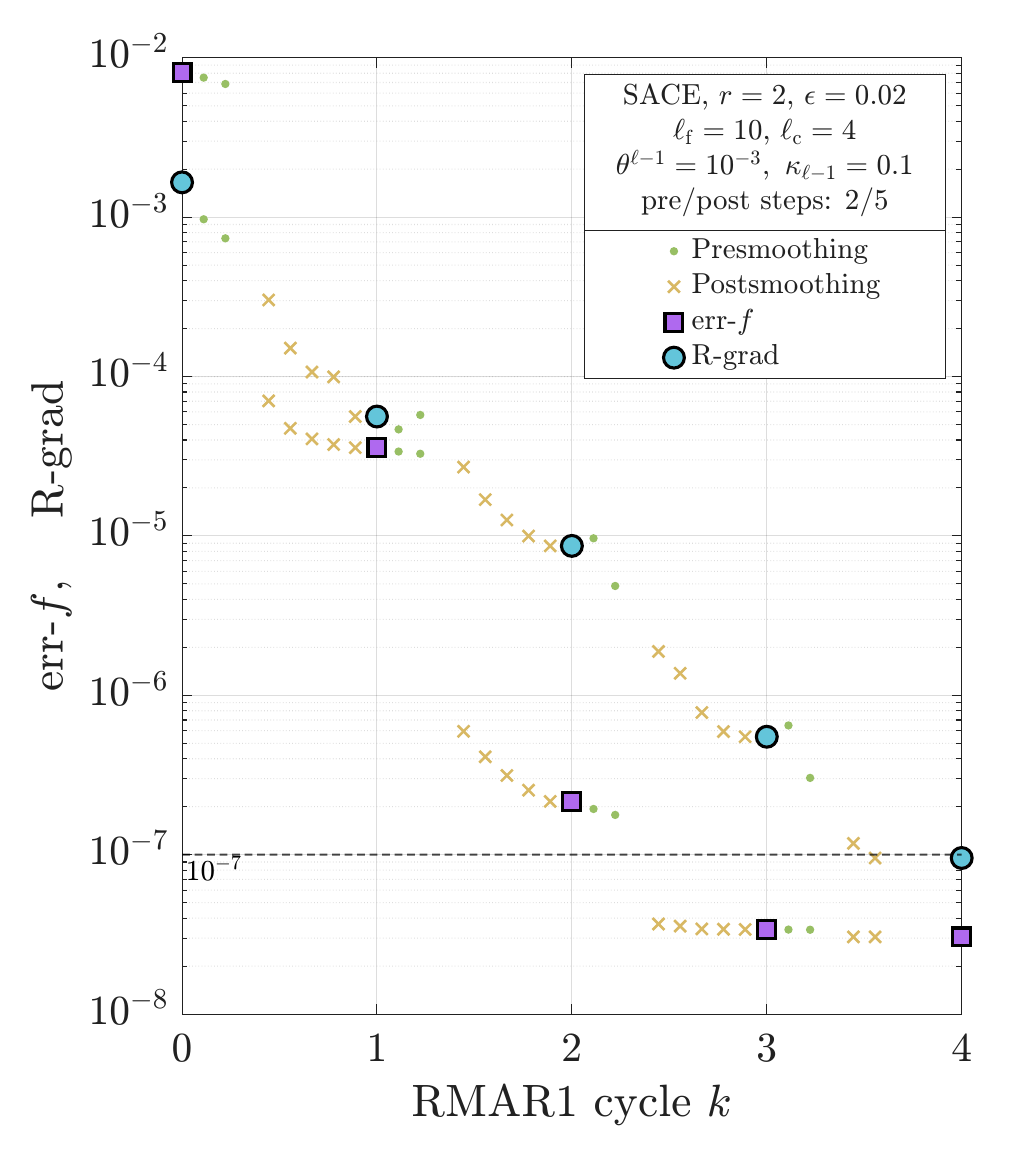}
     \par\smallskip (c)
  \end{minipage}
  \caption{Convergence histories of \algname{} for (a) LYAP at rank $5$, (b) NPDE at rank $5$, and (c) SACE at rank $2$.} \label{fig:rmar1_iteration_history}
\end{figure}

\Cref{fig:rmar1_iteration_history} shows the convergence histories of \algname{} for LYAP, NPDE, and SACE in panels (a), (b), and (c), respectively. All three experiments use $\ell_{\mathrm{f}} = 10$ and $\ell_{\mathrm{c}} = 4$. For LYAP and NPDE, the rank is fixed at $r = 5$, and each cycle allows up to two accepted presmoothing steps and ten accepted postsmoothing steps. We initialize each run by approximately solving the rank-constrained problem on the coarsest grid with ARC and prolonging the resulting point to the finest grid. In all three experiments, \algname{} applies the stopping criterion \eqref{eq:coarse_stop_crit_multilevel} to recursive coarse solves and uses RAR1 as the coarsest-level solver. The coarse stopping factor $\theta^{\ell - 1}$ is $0.02$ for LYAP and NPDE and $10^{-3}$ for SACE. The recursion threshold of \eqref{eq:go_down_cond_multilevel} is $\kappa_{\ell - 1} = 0.1$ in all three experiments. Both parameters are held fixed across the applicable levels.

We plot err-$f$ and R-grad as defined in \Cref{subsec:setup_and_quantities}. The squares and circles represent err-$f$ and R-grad, respectively, at the beginning and at the end of each \algname{} cycle. Dots and crosses show the values of both monitored quantities after accepted presmoothing and postsmoothing steps, respectively. Here, $k$ denotes the cycle index, and $X_k$ is the iterate at the end of cycle $k$. The positions of the intermediate markers within the cycle interval $(k - 1, k)$ schematically indicate the order of the steps.

The intermediate markers reveal changes that would remain hidden if only end-of-cycle values were displayed. In particular, R-grad sometimes increases between accepted smoothing steps even as err-$f$ decreases. Where both groups of markers appear, the change from the last presmoothing marker to the first postsmoothing marker includes the coarse-grid phase and the first accepted postsmoothing step. LYAP and NPDE both reach the absolute gradient tolerance of $10^{-7}$. LYAP reaches it after the second accepted presmoothing step of cycle 5. NPDE reaches it after the sixth accepted postsmoothing step of cycle 4. The final gradient norms are $6.98 \times 10^{-8}$ for LYAP and $9.17 \times 10^{-8}$ for NPDE.

For SACE, we use the diffusion coefficient $\epsilon = 0.02$ and rank $r = 2$. Each \algname{} cycle allows up to two accepted presmoothing steps and five accepted postsmoothing steps. As shown in panel~(c) of \Cref{fig:rmar1_iteration_history}, the run reaches the absolute gradient tolerance of $10^{-7}$ after the second accepted postsmoothing step of cycle 4. The final R-grad and err-$f$ are $9.53 \times 10^{-8}$ and $3.05 \times 10^{-8}$, respectively.

\begin{figure}[htb]
  \centering
  \begin{minipage}[t]{0.32\textwidth}
     \centering
     \includegraphics[width=\linewidth]{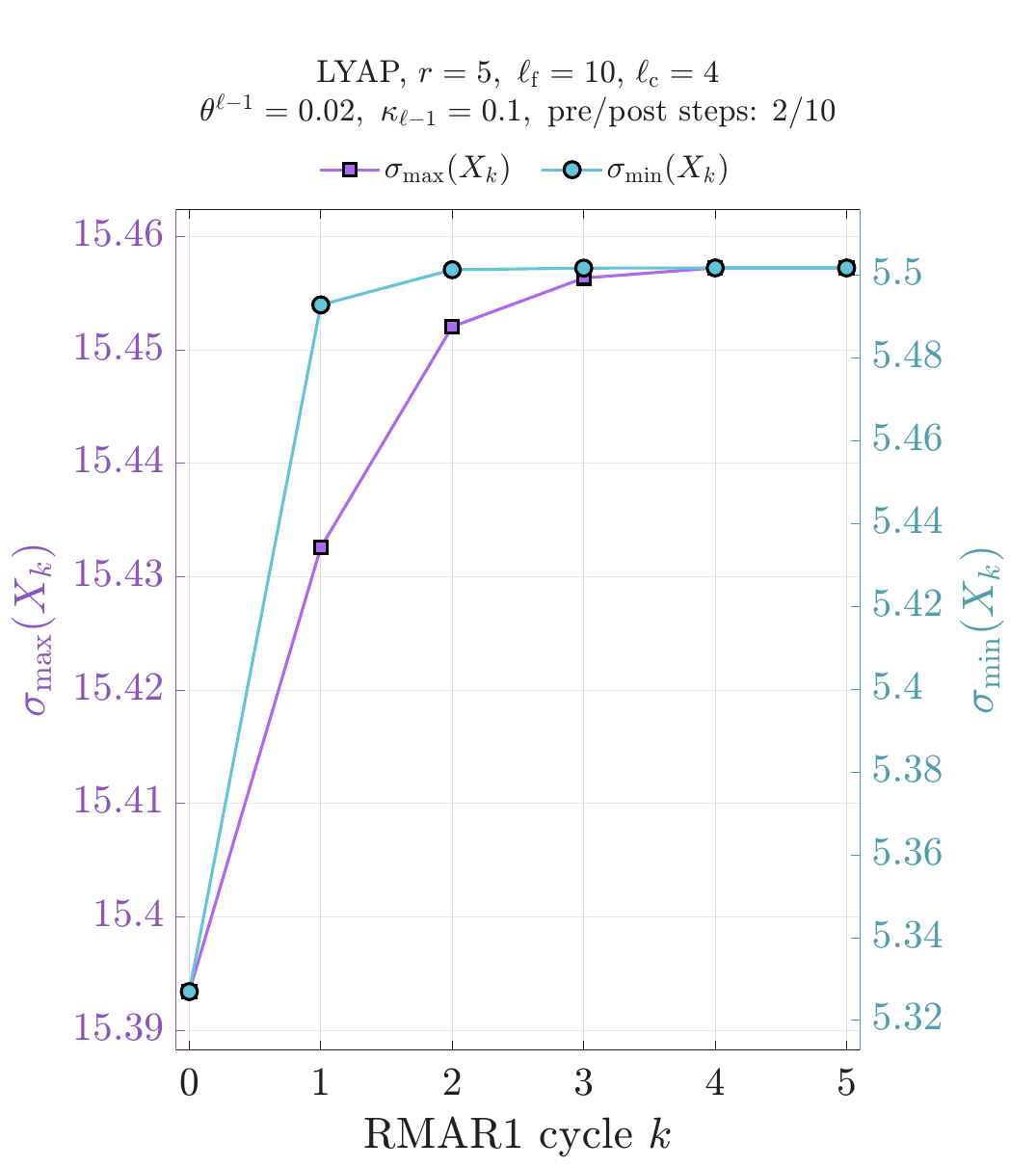}
     \par\smallskip (a)
  \end{minipage}\hfill
  \begin{minipage}[t]{0.32\textwidth}
     \centering
     \includegraphics[width=\linewidth]{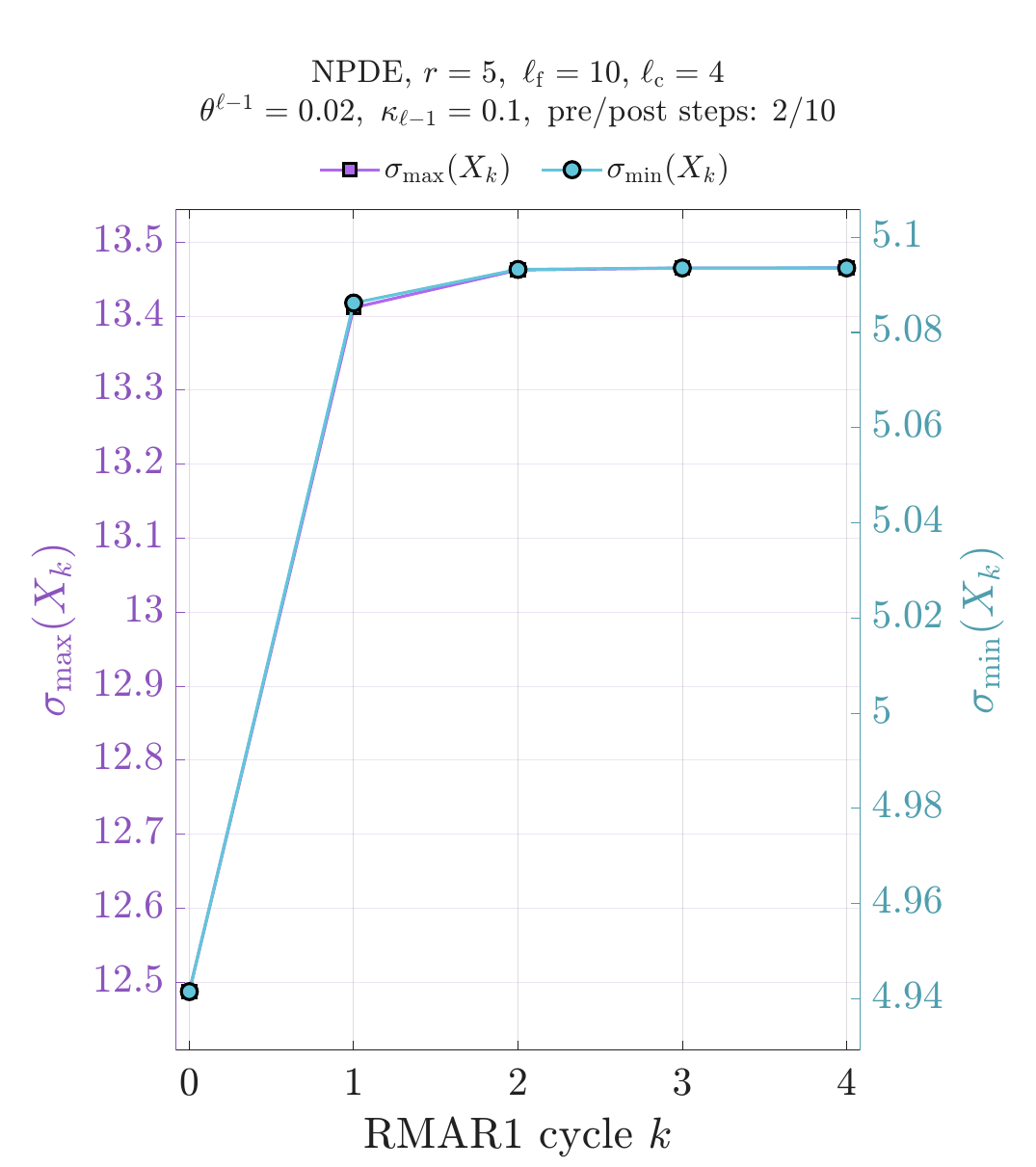}
     \par\smallskip (b)
  \end{minipage}\hfill
  \begin{minipage}[t]{0.32\textwidth}
     \centering
     \includegraphics[width=\linewidth]{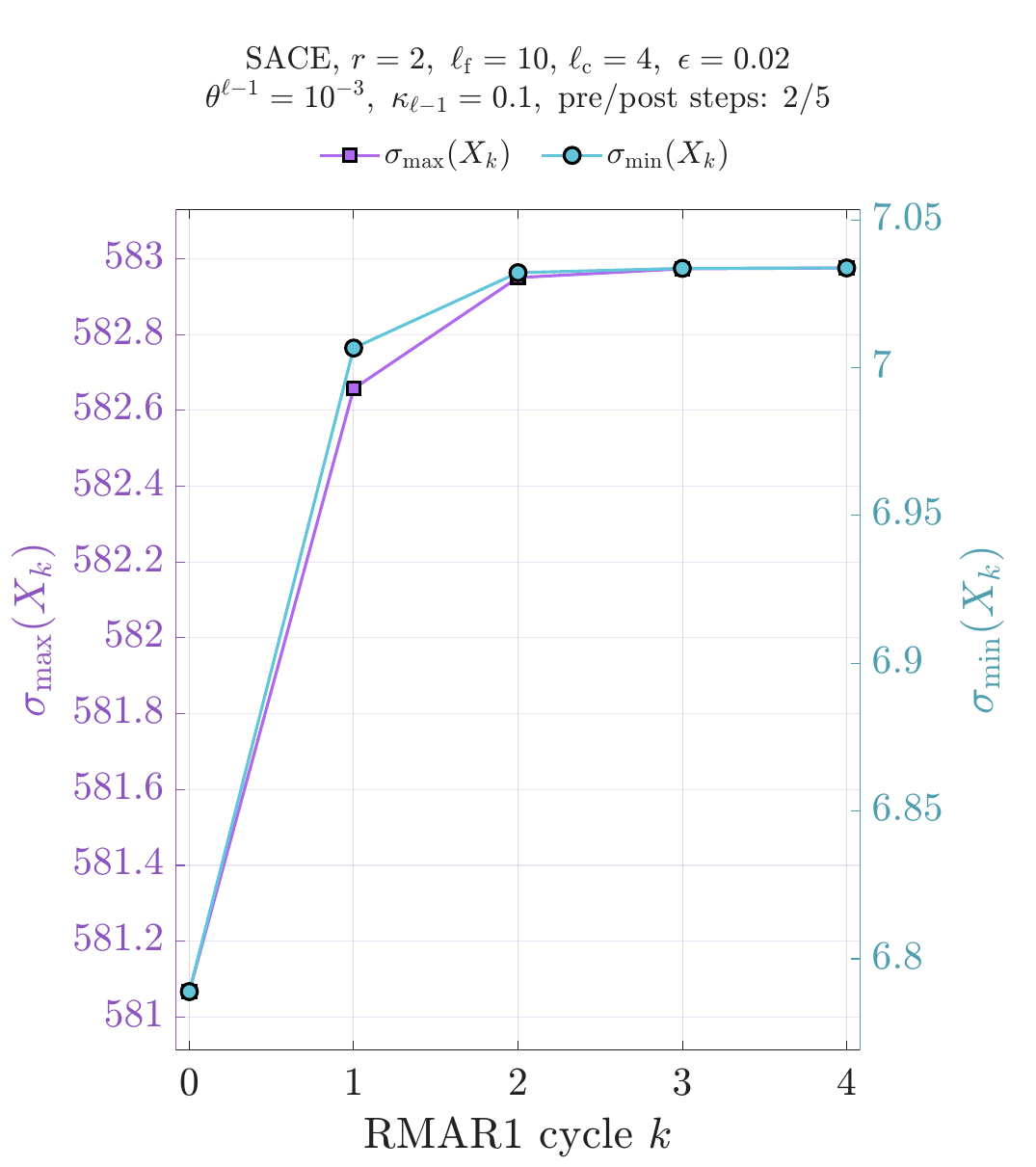}
     \par\smallskip (c)
  \end{minipage}
  \caption{Largest singular value $\sigma_{\max}(X_k)$ (left axis) and smallest nonzero singular value $\sigma_{\min}(X_k)$ (right axis) for (a) LYAP at rank $5$, (b) NPDE at rank $5$, and (c) SACE at rank $2$.} \label{fig:rmar1_singular_values_dual_axis}
\end{figure}

\Cref{fig:rmar1_singular_values_dual_axis} shows the evolution of $\sigma_{\max}(X_k)$ and the smallest nonzero singular value $\sigma_{\min}(X_k)$ for LYAP, NPDE, and SACE in panels (a), (b), and (c), respectively. Each panel uses independent vertical scales, with $\sigma_{\max}(X_k)$ on the left axis and $\sigma_{\min}(X_k)$ on the right. For LYAP and NPDE, both singular values increase overall and change by less than $10^{-5}$ in relative terms between cycles 4 and 5 for LYAP and between cycles 3 and 4 for NPDE. Over the recorded iterates, the smallest values of $\sigma_{\min}(X_k)$ are approximately $5.33$ for LYAP and $4.94$ for NPDE. The ratio $\sigma_{\max}(X_k)/\sigma_{\min}(X_k)$ ranges from approximately $2.81$ to $2.89$ for LYAP and from $2.53$ to $2.64$ for NPDE.
For SACE, $\sigma_{\max}(X_k)$ increases from approximately $581.07$ to $582.98$, while $\sigma_{\min}(X_k)$ increases from approximately $6.79$ to $7.03$. The ratio $\sigma_{\max}(X_k)/\sigma_{\min}(X_k)$ ranges from approximately $82.9$ to $85.6$. For a rank-$r$ iterate, $\sigma_{\min}(X_k) = \sigma_r(X_k)$ is the Frobenius distance from $X_k$ to the set of matrices of rank less than or equal to $r - 1$. Thus, for the five recorded SACE iterates $X_0,\ldots,X_4$, the minimum distance to the set of matrices of rank less than or equal to one is approximately $6.79$. \Cref{fig:rmar1_singular_values_dual_axis} also illustrates that in these experiments the sequence $\{X_k\}$ is bounded over the recorded iterates; see (ii) in \Cref{subsec:discussion}.

\subsection{Quasi-mesh-independent convergence} \label{subsec:mesh_study}
Classical multigrid and multilevel optimization methods aim to limit the dependence of convergence on mesh resolution. In classical multigrid, smoothing and coarse-grid corrections act on complementary error components, which can yield mesh-independent convergence. In this section, we investigate how the number of outer cycles that \algname{} needs to reach a fixed tolerance varies with mesh resolution.

\begin{figure}[H]
  \centering
  \begin{minipage}[t]{0.485\textwidth}
     \centering
     \includegraphics[width=\linewidth,height=0.50\textheight,keepaspectratio]{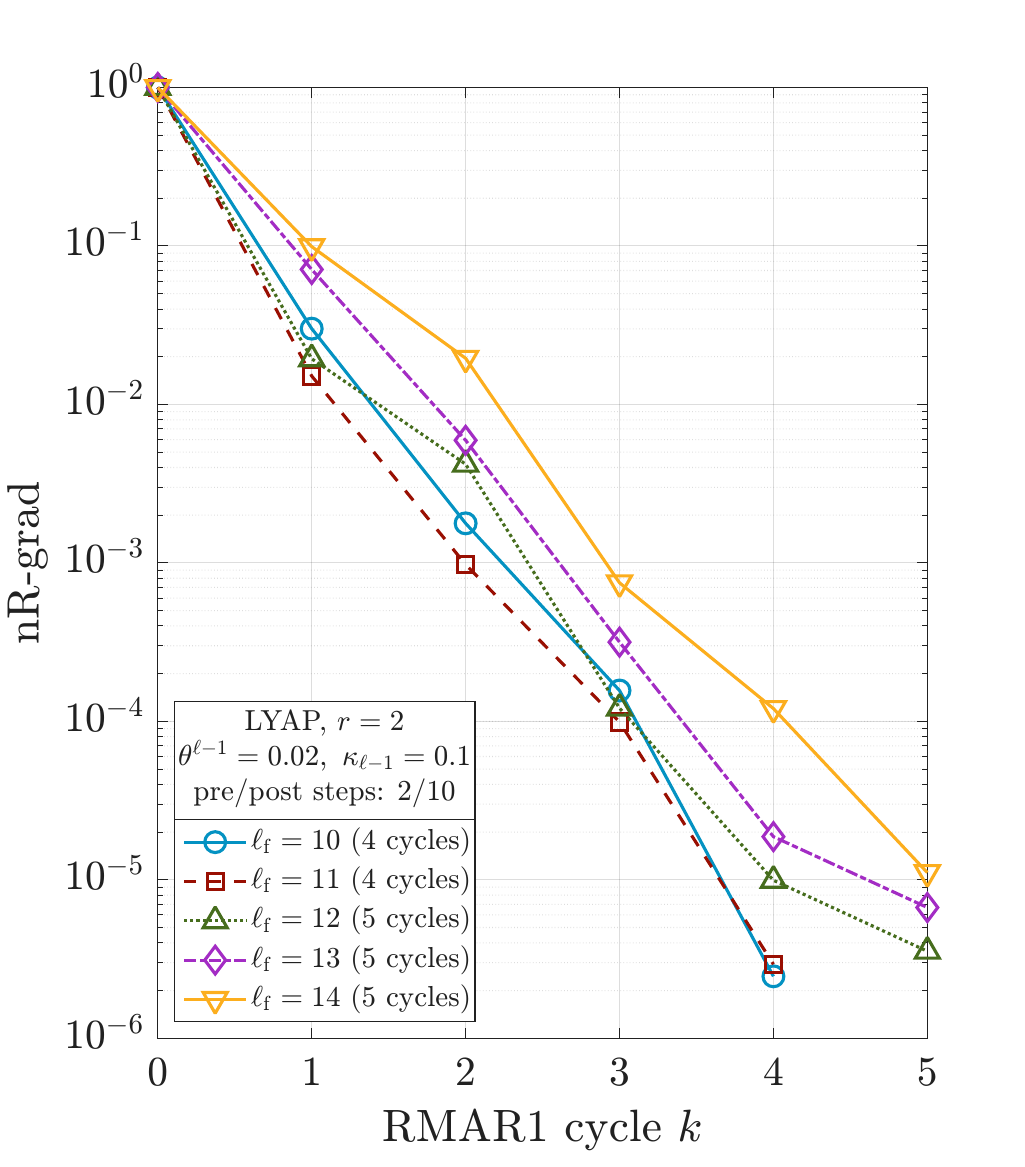}
     \caption{Mesh study for LYAP with rank $r = 2$.}
     \label{fig:rmar1_mesh_study_LYAP_K2}
  \end{minipage}\hfill
  \begin{minipage}[t]{0.485\textwidth}
     \centering
     \includegraphics[width=\linewidth,height=0.50\textheight,keepaspectratio]{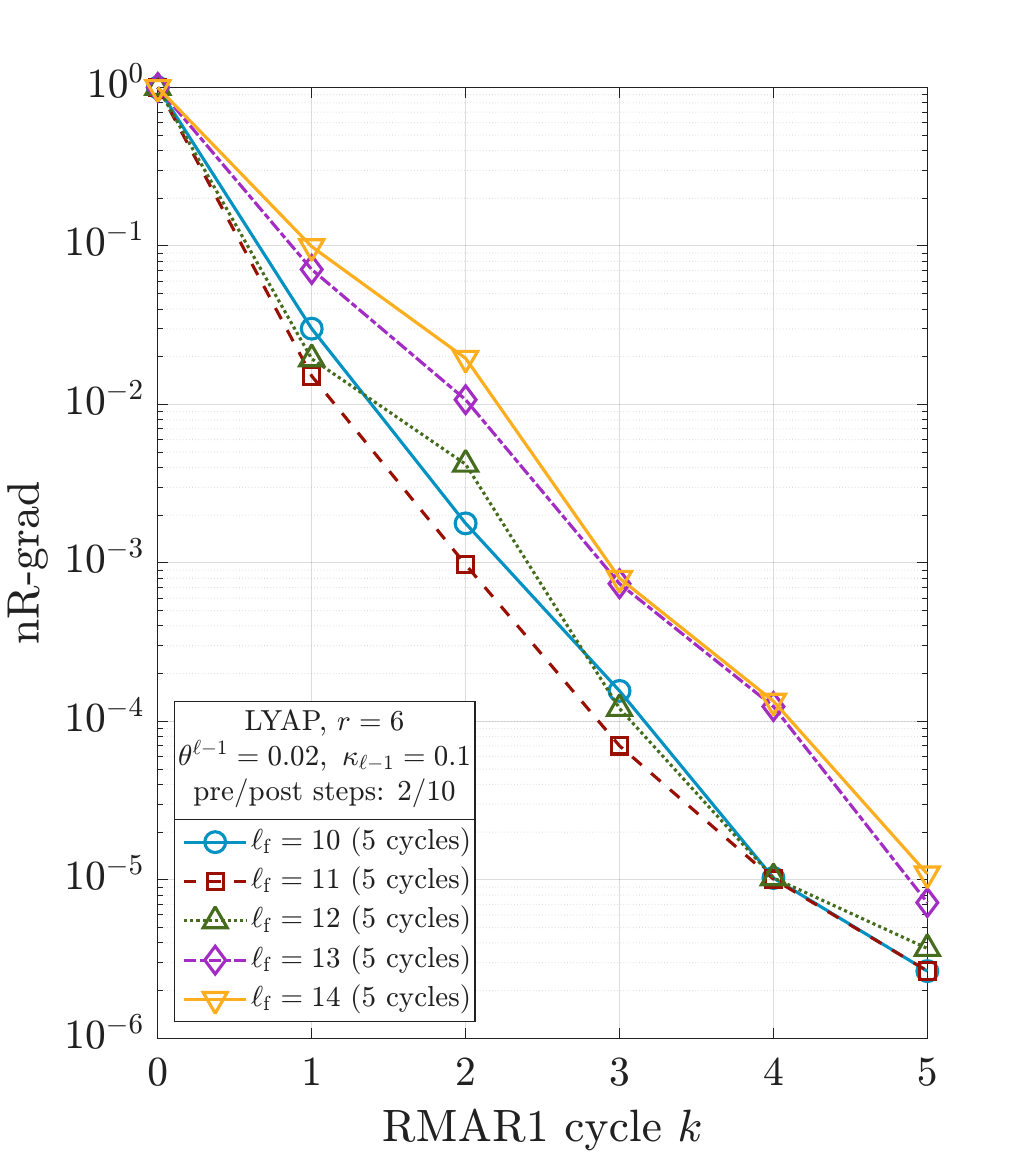}
     \caption{Mesh study for LYAP with rank $r = 6$.}
     \label{fig:rmar1_mesh_study_LYAP_K6}
  \end{minipage}
\end{figure}

\begin{figure}[H]
  \centering
  \begin{minipage}[t]{0.485\textwidth}
     \centering
     \includegraphics[width=\linewidth,height=0.50\textheight,keepaspectratio]{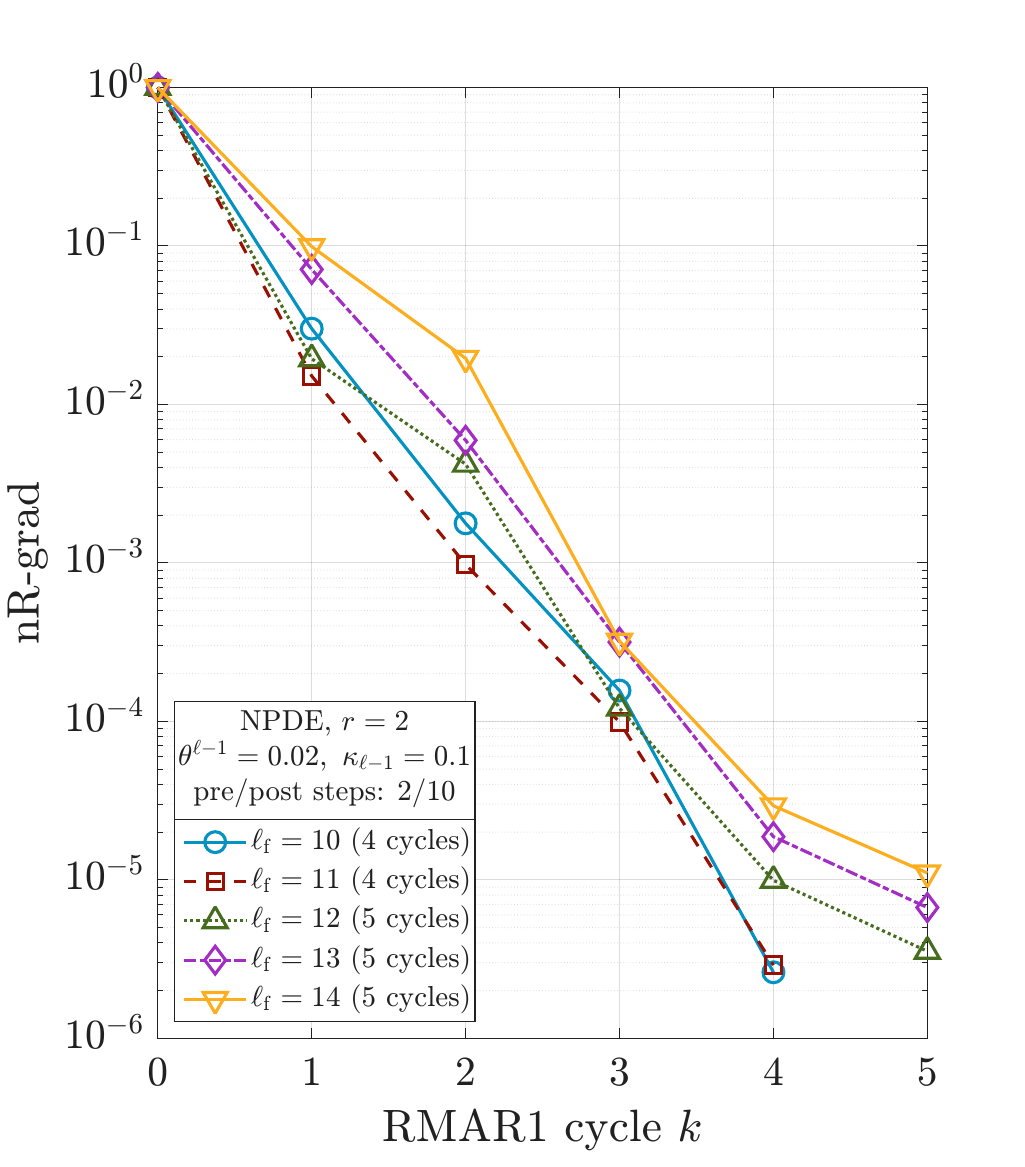}
     \caption{Mesh study for NPDE with rank $r = 2$.}
     \label{fig:rmar1_mesh_study_NPDE_K2}
  \end{minipage}\hfill
  \begin{minipage}[t]{0.485\textwidth}
     \centering
     \includegraphics[width=\linewidth,height=0.50\textheight,keepaspectratio]{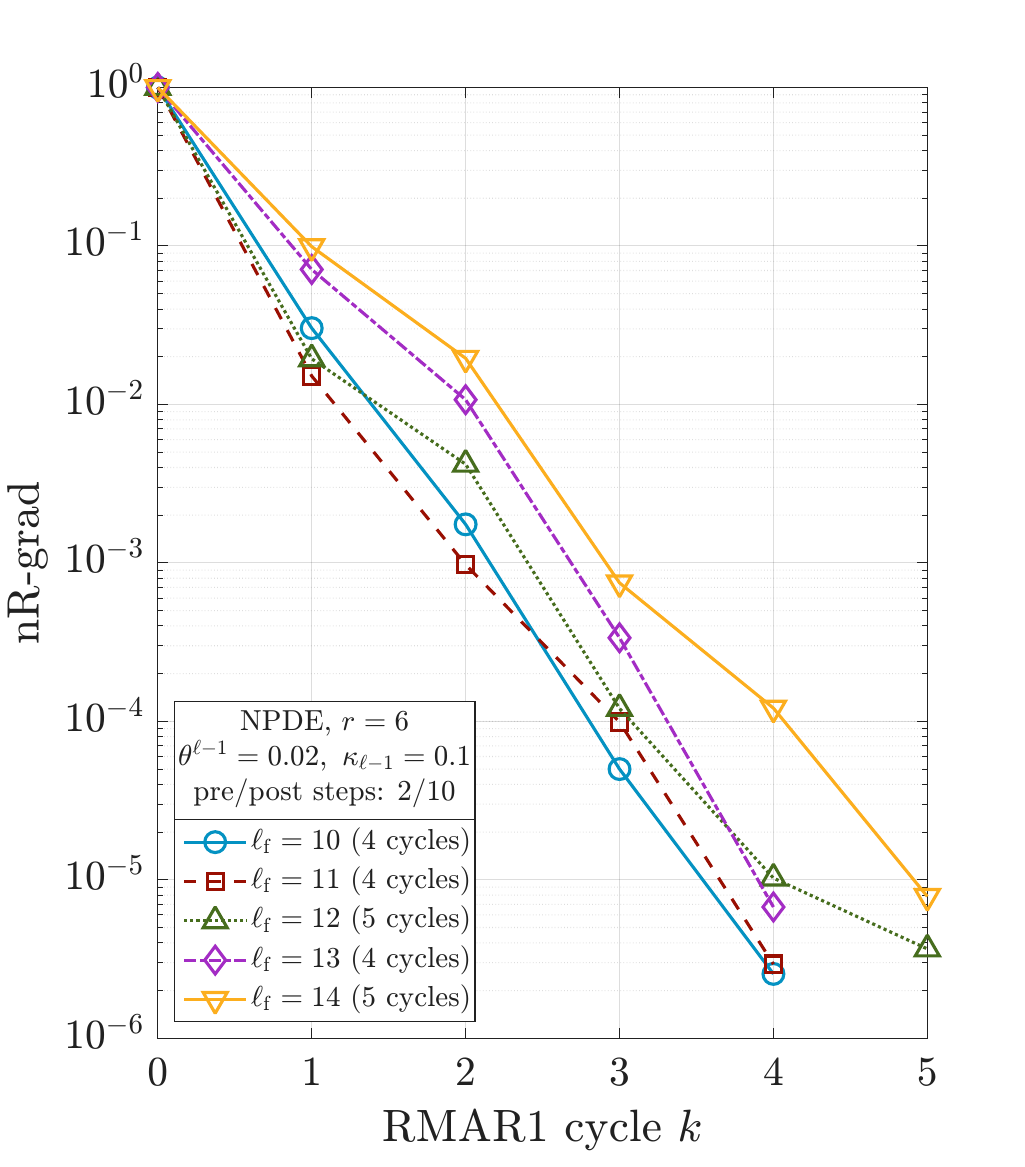}
     \caption{Mesh study for NPDE with rank $r = 6$.}
     \label{fig:rmar1_mesh_study_NPDE_K6}
  \end{minipage}
\end{figure}

\Cref{fig:rmar1_mesh_study_LYAP_K2,fig:rmar1_mesh_study_LYAP_K6,fig:rmar1_mesh_study_NPDE_K2,fig:rmar1_mesh_study_NPDE_K6} show the convergence histories of \algname{} for the LYAP and NPDE test problems, with ranks $r = 2$ and $r = 6$ and finest levels $\ell_{\mathrm{f}} = 10, \ldots, 14$. Each curve plots, against the outer-cycle index $k$, the normalized Riemannian gradient norm nR-grad, i.e., the Riemannian gradient norm divided by its value at the initial finest-level iterate $X_{0}$:
\[
    \text{nR-grad} \coloneqq \frac{\text{R-grad}_{k}}{\text{R-grad}_{0}}.
\]

Each cycle allows up to two accepted presmoothing steps and ten accepted postsmoothing steps in both problems. For each problem and rank, we keep all algorithmic settings fixed across mesh resolutions.

All twenty runs reach the prescribed tolerance. For LYAP (\Cref{fig:rmar1_mesh_study_LYAP_K2,fig:rmar1_mesh_study_LYAP_K6}), this takes 4--5 outer cycles at rank $r = 2$ and 5 at rank $r = 6$. For NPDE (\Cref{fig:rmar1_mesh_study_NPDE_K2,fig:rmar1_mesh_study_NPDE_K6}), it takes 4--5 outer cycles at both ranks.

With this initialization and stopping tolerance, the outer-cycle count varies by at most one across the tested meshes for each problem and rank. Clearly, the actual computational work per cycle depends on the mesh resolution.

\newpage


The four mesh studies below show R-grad for LYAP and NPDE at ranks $r = 2$ and $r = 6$ on finest levels 10--14. The gradient norm decreases rapidly at first and then levels off, typically above $10^{-14}$. This stagnation is consistent with the effects of round-off in MATLAB double-precision arithmetic, for which the machine epsilon is approximately $2.2204 \times 10^{-16}$. The dashed segments show the a diagnostic Taylor continuation after the coarse stopping test could no longer be certified.

\begin{figure}[H]
  \centering
  \begin{minipage}[t]{0.485\textwidth}
     \centering
     \includegraphics[width=\linewidth,height=0.34\textheight,keepaspectratio]{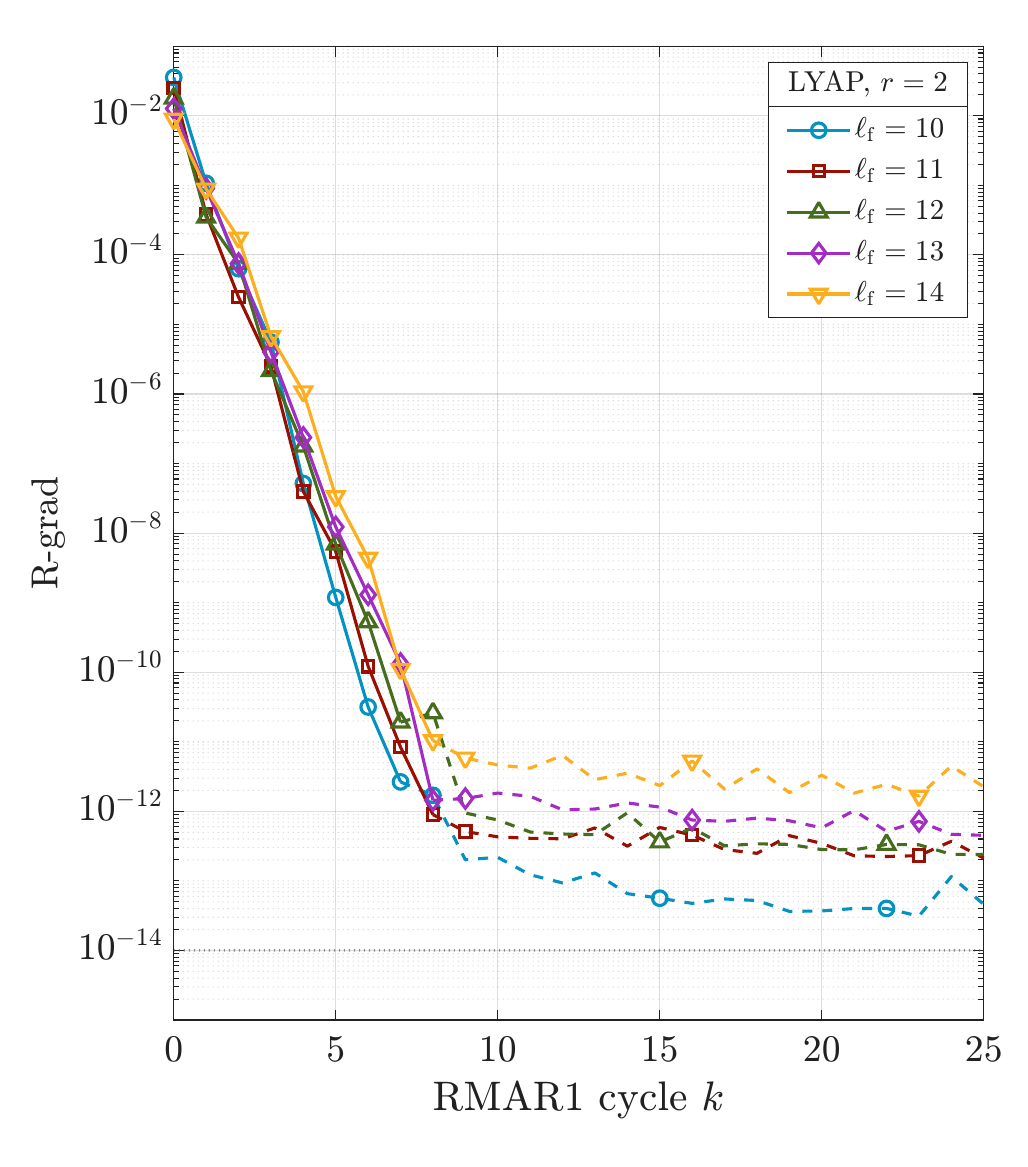}
     \caption{Mesh study for LYAP with rank $r = 2$.}
     \label{fig:rmar1_mesh_study_LYAP_K2_with_plateau}
  \end{minipage}\hfill
  \begin{minipage}[t]{0.485\textwidth}
     \centering
     \includegraphics[width=\linewidth,height=0.34\textheight,keepaspectratio]{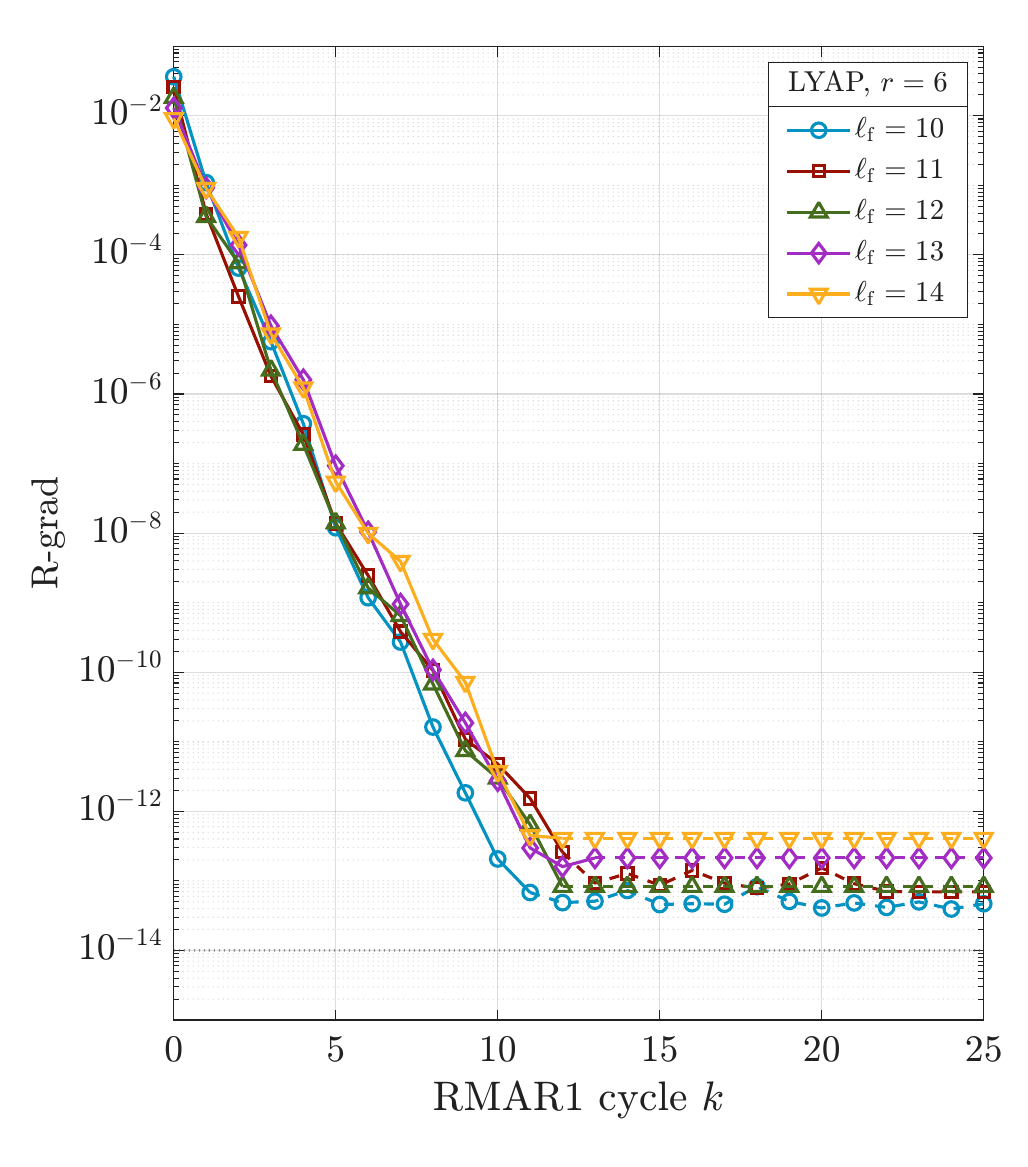}
     \caption{Mesh study for LYAP with rank $r = 6$.}
     \label{fig:rmar1_mesh_study_LYAP_K6_with_plateau}
  \end{minipage}
\end{figure}

\begin{figure}[H]
  \centering
  \begin{minipage}[t]{0.485\textwidth}
     \centering
     \includegraphics[width=\linewidth,height=0.34\textheight,keepaspectratio]{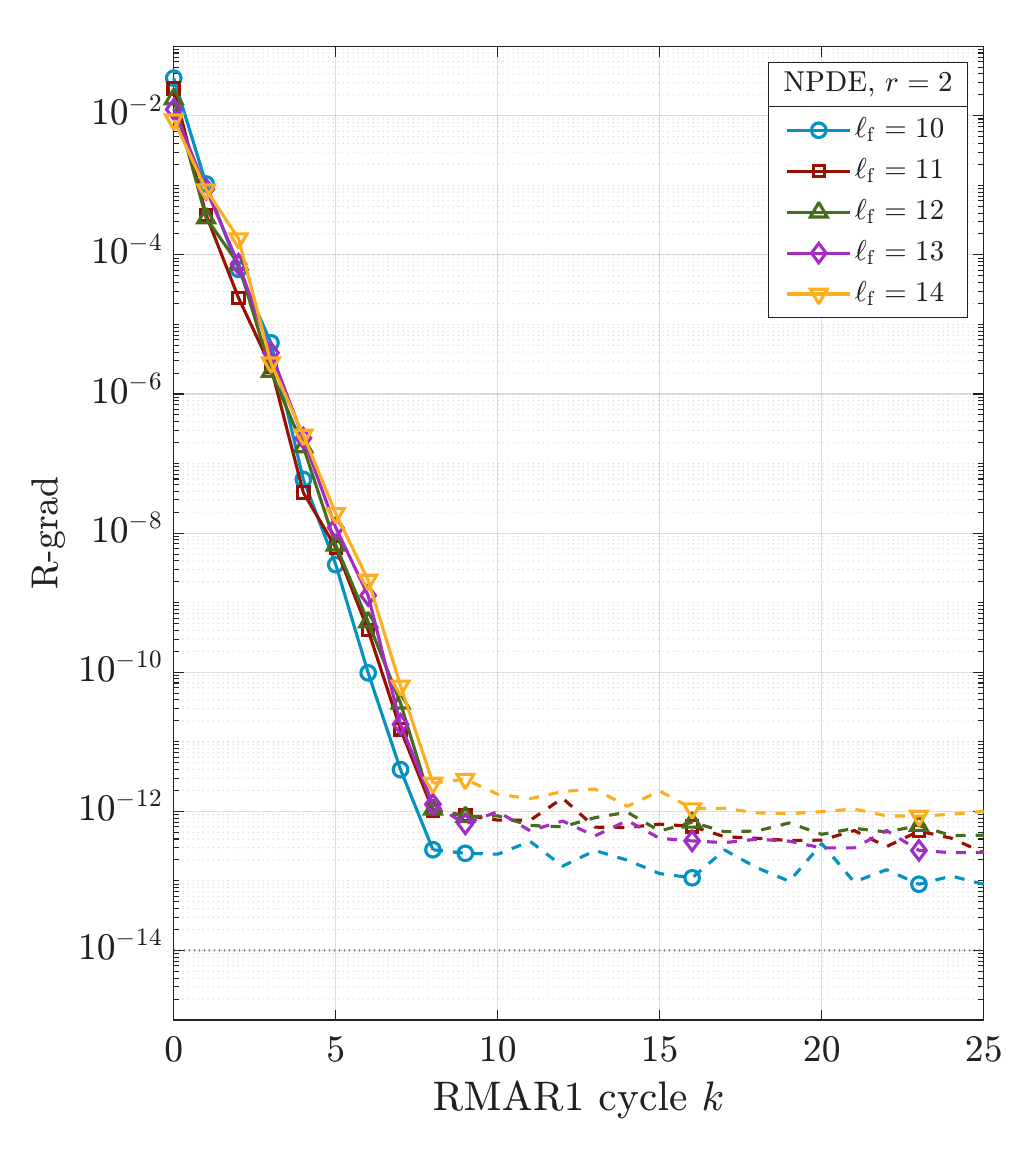}
     \caption{Mesh study for NPDE with rank $r = 2$.}
     \label{fig:rmar1_mesh_study_NPDE_K2_with_plateau}
  \end{minipage}\hfill
  \begin{minipage}[t]{0.485\textwidth}
     \centering
     \includegraphics[width=\linewidth,height=0.34\textheight,keepaspectratio]{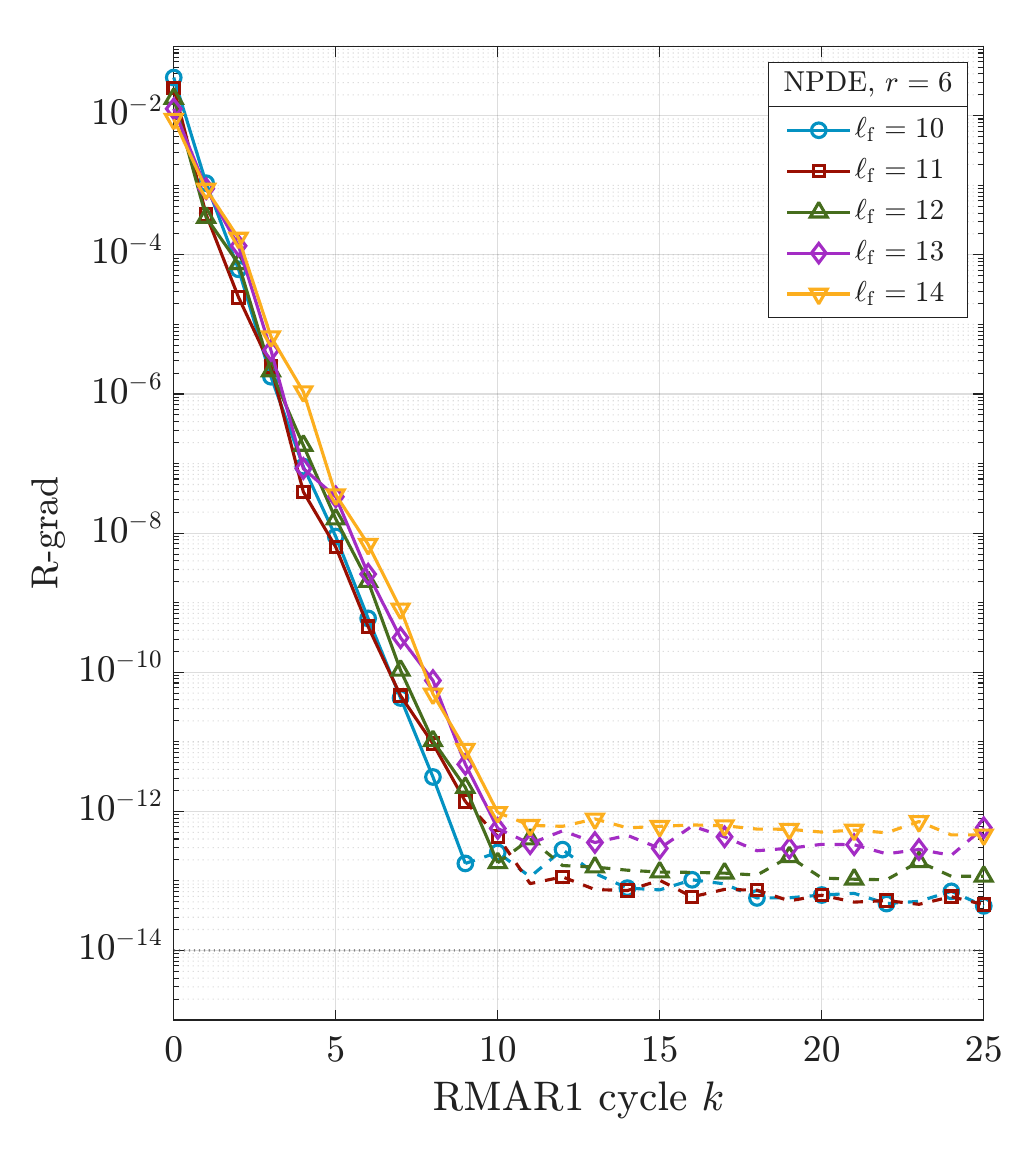}
     \caption{Mesh study for NPDE with rank $r = 6$.}
     \label{fig:rmar1_mesh_study_NPDE_K6_with_plateau}
  \end{minipage}
\end{figure}

These four plots again illustrate the nearly mesh-independent behavior of RMAR1 before stagnation. For each problem and rank, R-grad decreases at a similar rate across the tested meshes as the finest grid grows from $1025 \times 1025$ points at $\ell_{\mathrm{f}} = 10$ to $16385 \times 16385$ points at $\ell_{\mathrm{f}} = 14$.

\subsection{Comparisons with RAR1} \label{subsec:comparisons_rar1}
\begin{figure}[H]
   \centering
   \begin{minipage}[t]{0.49\textwidth}
      \centering
      \includegraphics[width=\linewidth,height=0.375\textheight,keepaspectratio]{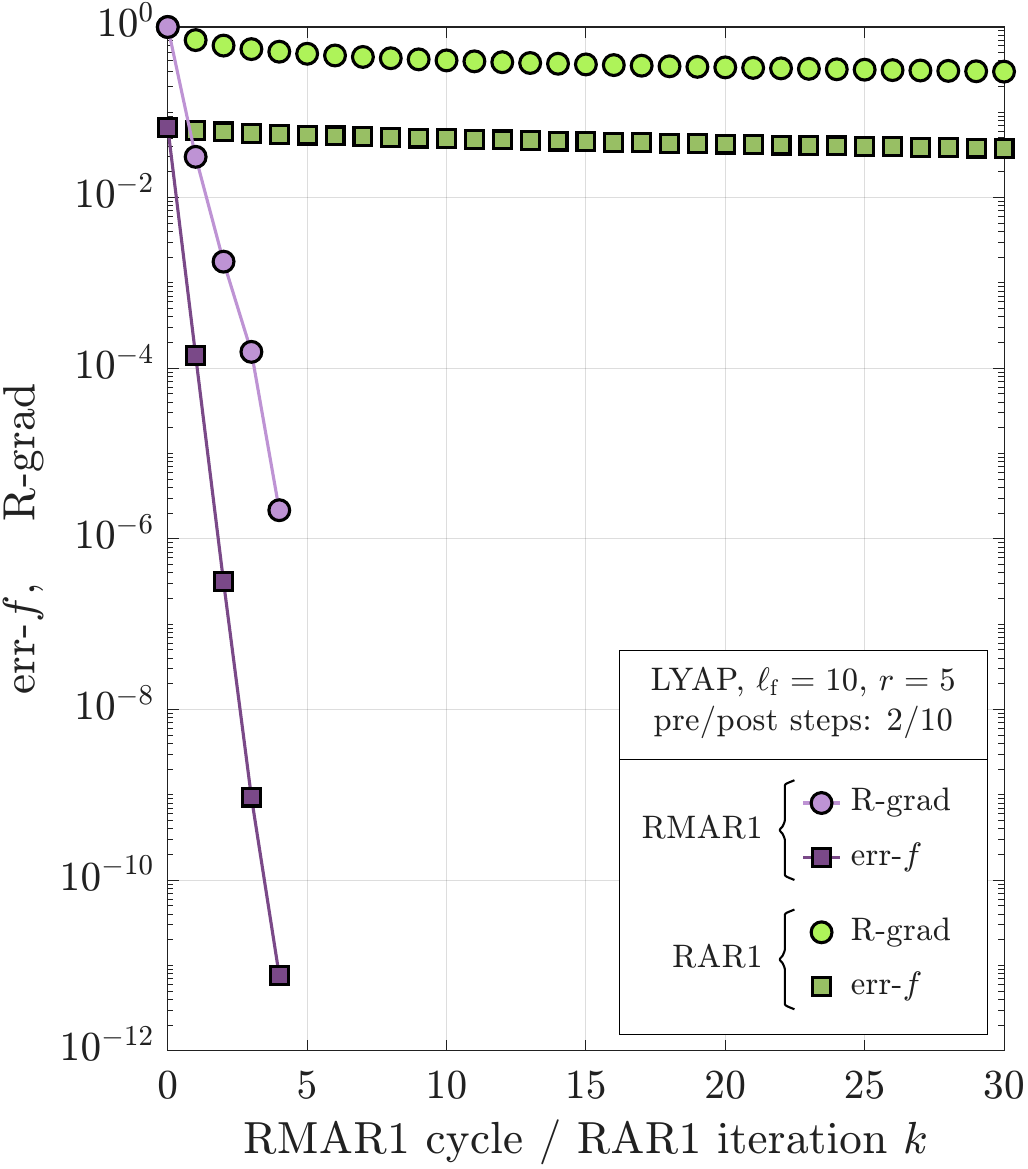}
      \caption{Convergence of RMAR1 and RAR1 versus iteration count for LYAP.} \label{fig:Comparison_RAR1_RMAR1_LYAP}
   \end{minipage}\hfill
   \begin{minipage}[t]{0.49\textwidth}
      \centering
      \includegraphics[width=\linewidth,height=0.375\textheight,keepaspectratio]{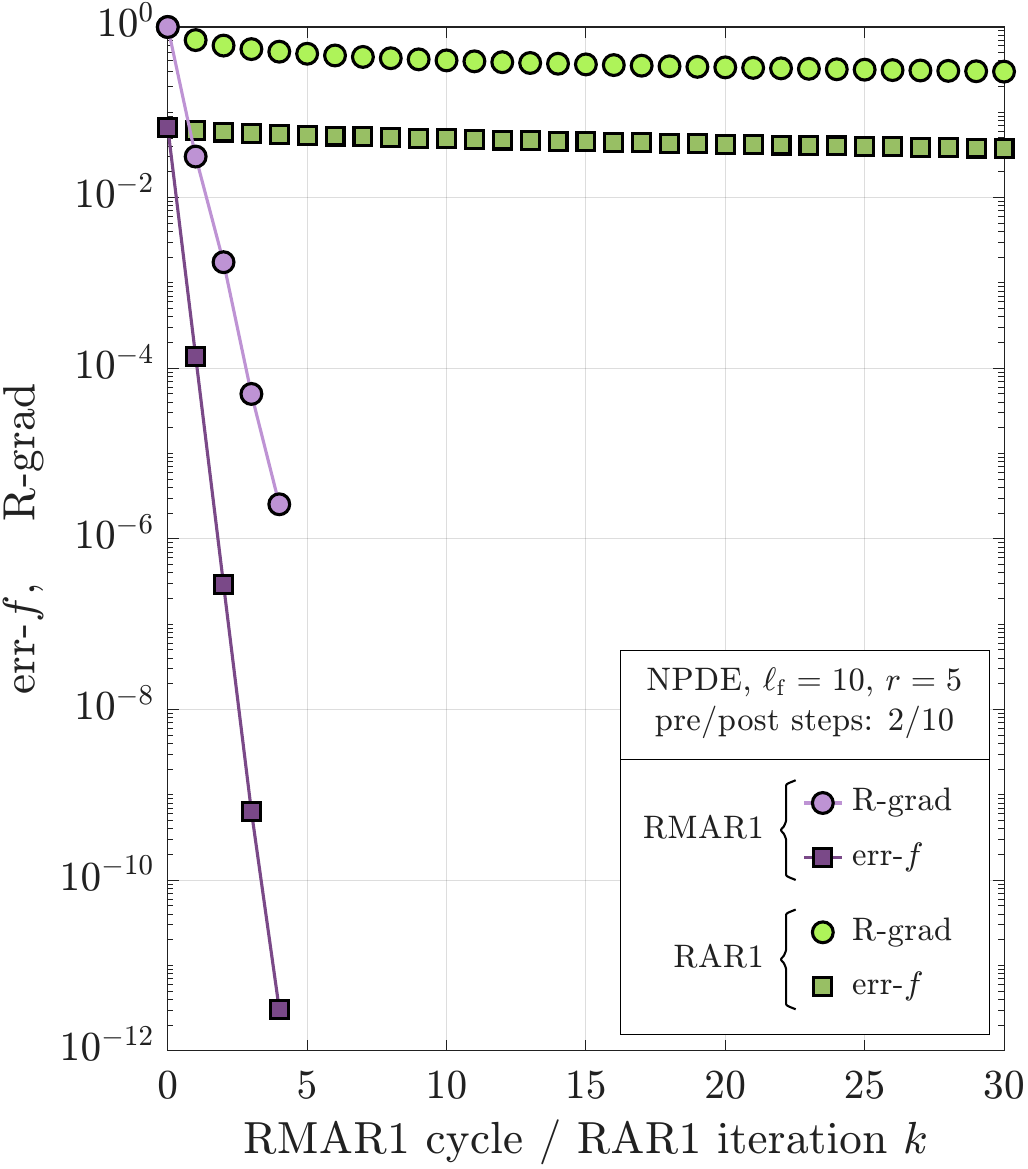}
      \caption{Convergence of RMAR1 and RAR1 versus iteration count for NPDE.} \label{fig:Comparison_RAR1_RMAR1_NPDE}
   \end{minipage}
\end{figure}

\Cref{fig:Comparison_RAR1_RMAR1_LYAP,fig:Comparison_RAR1_RMAR1_NPDE} compare \algname{} againts the plain one-level RAR1 for LYAP and NPDE, respectively. For \algname{}, we use the grid levels, rank, smoothing budgets, and coarse-solve strategy described for these problems in \Cref{subsec:convg_behavior}, together with the default regularization parameters in \Cref{subsec:setup_and_quantities}. For each problem, both methods start from the same point on the finest grid. The figures report err-$f$ and R-grad, with the same reference objective value and initial gradient norm used for both methods.

\algname{} reaches the absolute Riemannian gradient norm tolerance of $10^{-7}$ on the finest grid in four cycles for both problems. RAR1 reaches the same tolerance after $119\,899$ accepted steps for LYAP and $101\,120$ for NPDE. The horizontal axes count complete multilevel cycles for \algname{} and accepted steps for RAR1, so the two counts represent different units of work. For each problem, the iteration plot shows the first 30 accepted RAR1 steps. Within the displayed ranges, \algname{} reduces both err-$f$ and R-grad by several orders of magnitude, while RAR1 shows only a modest decrease in err-$f$.

\begin{figure}[H]
   \centering
   \begin{minipage}[t]{0.49\textwidth}
      \centering
      \includegraphics[width=\linewidth,height=0.375\textheight,keepaspectratio]{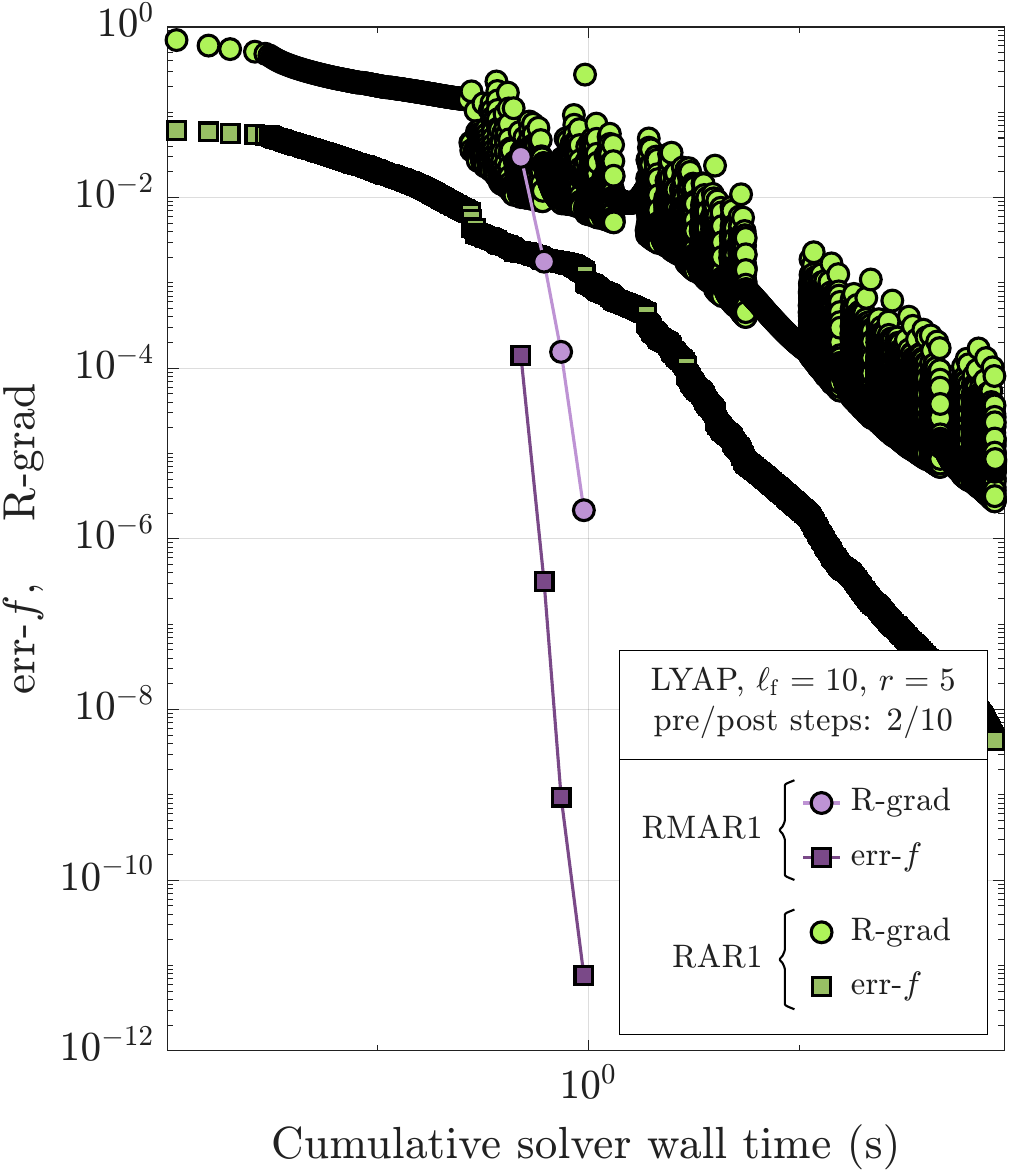}
      \caption{Convergence of RMAR1 and RAR1 versus cumulative solver wall time for LYAP.} \label{fig:Comparison_RAR1_RMAR1_LYAP_runtime}
   \end{minipage}\hfill
   \begin{minipage}[t]{0.49\textwidth}
      \centering
      \includegraphics[width=\linewidth,height=0.375\textheight,keepaspectratio]{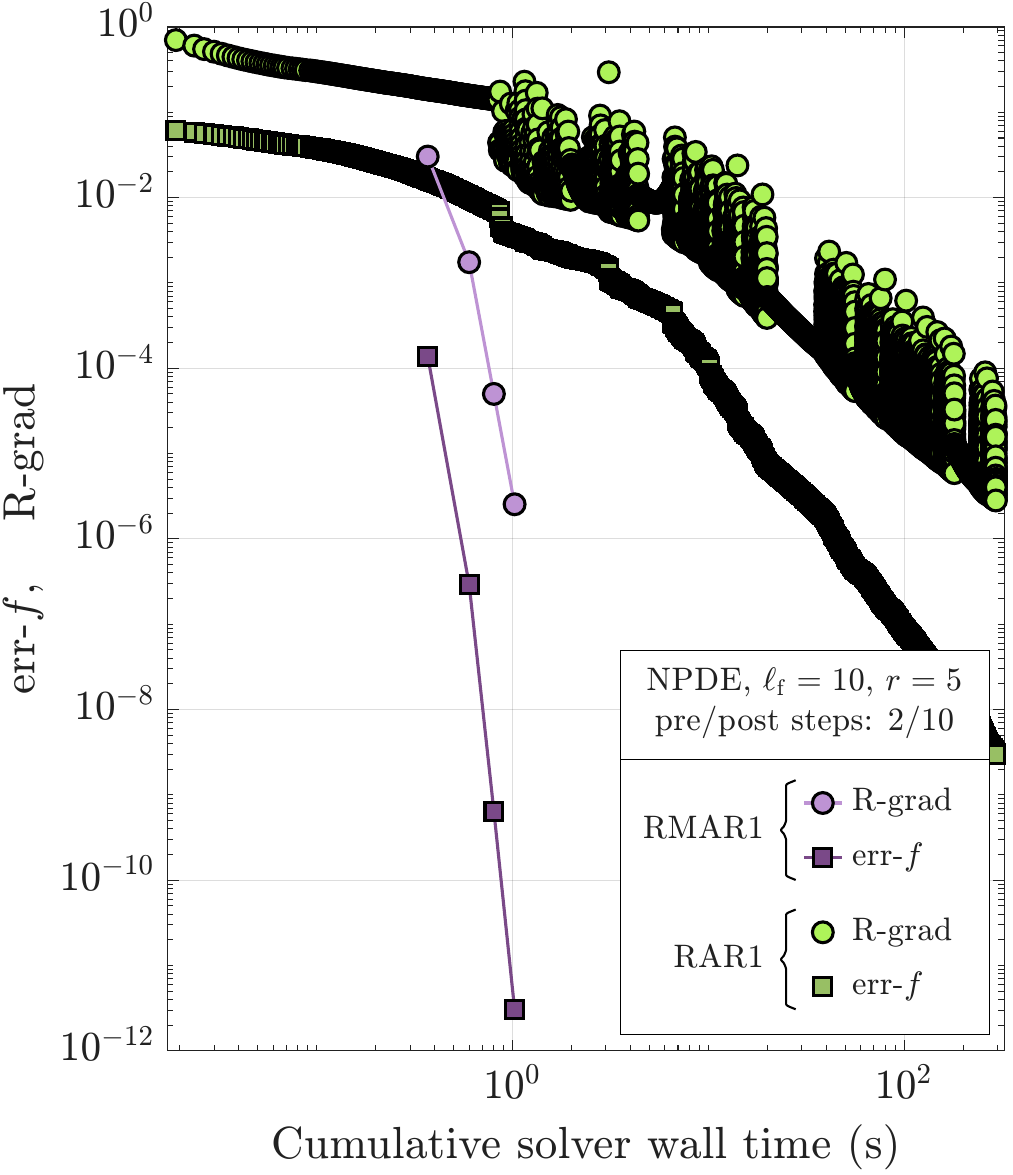}
      \caption{Convergence of RMAR1 and RAR1 versus cumulative solver wall time for NPDE.} \label{fig:Comparison_RAR1_RMAR1_NPDE_runtime}
   \end{minipage}
\end{figure}

\Cref{fig:Comparison_RAR1_RMAR1_LYAP_runtime,fig:Comparison_RAR1_RMAR1_NPDE_runtime} show the same quantities against cumulative solver wall time on a logarithmic horizontal axis. The RAR1 traces show all accepted steps, and the elapsed time includes rejected trials. The timings cover the solver calls and exclude problem setup, reference handling, and output generation. In these runs, \algname{} reaches the tolerance in $0.949\,\mathrm{s}$ for LYAP and $1.025\,\mathrm{s}$ for NPDE. The corresponding RAR1 times are $84.715\,\mathrm{s}$ and $294.417\,\mathrm{s}$.

\subsection{Comparisons with RMGLS} \label{subsec:comparisons_rmgls}
In this section, we compare \algname{} and RMGLS \cite{sutti2021riemannian} on the test problems, using a fixed rank of $r = 5$. The finest grid level ranges from $\ell_{\mathrm{f}} = 10$ to $16$, while the coarsest level is fixed at $\ell_{\mathrm{c}} = 4$. Both methods use the same initial iterate, obtained as described in \Cref{subsec:setup_and_quantities}, and use a maximum of 100 outer cycles.

RMGLS uses budgets of two presmoothing and ten postsmoothing steepest-descent steps, along with a trust-region solver on the Euclidean space of coarse-grid matrices that allows at most 40 iterations per coarsest-grid call, using the Riemannian trust-region method of \cite{absil2007trust} as implemented in Manopt \cite{manopt}. Thus, the comparison uses a common initialization and finest-grid stopping tolerance, with the same smoothing budgets and different coarsest-grid solvers.

For each problem and finest level, we report the mean and median solver times over six measured runs per method, including stopping tests but excluding problem setup, coarse initialization, and result serialization.
The second column of \Cref{tab:runtime_comparisons} reports the number of entries $(2^{\ell_{\mathrm{f}}} + 1)^2$ in the corresponding finest-grid dense matrix. All 84 measured runs for each problem reached the prescribed convergence tolerance. The last column of \Cref{tab:runtime_comparisons} reports the final R-grad, denoted R-grad$_\ast$, corresponding to the final iterate of the last measured runs for each method and finest-grid level.
All measured runs reached the tolerance. Bold text marks the lower time for each statistic at a given finest level.


\Cref{tab:runtime_comparisons} shows that \algname{} has lower mean and median solver times than RMGLS at every tested finest-grid level for LYAP and NPDE. At $\ell_{\mathrm{f}} = 16$, the ratio of the RMGLS median time to the \algname{} median time is approximately $4.17$ for LYAP and $1.62$ for NPDE. The size of the speedup varies with the problem and mesh resolution.

As a final remark, we wish to emphasize that for $\ell_{\mathrm{f}} = 16$, \algname{} and RMGLS are addressing very large-scale problems with matrices of size $65\,537 \times 65\,537$, i.e., containing $(2^{16} + 1)^{2} = 4\,295\,098\,369$ entries. Storing such a matrix in dense double-precision format requires about $34.36\ \mathrm{GB}$. In contrast, for rank $r = 5$, storing the two low-rank factors and the singular values requires only $(2 \times 65\,537 \times 5 + 5) \times 8 / 10^{6} \approx 5.24\ \mathrm{MB}$, substantially reducing the storage required to represent the matrix.

\begin{table}[htbp]
    \providecommand{\algname}{RMAR1}
    \caption{Runtime comparison of \algname{} and RMGLS for the LYAP, NPDE, and SACE problems. The column $\text{R-grad}_{\ast}$ reports the final Riemannian gradient norm.}
    \label{tab:runtime_comparisons}
    \begin{center}
        \setlength{\tabcolsep}{4pt}
        \begin{tabular}{crcccccc}
            \toprule
            & & \multicolumn{3}{c}{\algname{}} & \multicolumn{3}{c}{RMGLS} \\
            \cmidrule(lr){3-5}\cmidrule(lr){6-8}
            $\ell_{\mathrm{f}}$ & $(2^{\ell_{\mathrm{f}}} + 1)^{2}$
            & \shortstack{time\\(mean)} & \shortstack{time\\(median)} & $\text{R-grad}_{\ast}$
            & \shortstack{time\\(mean)} & \shortstack{time\\(median)} & $\text{R-grad}_{\ast}$ \\
            \midrule
            \multicolumn{8}{c}{\textbf{LYAP}} \\
            \midrule
            10 &      1\,050\,625 &  \textbf{0.39} &  \textbf{0.39} & $7.890 \times 10^{-8}$ &   0.79 &  0.79 & $9.951 \times 10^{-8}$ \\
            11 &      4\,198\,401 &  \textbf{0.60} &  \textbf{0.61} & $6.711 \times 10^{-8}$ &   1.12 &  1.12 & $5.197 \times 10^{-8}$ \\
            12 &     16\,785\,409 &  \textbf{0.84} &  \textbf{0.84} & $6.746 \times 10^{-8}$ &   1.78 &  1.78 & $6.594 \times 10^{-8}$ \\
            13 &     67\,125\,249 &  \textbf{1.39} &  \textbf{1.39} & $9.239 \times 10^{-8}$ &   2.77 &  2.77 & $8.440 \times 10^{-8}$ \\
            14 &    268\,468\,225 &  \textbf{2.63} &  \textbf{2.62} & $9.975 \times 10^{-8}$ &   5.20 &  5.30 & $7.195 \times 10^{-8}$ \\
            15 & 1\,073\,807\,361 &  \textbf{7.34} &  \textbf{7.27} & $6.133 \times 10^{-8}$ &  19.15 & 19.59 & $6.307 \times 10^{-8}$ \\
            16 & 4\,295\,098\,369 & \textbf{14.54} & \textbf{14.54} & $8.790 \times 10^{-8}$ &  60.67 & 60.69 & $8.627 \times 10^{-8}$ \\
            \midrule
            \multicolumn{8}{c}{\textbf{NPDE}} \\
            \midrule
            10 &      1\,050\,625 &  \textbf{1.21} &  \textbf{1.19} & $9.021 \times 10^{-8}$ &   2.48 &  2.47 & $8.515 \times 10^{-8}$ \\
            11 &      4\,198\,401 &  \textbf{1.46} &  \textbf{1.42} & $9.783 \times 10^{-8}$ &   2.27 &  2.22 & $9.825 \times 10^{-8}$ \\
            12 &     16\,785\,409 &  \textbf{2.83} &  \textbf{2.83} & $6.526 \times 10^{-8}$ &   5.02 &  5.02 & $5.550 \times 10^{-8}$ \\
            13 &     67\,125\,249 &  \textbf{5.35} &  \textbf{5.34} & $8.491 \times 10^{-8}$ &   7.81 &  7.40 & $8.792 \times 10^{-8}$ \\
            14 &    268\,468\,225 & \textbf{12.47} & \textbf{12.34} & $9.375 \times 10^{-8}$ &  15.23 & 15.27 & $6.774 \times 10^{-8}$ \\
            15 & 1\,073\,807\,361 & \textbf{30.54} & \textbf{30.08} & $7.263 \times 10^{-8}$ &  36.11 & 36.07 & $8.643 \times 10^{-8}$ \\
            16 & 4\,295\,098\,369 & \textbf{71.23} & \textbf{60.62} & $4.942 \times 10^{-8}$ & 102.93 & 98.36 & $9.608 \times 10^{-8}$ \\
            \midrule
            \multicolumn{8}{c}{\textbf{SACE}} \\
            \midrule
            10 &      1\,050\,625 &  \textbf{0.81} &  \textbf{0.80} & $9.021 \times 10^{-8}$ &   1.56 &  1.56 & $8.515 \times 10^{-8}$ \\
            11 &      4\,198\,401 &  \textbf{1.32} &  \textbf{1.33} & $9.783 \times 10^{-8}$ &   2.04 &  2.04 & $9.825 \times 10^{-8}$ \\
            12 &     16\,785\,409 &  \textbf{2.66} &  \textbf{2.67} & $6.526 \times 10^{-8}$ &   4.54 &  4.51 & $5.550 \times 10^{-8}$ \\
            13 &     67\,125\,249 &  \textbf{4.80} &  \textbf{4.80} & $8.491 \times 10^{-8}$ &   6.49 &  6.49 & $8.792 \times 10^{-8}$ \\
            14 &    268\,468\,225 & \textbf{11.49} & \textbf{11.49} & $9.375 \times 10^{-8}$ &  12.91 & 12.76 & $6.774 \times 10^{-8}$ \\
            15 & 1\,073\,807\,361 & \textbf{27.78} & \textbf{27.87} & $7.263 \times 10^{-8}$ &  32.85 & 32.51 & $8.643 \times 10^{-8}$ \\
            16 & 4\,295\,098\,369 & \textbf{58.86} & \textbf{58.30} & $4.942 \times 10^{-8}$ &  94.38 & 93.85 & $9.608 \times 10^{-8}$ \\
            \bottomrule
        \end{tabular}
    \end{center}
\end{table}

\subsection{Rank-adaptive \algname{}} \label{subsec:rank_adaptivity} 
As discussed at the end of \Cref{subsec:test_problems}, the rank needed to achieve a prescribed accuracy is typically not known a priori. We address this scenario by proposing a rank-adaptive version of \algname{} able to increase the rank during the solve. We measure the accuracy of the current iterate $X_k$ by its relative Frobenius error
\[
    \text{err-}X \coloneqq \frac{\|X_k - X_{\mathrm{ref}}\|_{\mathrm{F}}}{\|X_{\mathrm{ref}}\|_{\mathrm{F}}},
\]
where $X_{\mathrm{ref}}$ is a highly accurate solution of the same discretized problem computed without a rank constraint\footnote{Compare, in contrast, with the experiments in \Cref{subsec:convg_behavior}, where the reference for err-$f$ in the fixed-rank convergence study is the objective value at a rank-constrained solution.}. At each rank $r$, we report two relative Frobenius errors, both measured against $X_{\mathrm{ref}}$ and divided by $\|X_{\mathrm{ref}}\|_{\mathrm{F}}$. The final relative error is err-$X$ evaluated at the last iterate of the rank-$r$ stage. The best rank-$r$ relative error is
\[
    \frac{\|X_{\mathrm{ref}}^{(r)} - X_{\mathrm{ref}}\|_{\mathrm{F}}}{\|X_{\mathrm{ref}}\|_{\mathrm{F}}},
\]
where $X_{\mathrm{ref}}^{(r)}$ is the rank-$r$ truncated SVD of $X_{\mathrm{ref}}$. This is the smallest relative error achievable by a rank-$r$ approximation to $X_{\mathrm{ref}}$. The reference solution is used only to evaluate these errors, and it is not employed in the rank-adaptive procedure.

For both LYAP and NPDE, we set $\ell_{\mathrm{f}} = 10$ and $\ell_{\mathrm{c}} = 4$, start at rank $r_0 = 5$. Each \algname{} cycle uses up to three presmoothing steps and one postsmoothing step, with at most ten RAR1 attempts per coarsest-grid visit.

On the finest grid, we monitor the Frobenius norms of the ambient gradient $\nabla\bar{f}(X)$ and the Riemannian gradient $\operatorname{grad} f(X)$. A stage is ready for a rank decision when R-grad is less than or equal to $10^{-12}$ or the stagnation test is satisfied.  Let $Z_k$ denote the iterate after $k$ cycles in the current stage, with $Z_0$ its starting point. From cycle $k = 8$ onward, a cycle qualifies for the stagnation test if the Riemannian gradient norm is less than or equal to $10^{-8}$ and the relative iterate change $\|Z_k - Z_{k - 1}\|_{\mathrm{F}}/\|Z_{k - 1}\|_{\mathrm{F}}$ is less than or equal to a prescribed fraction of the positive running maximum of that change through cycle $k$. When the current ambient gradient norm exceeds $\tau_{\mathrm{amb}} = 10^{-10}$, the test requires seven consecutive qualifying cycles with a fraction of $1\%$. For LYAP, this branch also requires the ratio of Riemannian to ambient gradient norms to be less than or equal to $0.2$ at the current iterate. When the current ambient gradient norm is less than or equal to $\tau_{\mathrm{amb}}$, the test requires three consecutive qualifying cycles with a fraction of $0.2\%$.

At the end of a stage, we stop if the ambient gradient norm is less than or equal to $\tau_{\mathrm{amb}}$; otherwise, we increase the rank by one. To warm-start the higher-rank stage, we use the rank-one normal correction of Gao and Absil~\cite[Alg.~3]{GaoAbsil:2022}. At the current rank-$r$ stage, let the final iterate be $X = U\Sigma V^{\top}$, and let $(\sigma_{\perp},u_{\perp},v_{\perp})$ be a leading singular triplet of the negative normal gradient $G_{\perp} = - P_{U}^{\perp} \bigl( \nabla\bar{f}(X) \bigr) P_{V}^{\perp}$. We set
\begin{equation}\label{eq:rank_one_correction}
    X_{\mathrm{new}} = X + \alpha \, \sigma_{\perp} u_{\perp}v_{\perp}^{\top},
\end{equation}
where $\alpha > 0$ is chosen by Armijo backtracking on the finest-grid objective. For $\sigma_{\perp} > 0$, the singular vectors $u_{\perp}$ and $v_{\perp}$ are orthogonal to the columns of $U$ and $V$, respectively, so the correction \eqref{eq:rank_one_correction} raises the rank to $r + 1$.

\begin{figure}[H]
   \centering
   \begin{minipage}[t]{0.49\textwidth}
      \centering
      \includegraphics[width=\linewidth,height=0.375\textheight,keepaspectratio]{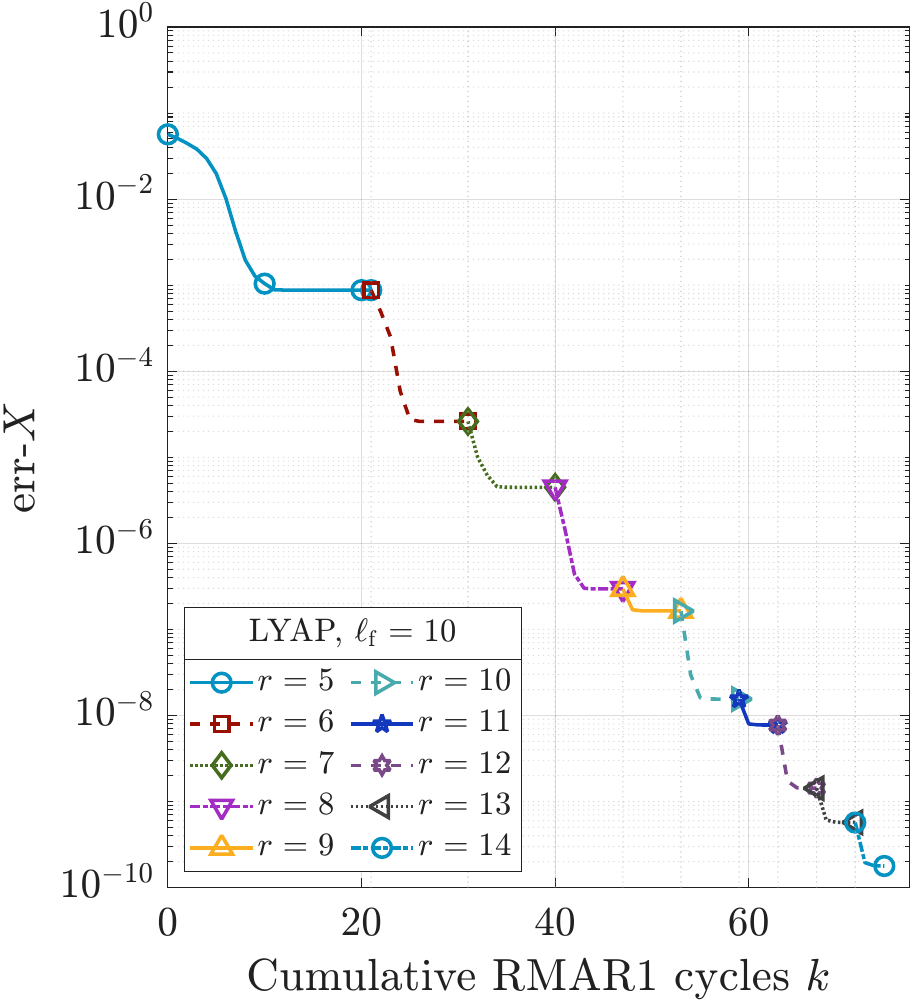}
      \caption{Rank-adaptive \algname{} for LYAP.} \label{fig:RMAR1_ErrX_LYAP_lev10_autoRankStudy_arcCoarseStartRank1Ambient1e10}
   \end{minipage}\hfill
   \begin{minipage}[t]{0.49\textwidth}
      \centering
      \includegraphics[width=\linewidth,height=0.375\textheight,keepaspectratio]{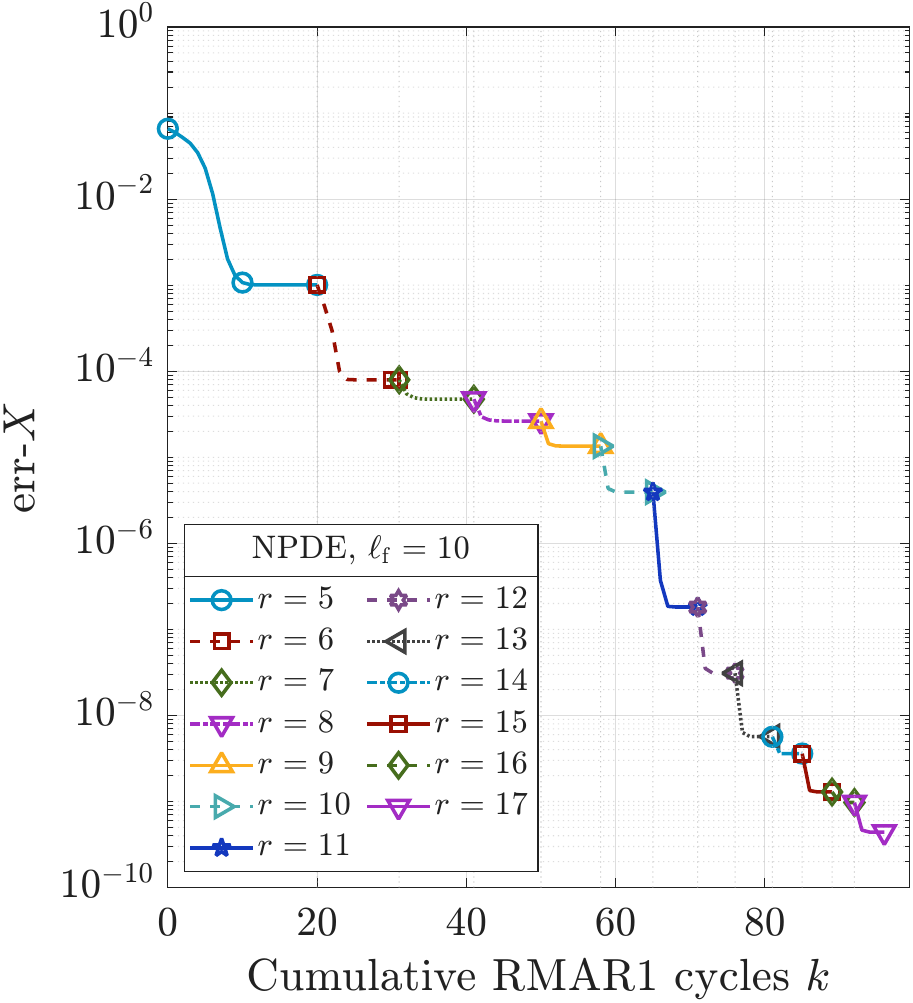}
      \caption{Rank-adaptive \algname{} for NPDE.} \label{fig:RMAR1_ErrX_NPDE_lev10_autoRankStudy_arcCoarseStartRank1Ambient1e10RankCap20}
   \end{minipage}
\end{figure}

\Cref{fig:RMAR1_ErrX_LYAP_lev10_autoRankStudy_arcCoarseStartRank1Ambient1e10,fig:RMAR1_ErrX_NPDE_lev10_autoRankStudy_arcCoarseStartRank1Ambient1e10RankCap20} show $\text{err-}X$ versus cumulative \algname{} cycles. 
\Cref{tab:rank_adaptivity_LYAP,tab:rank_adaptivity_NPDE} report the cycle count and errors for each rank. LYAP satisfies the ambient-gradient tolerance $\tau_{\mathrm{amb}}$ at rank $14$ after $74$ cycles, with an ambient gradient norm of $6.33 \times 10^{-11}$. Its final err-$X$ is $1.78 \times 10^{-10}$, compared with a best rank-$14$ error of $1.38 \times 10^{-10}$. NPDE satisfies the same tolerance at rank $17$ after $96$ cycles, with an ambient gradient norm of $9.98 \times 10^{-11}$. Its final err-$X$ is $4.42 \times 10^{-10}$, compared with a best rank-$17$ error of $4.15 \times 10^{-10}$.

\begin{table}[htbp]
   \centering
   \caption{LYAP rank-adaptivity results with unit rank increments, and $\tau_{\mathrm{amb}} = 10^{-10}$.} \label{tab:rank_adaptivity_LYAP}
   \begin{tabular}{crcc}
      \toprule
      Rank $r$ & Cycles & Final err-$X$ & Best rank-$r$ relative error \\
      \midrule
       5 & 21 & $8.75\times10^{-4}$ & $8.73\times10^{-4}$ \\
       6 & 10 & $2.61\times10^{-5}$ & $2.60\times10^{-5}$ \\
       7 & 9 & $4.48\times10^{-6}$ & $4.44\times10^{-6}$ \\
       8 & 7 & $2.95\times10^{-7}$ & $2.87\times10^{-7}$ \\
       9 & 6 & $1.64\times10^{-7}$ & $1.53\times10^{-7}$ \\
      10 & 6 & $1.55\times10^{-8}$ & $1.41\times10^{-8}$ \\
      11 & 4 & $7.74\times10^{-9}$ & $7.22\times10^{-9}$ \\
      12 & 4 & $1.43\times10^{-9}$ & $1.20\times10^{-9}$ \\
      13 & 4 & $5.72\times10^{-10}$ & $5.00\times10^{-10}$ \\
      14 & 3 & $1.78\times10^{-10}$ & $1.38\times10^{-10}$ \\
      \bottomrule
   \end{tabular}
\end{table}

\begin{table}[htbp]
   \centering
   \caption{NPDE rank-adaptivity results with unit rank increments, and $\tau_{\mathrm{amb}} = 10^{-10}$.} \label{tab:rank_adaptivity_NPDE}
   \begin{tabular}{crcc}
      \toprule
      Rank $r$ & Cycles & Final err-$X$ & Best rank-$r$ relative error \\
      \midrule
       5 & 20 & $1.01\times10^{-3}$ & $1.00\times10^{-3}$ \\
       6 & 11 & $7.95\times10^{-5}$ & $7.94\times10^{-5}$ \\
       7 & 10 & $4.74\times10^{-5}$ & $4.70\times10^{-5}$ \\
       8 & 9 & $2.63\times10^{-5}$ & $2.61\times10^{-5}$ \\
       9 & 8 & $1.35\times10^{-5}$ & $1.34\times10^{-5}$ \\
      10 & 7 & $3.94\times10^{-6}$ & $3.90\times10^{-6}$ \\
      11 & 6 & $1.82\times10^{-7}$ & $1.80\times10^{-7}$ \\
      12 & 5 & $3.10\times10^{-8}$ & $3.08\times10^{-8}$ \\
      13 & 5 & $5.68\times10^{-9}$ & $5.52\times10^{-9}$ \\
      14 & 4 & $3.61\times10^{-9}$ & $3.47\times10^{-9}$ \\
      15 & 4 & $1.29\times10^{-9}$ & $1.26\times10^{-9}$ \\
      16 & 3 & $9.74\times10^{-10}$ & $8.59\times10^{-10}$ \\
      17 & 4 & $4.42\times10^{-10}$ & $4.15\times10^{-10}$ \\
      \bottomrule
   \end{tabular}
\end{table}

\section{Conclusion and outlook} \label{sec:conclusions}
In this work, we introduced \algname{}, the first multilevel Riemannian adaptive regularization algorithm for optimization problems on manifolds endowed with a hierarchical structure. The method alternates between minimization steps of regularized first-order models on the fine manifold and steps computed on the coarse manifold. It generalizes the RAR1 method, to which it reduces when only fine steps are performed. For general objective functions, we established a worst-case iteration complexity bound and proved convergence to first-order $\varepsilon$-approximate stationary points, thereby extending the Euclidean multilevel AR$q$ theory (with $q=1$) of Calandra et al. \cite{Calandra:2021} to the Riemannian setting. We then specialized the framework to the manifold of fixed-rank matrices, following the implementation of RMGLS in \cite{sutti2021riemannian}. Numerical experiments on large-scale low-rank matrix optimization problems showed moderate mesh dependence and consistently outperformed RMGLS, the state-of-the-art multilevel Riemannian solver, while providing the convergence guarantees that RMGLS lacks.
As a final remark, we wish to stress that combining multilevel corrections with low-rank representations allows us to handle very large-scale problems that would otherwise be prohibitive on a laptop: for instance, in \Cref{subsec:comparisons_rmgls}, \algname{} solved problems corresponding to uniform square grids of side $65\,537$ grid-points, reducing storage per matrix from $34.36\,\mathrm{GB}$ in dense form to $5.24\,\mathrm{MB}$ at rank $5$.

{\footnotesize
\subsection*{\footnotesize Acknowledgments}
We would like to thank Filippo Marini and the master's students of ENS de Lyon, Alexandre Lelong and Naomi Bronstein, for their contribution to the preliminary development of this project.
We also thank Pierre-Antoine Absil for useful discussions about retractions on the manifold of fixed-rank matrices.

\subsection*{\footnotesize Ethics approval and consent to participate}
Not applicable.
\subsection*{\footnotesize Consent for publication}
Not applicable.
\subsection*{\footnotesize Funding}
{\footnotesize
The work of M.P. was partially supported by INdAM-GNCS under the project CUP E53C25002010001. The research of M.P. and M.S. was partially funded by the Italian Ministry of University and Research (MUR) through the PRIN 2022 ``MOLE: Manifold constrained Optimization and LEarning'', code: 2022ZK5ME7, MUR D.D. financing decree n. 20428 of November 6th, 2024 (CUP B53C24006410006). The work of E.R. was partially supported by the ANR-24-CE23-7039
MEPHISTO project and the Simone and Cino Del Duca Foundation. 

}
\subsection*{\footnotesize Availability of data and materials}
The data presented in this manuscript are reproducible through the implementation publicly available at \url{https://github.com/MarcoSutti/RMAR1}.

\subsection*{\footnotesize Competing interests}
The authors declare that they have no competing interests.
\subsection*{\footnotesize Authors' contributions}
All authors contributed equally to the writing of this article. All authors reviewed the manuscript.}

\subsection*{\footnotesize AI declaration}
The authors assume responsibility for all content of the present paper. The research direction, method, and approach to the theoretical analysis, including proofs, as well as the core numerical implementation of the proposed method, were developed by the authors without using any AI tool. Anthropic Claude was prompted to review the draft and check the correctness of the mathematical arguments. OpenAI Codex was used to debug and optimize MATLAB code, run benchmarks and parameter studies, format figures and tables, and align the pseudocodes with the actual implementation.

\appendix

\section{The Riemannian Adaptive Regularized Gradient method (RAR1)}\label[appendix]{app:RAR1}

Let us consider the following minimization problem on a given Riemannian manifold $\cM$ 
\[
    \underset{x \in \cM}{\min} f(x),
\]
where $ f \colon \cM \mapsto \R $ is a (sufficiently smooth) cost function.

We review the main steps and properties of the Riemannian version of the AR1 (RAR1), which can be interpreted as a first-order version of the second-order method proposed in \cite{agarwal2021adaptive}. At each iteration $k$, given the current iterate $x_{k} \in \cM$, the following Taylor model is built:
\begin{equation}\label{app:eq:taylor}
    T_{k}(s) \coloneqq f(x_{k}) + \langle \grad f(x_{k}), s \rangle_{x_{k}},
\end{equation}
and a step $s_{k}$ in the tangent space is computed by minimizing the regularized model $m_{k} \colon \T_{x_{k}}\cM \to \R$
\[
    m_{k}(s) \coloneqq  T_{k}(s)  + \frac{\lambda_{k}}{2} \| s \|_{x_{k}}^{2},
\]
where $\lambda_k$ is a positive value called the regularization parameter that determines the step length. Minimizing the model in this case indeed amounts to choosing
\[
   s_{k} = -\frac{1}{\lambda_k} \grad f(x_{k}).    
\]
Given a retraction $\operatorname{R}$ on the manifold $\cM$, the step $s_{k}$ is used to define a trial point as 
\[
    \tilde x_{k+1} = \operatorname{R}_{x_{k}}(s_{k}).
\] 
At each iteration, we must decide whether to accept or reject this trial point. This decision is based on the agreement between the actual decrease in the function and the model's predicted decrease. More precisely, at each iteration, both the decrease achieved in the model, which we call predicted reduction, $ \mathrm{pred} =  T_{k}(0) -  T_{k}(s_{k})$, and that achieved in the objective function, called actual reduction, $\mathrm{ared} = f(x_{k})-f(\tilde x_{k+1})$, are computed. The step acceptance is then based
on the ratio:
\begin{equation}\label{app:eq:ratio}
    \rho_k = \frac{\mathrm{ared}}{\mathrm{pred}} = \frac{f(x_{k})-f(\tilde x_{k+1})}{T_{k}(0) -  T_{k}(s_{k})} = \lambda_k \, \frac{f(x_{k})-f(\tilde x_{k+1})}{\| \grad f(x_{k})\|^2_{x_{k}}}.
\end{equation}
If the model is accurate, $\rho_k$ will be close to one. Then, the trial step $s_{k}$ is accepted if $\rho_k$ is larger than or equal to a chosen threshold $\eta_1\in (0,1)$, and is rejected otherwise. In the first case, the step is successful; otherwise, it is unsuccessful.
After the step acceptance, the regularization parameter is updated for the next iteration. The update is still based on the ratio \eqref{app:eq:ratio}. If the step is successful, the regularization parameter decreases; otherwise, it increases. The procedure stops when it reaches a minimizer of $f$. Usually, the stopping criterion is based on the norm of the Riemannian gradient; i.e., given an absolute accuracy level $\varepsilon> 0$, the iterations stop as soon as $\|\grad f(x_{k})\|_{x_k} < \varepsilon$. \Cref{app:alg:RAR1} sketches the full procedure.

\subsection{Complexity analysis of RAR1}\label{sec:compRAR1}
 The complexity analysis of \Cref{app:alg:RAR1} is carried out under the following standard regularity assumptions.
\begin{enumerate}
    \item[\textbf{A1}] The cost function $f$ is continuously differentiable.
    \item[\textbf{A2}]  The cost function $f$ is bounded below, that is there exists $ f_{\text{low}}$ such that $f(x) \geq f_{\text{low}}$ for all $x \in \cM$.
    
    \item[\textbf{A3}]  For each iterate $x_{k}$ and step $s_{k}$ generated by~\Cref{app:alg:RAR1}, there exists a constant $L$ such that the pullback 
    $\hat{f}_k = f \circ \operatorname{R}_{x_{k}}$ satisfies 
    \[
        \hat{f}_k(s_{k}) - \left [ f(x_{k}) + \innerproduct{s_{k}}{\nabla  \hat{f}(0_{x_{k}})}_{x_{k}}\right] \leq \frac{L}{2} \|s_{k}\|_{x_{k}}^2.
    \] 
\end{enumerate}
We thus have the following property on the number of function and gradient evaluations.
\begin{theorem}
Let $\cM$ be a Riemannian manifold equipped with a retraction $\operatorname{R}$. 
Under \textbf{A1}, \textbf{A2} and \textbf{A3}, for an arbitrary $x_0 \in \cM$, let $x_0,x_1, x_2 \dots$ be the iterates produced by~\Cref{app:alg:RAR1}. For any $\varepsilon > 0$,
the RAR1 algorithm requires at most $\mathcal{O}(\varepsilon^{-2})$ function evaluations, and at most $\mathcal{O}(\varepsilon^{-2})$ gradient evaluations to compute an iterate $x_\varepsilon$ such that $x_\varepsilon$ is an $\varepsilon$-approximate first-order minimizer of $f$, i.e., $\|{\grad  f(x_\varepsilon)}\|_{x_\varepsilon}\leq \varepsilon$.
\end{theorem}
\begin{proof}
The bounds on the number of function evaluations and on the number of gradient evaluations can be easily computed following \cite[Cor. 4.8]{boumal2023introduction} and \cite[Thm. 2.4.4]{book_compl}.
\end{proof}

\begin{algorithm}
    \caption{RAR1$(x_0,\lambda_0,\varepsilon)$ The Riemannian AR1 method.}
    \label{app:alg:RAR1}
    \DontPrintSemicolon
    \KwIn{Initial point $x_0 \in \cM$, initial regularization parameter $\lambda_0 > 0$, and gradient tolerance $\varepsilon > 0$.}
    \KwOut{An iterate $x_\ast \in \cM$ satisfying $\lVert \grad f(x_\ast) \rVert_{x_\ast} \leq \varepsilon$.}
    \textbf{Parameters.} The constants $\eta_1$, $\eta_2$, $\gamma_1$, $\gamma_2$, $\gamma_3$, and $\lambda_{\mathrm{min}}$ satisfy
    \begin{algomathdisplay}
        \lambda_{\mathrm{min}} \in (0,\lambda_0], \quad 0 < \eta_1 \leq \eta_2 < 1, \quad 0 < \gamma_2 \leq \gamma_1 < 1 < \gamma_3.
    \end{algomathdisplay}
    Compute $f(x_0)$ and set $k = 0$.\;
    \While{$\lVert \grad f(x_k) \rVert_{x_k} > \varepsilon$}{
        \textbf{Initialization.} Define the model $T_k$ as in \eqref{app:eq:taylor}.\;
        \textbf{Model minimization.} Define the step $s_k = - \frac{1}{\lambda_k}\grad f(x_k) \in \T_{x_k}\cM$ and compute the trial point $\tilde{x}_{k + 1} = \operatorname{R}_{x_k}(s_k)$.\;
        \textbf{Acceptance of the trial point.} Compute
        \begin{algomathdisplay}
            \rho_k = \frac{f(x_k) - f(\tilde{x}_{k + 1})}{- \innerproduct{\grad f(x_k)}{s_k}_{x_k}}.
        \end{algomathdisplay}
        \eIf{$\rho_k \geq \eta_1$}{
            $x_{k + 1} = \tilde{x}_{k + 1}$.\;
        }{
            $x_{k + 1} = x_k$.\;
        }
        \textbf{Regularization parameter update.}\;
        \eIf{$\rho_k \geq \eta_1$}{
            \begin{algomathdisplay}
                \lambda_{k + 1} = \begin{cases}
                    \max\{\lambda_{\mathrm{min}}, \gamma_2\lambda_k\} & \text{if } \rho_k \geq \eta_2,\\
                    \max\{\lambda_{\mathrm{min}}, \gamma_1\lambda_k\} & \text{if } \rho_k < \eta_2.
                \end{cases}
            \end{algomathdisplay}
        }{
            $ \lambda_{k + 1} = \gamma_3\lambda_k $.\;
        }
        $k = k + 1$.\;
    }
    Set $x_\ast \gets x_k$.\;
    \Return{$x_\ast$}.\;
\end{algorithm}

\input{supplementary}

\bibliographystyle{aomalpha}
\bibliography{RMAR_2026}

\end{document}

%% file: supplementary.tex
%
%
%



\section{Detailed calculations}\label[appendix]{app:suppl}

In this section, we detail computations used in this work and specify the derivation of quantities used in \Cref{sec:two-level,sec:theory} for the specific choice of the manifold of fixed-rank matrices $\mathcal{M}_{r}$.

Throughout this section, $\cMh=\mathcal{M}_{r}^{h}$ and $\cMH=\mathcal{M}_{r}^{H}\subset\R^{m\times n}$ denote the fine and coarse fixed-rank matrix manifolds, equipped with their induced Frobenius metrics, and $f^{h}\colon\cMh\to\R$ and $f^{H}\colon\cMH\to\R$ are smooth objectives from the discretization hierarchy of $f$. We fix a point $X_{k}^{h}\in\cMh$ and assume that $X_{k}^{H}=\mathcal{I}_{h}^{H}X_{k}^{h}\in\cMH$. The tangent transfer operators are fixed linear maps (see also discussion in \Cref{sec:preliminaries})
\[
    \widetilde{\mathcal{I}}_{h}^{H}\colon\mathrm{T}_{X_{k}^{h}}\cMh\to\mathrm{T}_{X_{k}^{H}}\cMH, \qquad \widetilde{\mathcal{I}}_{H}^{h}\colon\mathrm{T}_{X_{k}^{H}}\cMH\to\mathrm{T}_{X_{k}^{h}}\cMh.
\]
Consider any smooth retraction $\operatorname{R}_{X_k^H} \colon \mathrm{T}_{X^{H}_{k}}\cMH \to \cMH $. Since $\operatorname{D} \! \operatorname{R}_{X_k^H}(0_{X_k^H})=\operatorname{Id}_{\mathrm{T}_{X_k^H}\mathcal{M}^{H}}$, the inverse function theorem guarantees the existence of an open neighborhood $\mathcal U_k^H$ of $0_{X_k^H}$ in $\mathrm{T}_{X_k^H}\cMH$ on which $\operatorname{R}_{X_k^H}$ is a diffeomorphism onto its image $\mathcal{V}_{k}^{H}\coloneqq\operatorname{R}_{X_{k}^{H}}(\mathcal{U}_{k}^{H})$; see~\cite[Thm.~4.16 and Cor.~4.17]{boumal2023introduction}.
We restrict the retraction $\operatorname{R}_{X_{k}^{H}}$ to an open neighborhood $\mathcal{U}_{k}^{H}\subset\mathrm{T}_{X_{k}^{H}}\cMH$ of $0_{X_{k}^{H}}$ on which it is a diffeomorphism. The coarse models below are defined for $X^{H}\in\mathcal{V}_{k}^{H}$, while  $\lambda_{k}>0$ and the base point $X^{H}_{k}$ \textit{are fixed}. 

We first recall the definition of the models used by our algorithm:
\begin{itemize}
\item \textbf{Taylor model}
\begin{equation}\label{app:eq:taylor_model}
	T(S^{h}) = f^{h}(X_{k}^{h}) + \innerproduct{\operatorname{grad} f^{h}(X_{k}^{h})}{S^{h}}_{\mathrm{F}}, \qquad S^{h}\in\mathrm{T}_{X_{k}^{h}}\cMh,
\end{equation}
where $\langle A,B\rangle_{\mathrm{F}} \coloneqq \trace(A^\top B)$ is the Frobenius inner product. 
\item \textbf{Coarse model}
\begin{equation}\label{app:eq:psiHk}
	\psi^{H}_{k}(X^{H}) = f^{H}(X^{H}) + \innerproduct{\operatorname{R}^{-1}_{X^{H}_{k}}(X^{H})}{V^{H}_{k}}_{\mathrm{F}},
\end{equation}
where $ V_{k}^{H} $ is defined as
\begin{equation}\label{app:eq:vHk}
	V_{k}^{H} = \widetilde{\mathcal{I}}_{h}^{H} \bigl(\operatorname{grad} f^{h}(X_{k}^{h}) \bigr) - \operatorname{grad} f^{H}(X_{k}^{H}).
\end{equation}
\item \textbf{Regularized coarse model}
\begin{equation}\label{app:eq:mHk}
	m_{k}^{H}(X^{H}) = \psi^{H}_{k}(X^{H}) + \frac{1}{2}\lambda_{k} \| \operatorname{R}^{-1}_{X^{H}_{k}}(X^{H}) \|^{2}_{\mathrm{F}}.
\end{equation}
\end{itemize}

\subsection{Euclidean gradient of the pullback of the coarse model}\label{app:pullback_coarse_gradient}
For a fixed $V_k^H\in\mathrm{T}_{X_k^H}\cMH$, the pullback of the coarse model satisfies
\begin{align*}
    \hat\psi_k^H(S^H) &\coloneqq \psi_k^H\bigl(\operatorname{R}_{X_k^H}(S^H)\bigr) = \hat f^H(S^H)+\innerproduct{S^H}{V_k^H}_{\mathrm{F}},
\end{align*}
where $X^H\in\operatorname{R}_{X_k^H}(\mathcal U_k^H)$, $S^H\in\mathcal U_k^H$, and $\hat f^H\coloneqq f^H\circ\operatorname{R}_{X_k^H}$. 

We compute Euclidean gradients on the vector space $\mathrm{T}_{X_k^H}\cMH$ equipped with the inner product $\innerproduct{\cdot}{\cdot}_{\mathrm{F}}$, identifying its tangent spaces canonically with the vector space itself. For every direction $E\in\mathrm{T}_{X_k^H}\cMH$,
\[
    \operatorname{D} \hat\psi_k^H(S^H)[E] = \operatorname{D} \! \hat f^H(S^H)[E]+\innerproduct{E}{V_k^H}_{\mathrm{F}} = \innerproduct{\nabla\hat f^H(S^H)+V_k^H}{E}_{\mathrm{F}}.
\]
Therefore,
\begin{equation}\label{app:eq:egrad_pullback_psiHk}
  \nabla\hat\psi_k^H(S^H)=\nabla\hat f^H(S^H)+V_k^H.
\end{equation}
Equivalently, applying the chain rule,
\[
    \nabla\hat\psi_k^H(S^H)=\Bigl(\operatorname{D} \! \operatorname{R}_{X_k^H}(S^H)\Bigr)^\ast\Bigl[\operatorname{grad}f^H\bigl(\operatorname{R}_{X_k^H}(S^H)\bigr)\Bigr]+V_k^H.
\]
Here, $\bigl(\operatorname{D} \! \operatorname{R}_{X_k^H}(S^H)\bigr)^\ast$ is the adjoint with respect to the Riemannian inner products at $X_k^H$ and $\operatorname{R}_{X_k^H}(S^H)$.

\subsection{Euclidean and Riemannian gradient of \texorpdfstring{$\psi^{H}_{k}(X^{H})$}{psiHk(XH)} when using the orthographic retraction}\label{app:sec:egrad_psiHk}

Let's call $\phi_{1}(X^{H})$ the correction term that appears in $\psi^{H}_{k}(X^{H})$ (see~\eqref{app:eq:psiHk}), namely,
\begin{equation}\label{app:eq:coarse_correction_term}
    \phi_{1}(X^{H}) \coloneqq \innerproduct{\operatorname{R}^{-1}_{X^{H}_{k}}(X^{H})}{V^{H}_{k}}_{\mathrm{F}}.
\end{equation}
Here, we consider the case where $\operatorname{R}^{-1}_{X^{H}_{k}}(\cdot)$ is the \emph{inverse orthographic retraction}\footnote{See, e.g., \cite[\S 3.2.3]{AbsilOseledets:2015} or~\cite[(2.9)]{sutti2021riemannian}.} $\operatorname{P}_{X^{H}_{k}}$, given on its local domain by\footnote{Cf. \Cref{subsec:retraction_choice}.}
\begin{equation}\label{app:eq:def_inV_orth_retr}
    \operatorname{R}^{-1}_{X^{H}_{k}}(X^{H}) = \operatorname{P}_{X^{H}_{k}}(X^{H}-X^{H}_{k}).
\end{equation}
The right-hand side is affine in $X^{H}$, since the base point $X^{H}_{k}$ \textit{is fixed}, and $\operatorname{P}_{X^{H}_{k}}$ is a fixed linear operator\footnote{For brevity, we use the notation $\operatorname{P}_{X^{H}_{k}}$ to indicate $\operatorname{P}_{\mathrm{T}_{X^{H}_{k}}\cMH}$, the orthogonal projector onto $\mathrm{T}_{X^{H}_{k}}\cMH$.}.

For the purpose of computing the Euclidean gradient of~\eqref{app:eq:coarse_correction_term}, we explicitly introduce the smooth ambient extension $\bar\phi_{1}\colon\R^{m\times n}\to\R$, defined by
\begin{equation}\label{app:eq:ambient_ext_phi1}
    \bar\phi_{1}(X) \coloneqq \left\langle\operatorname{P}_{X^{H}_{k}}(X-X^{H}_{k}),V^{H}_{k}\right\rangle_{\mathrm{F}}, \qquad X\in\R^{m\times n}.
\end{equation}
Consequently, $\bar\phi_{1}$ agrees with $\phi_{1}$ on the domain of the inverse orthographic retraction. Moreover, $\bar\phi_{1}$ is affine and therefore smooth on the entire ambient space.

Here, we aim to show that the Euclidean gradient of this extension with respect to $X$ is
\[
    \nabla\bar\phi_{1}(X) = V^{H}_{k}.
\]
By the definition of the Euclidean gradient,
\[
    \D\bar\phi_{1}(X)[E] = \left\langle\nabla\bar\phi_{1}(X),E\right\rangle_{\mathrm{F}}, \qquad \forall E\in\R^{m\times n}.
\]
We recall the definition of the directional derivative:
\begin{equation}\label{app:eq:def_directional_deriv}
    \D\bar\phi_{1}(X)[E] \coloneqq \left.\frac{\mathrm{d}}{\mathrm{d}t}\,\bar\phi_{1}(X+tE)\right|_{t=0}.
\end{equation}
We develop $\bar\phi_{1}(X+tE)$ using~\eqref{app:eq:ambient_ext_phi1}:
\begin{align*}
    \bar\phi_{1}(X+tE) &= \left\langle\operatorname{P}_{X^{H}_{k}}(X+tE-X^{H}_{k}),V^{H}_{k}\right\rangle_{\mathrm{F}}\\
    &= \left\langle\operatorname{P}_{X^{H}_{k}}(X-X^{H}_{k}),V^{H}_{k}\right\rangle_{\mathrm{F}} + t\left\langle\operatorname{P}_{X^{H}_{k}}(E),V^{H}_{k}\right\rangle_{\mathrm{F}}.
\end{align*}
Then, the derivative of $\bar\phi_{1}(X+tE)$ with respect to $t$ is
\[
    \frac{\mathrm{d}}{\mathrm{d}t}\,\bar\phi_{1}(X+tE) = \left\langle\operatorname{P}_{X^{H}_{k}}(E),V^{H}_{k}\right\rangle_{\mathrm{F}}.
\]

Now, to identify the gradient of $\bar\phi_{1}$, we need to isolate $E$. To do this, we need to ``carry $\operatorname{P}_{X^{H}_{k}}$ over $V^{H}_{k}$''. This is achieved through the adjoint operator of $\operatorname{P}_{X^{H}_{k}}$, viewed as an operator on the ambient matrix space:
\[
    \left\langle\operatorname{P}_{X^{H}_{k}}(E),V^{H}_{k}\right\rangle_{\mathrm{F}} = \left\langle E,\operatorname{P}_{X^{H}_{k}}^\ast(V^{H}_{k})\right\rangle_{\mathrm{F}}.
\]
Therefore, we identify the gradient as
\[
    \nabla\bar\phi_{1}(X) = \operatorname{P}_{X^{H}_{k}}^\ast(V^{H}_{k}).
\]
Since $\operatorname{P}_{X^{H}_{k}}$ is an orthogonal projector, it is self-adjoint. Then,
\[
    \nabla\bar\phi_{1}(X) = \operatorname{P}_{X^{H}_{k}}(V^{H}_{k}).
\]
Finally, since $V^{H}_{k}\in\mathrm{T}_{X^{H}_{k}}\cMH$, we simply have
\begin{equation}\label{app:eq:egrad_linear_correction_orth_retr}
    \nabla\bar\phi_{1}(X) = V^{H}_{k}.
\end{equation}
Note that the resulting gradient is independent of the variable $X$, since $V^{H}_{k}$ stays fixed when differentiating with respect to $X$.

Now, let $\bar f^{H}$ be a smooth ambient extension of $f^{H}$ to an ambient neighborhood of the points under consideration, and define the corresponding extension of $\psi^{H}_{k}$ by
\[
    \bar\psi^{H}_{k}(X) \coloneqq \bar f^{H}(X)+\bar\phi_{1}(X), \qquad X \in \mathbb{R}^{m \times n}.
\]
Wrapping up this calculation, the Euclidean gradient of the ambient extension of $\psi^{H}_{k}$ is
\begin{equation}\label{app:eq:egrad_psiHk_orth_retr}
    \nabla\bar\psi^{H}_{k}(X) = \nabla\bar f^{H}(X)+V^{H}_{k}.
\end{equation}
Then, since the manifold of fixed-rank matrices is an embedded submanifold equipped with the induced Frobenius metric, the Riemannian gradient of $\psi^{H}_{k}$ at $X$ is given by the orthogonal projection of~\eqref{app:eq:egrad_psiHk_orth_retr} onto $\mathrm{T}_{X}\mathcal{M}^{H}$ (see, e.g., \cite[(3.37)]{AMS:2008}), namely,
\begin{align*}
    \operatorname{grad}\psi^{H}_{k}(X) &= \operatorname{P}_{X}\bigl(\nabla\bar\psi^{H}_{k}(X)\bigr)\\
    &= \operatorname{P}_{X}\bigl(\nabla\bar f^{H}(X)\bigr)+\operatorname{P}_{X}(V^{H}_{k})\\
    &= \operatorname{grad}f^{H}(X)+\operatorname{P}_{X}(V^{H}_{k}).
\end{align*}
Thus,
\begin{equation}\label{app:eq:rgrad_psiHk_orth_retr}
    \operatorname{grad}\psi^{H}_{k}(X) = \operatorname{grad}f^{H}(X)+\operatorname{P}_{X}(V^{H}_{k}).
\end{equation}
In particular, at the base point $X = X^{H}_{k}$, since $\operatorname{P}_{X^{H}_{k}}(V^{H}_{k})=V^{H}_{k}$, we obtain
\[
    \operatorname{grad}\psi^{H}_{k}(X^{H}_{k}) = \operatorname{grad}f^{H}(X^{H}_{k})+V^{H}_{k}.
\]

\subsection{Euclidean gradient of the regularization term}

Next, we calculate the Euclidean gradient with respect to $X$ of the regularization term in~\eqref{app:eq:mHk}, omitting the constant prefactor $\lambda_{k}/2$, i.e.,
\[
    \phi_{2}(X) \coloneqq \|\operatorname{R}^{-1}_{X^{H}_{k}}(X)\|^{2}_{\mathrm{F}} = \innerproduct{\operatorname{R}^{-1}_{X^{H}_{k}}(X)}{\operatorname{R}^{-1}_{X^{H}_{k}}(X)}_{\mathrm{F}}.
\]
As before, we explicitly introduce the smooth ambient extension $\bar\phi_{2}\colon\R^{m\times n}\to\R$, defined by
\begin{equation}\label{app:eq:ambient_ext_phi2}
    \bar\phi_{2}(X) \coloneqq \|\operatorname{P}_{X^{H}_{k}}(X-X^{H}_{k})\|_{\mathrm{F}}^{2}, \qquad X\in\R^{m\times n}.
\end{equation}
Since the Riemannian metric is induced by the Frobenius inner product, $\bar\phi_{2}$ agrees with $\phi_{2}$ on the domain of the inverse orthographic retraction. Moreover, $\bar\phi_{2}$ is a quadratic function on the ambient space and is therefore smooth.

The directional derivative is
\[
    \D\bar\phi_{2}(X)[E] \coloneqq \left.\frac{\mathrm{d}}{\mathrm{d}t}\,\bar\phi_{2}(X+tE)\right|_{t=0}, \qquad \forall E\in\R^{m\times n}.
\]
Developing $\bar\phi_{2}(X+tE)$,
\begin{align*}
    \bar\phi_{2}(X+tE) &= \|\operatorname{P}_{X^{H}_{k}}(X-X^{H}_{k})+t\operatorname{P}_{X^{H}_{k}}(E)\|_{\mathrm{F}}^{2}\\
    &= \|\operatorname{P}_{X^{H}_{k}}(X-X^{H}_{k})\|_{\mathrm{F}}^{2} + 2t\left\langle\operatorname{P}_{X^{H}_{k}}(X-X^{H}_{k}),\operatorname{P}_{X^{H}_{k}}(E)\right\rangle_{\mathrm{F}} + t^{2}\|\operatorname{P}_{X^{H}_{k}}(E)\|_{\mathrm{F}}^{2}.
\end{align*}
Then,
\[
    \frac{\mathrm{d}}{\mathrm{d}t}\,\bar\phi_{2}(X+tE) = 2\left\langle\operatorname{P}_{X^{H}_{k}}(E),\operatorname{P}_{X^{H}_{k}}(X-X^{H}_{k})\right\rangle_{\mathrm{F}} + 2t\|\operatorname{P}_{X^{H}_{k}}(E)\|_{\mathrm{F}}^{2}.
\]
Evaluating at $t=0$, we obtain
\[
    \D\bar\phi_{2}(X)[E] = 2\left\langle\operatorname{P}_{X^{H}_{k}}(E),\operatorname{P}_{X^{H}_{k}}(X-X^{H}_{k})\right\rangle_{\mathrm{F}} = 2\left\langle E,\operatorname{P}_{X^{H}_{k}}(X-X^{H}_{k})\right\rangle_{\mathrm{F}},
\]
where in the last equality we used the fact that $\operatorname{P}_{X^{H}_{k}}$ is self-adjoint and idempotent, i.e., $\operatorname{P}_{X^{H}_{k}}^{\ast}=\operatorname{P}_{X^{H}_{k}}$ and $\operatorname{P}_{X^{H}_{k}}^{2}=\operatorname{P}_{X^{H}_{k}}$.

Identifying the Euclidean gradient, we obtain
\begin{equation}\label{app:eq:egrad_quadratic_regularization_orth_retr}
    \nabla\bar\phi_{2}(X) = 2\operatorname{P}_{X^{H}_{k}}(X-X^{H}_{k}) = 2\operatorname{R}^{-1}_{X^{H}_{k}}(X).
\end{equation}
The last equality holds for $X$ in the domain of the inverse orthographic retraction.

\subsection{Euclidean and Riemannian gradients of the regularized coarse model}\label{app:sec:egrad_mHk}

Using the ambient extensions introduced above, we define the corresponding smooth ambient extension of $m_{k}^{H}$ by
\[
    \bar m_{k}^{H}(X) \coloneqq \bar\psi_{k}^{H}(X)+\frac{1}{2}\lambda_{k} \, \bar\phi_{2}(X).
\]
Since $\lambda_{k}$ is fixed when differentiating with respect to $X$, its Euclidean gradient is
\[
    \nabla\bar m_{k}^{H}(X) = \nabla\bar\psi_{k}^{H}(X)+\frac{1}{2}\lambda_{k}\nabla\bar\phi_{2}(X).
\]
Combining~\eqref{app:eq:egrad_psiHk_orth_retr} and~\eqref{app:eq:egrad_quadratic_regularization_orth_retr}, we obtain
\begin{equation}\label{app:eq:egrad_mHk_orth_retr}
    \nabla\bar m_{k}^{H}(X) = \nabla\bar f^{H}(X)+V_{k}^{H}+\lambda_{k}\operatorname{P}_{X^{H}_{k}}(X-X^{H}_{k}).
\end{equation}

Then, the Riemannian gradient of $m_{k}^{H}$ at $X$ is obtained by projecting its ambient Euclidean gradient onto $\mathrm{T}_{X}\mathcal{M}_{r}^{H}$:
\[
    \operatorname{grad}m_{k}^{H}(X) = \operatorname{P}_{X}\bigl(\nabla\bar m_{k}^{H}(X)\bigr).
\]
Therefore,
\begin{equation}\label{app:eq:rgrad_mHk_orth_retr}
    \operatorname{grad}m_{k}^{H}(X) = \operatorname{grad}f^{H}(X)+\operatorname{P}_{X}\left(V_{k}^{H}+\lambda_{k}\operatorname{P}_{X^{H}_{k}}(X-X^{H}_{k})\right).
\end{equation}
Equivalently, using~\eqref{app:eq:rgrad_psiHk_orth_retr},
\[
    \operatorname{grad}m_{k}^{H}(X) = \operatorname{grad}\psi_{k}^{H}(X)+\lambda_{k}\operatorname{P}_{X}\bigl(\operatorname{R}^{-1}_{X^{H}_{k}}(X)\bigr).
\]
The outer projection $\operatorname{P}_{X}$ is needed because $\operatorname{R}^{-1}_{X^{H}_{k}}(X)$ belongs to $\mathrm{T}_{X^{H}_{k}}\mathcal{M}^{H}$, whereas the Riemannian gradient at $X$ must belong to $\mathrm{T}_{X}\mathcal{M}^{H}$.

Finally, at the base point $X = X^{H}_{k}$, the inverse retraction vanishes and $\operatorname{P}_{X^{H}_{k}}(V_{k}^{H})=V_{k}^{H}$. Hence,
\[
    \operatorname{grad}m_{k}^{H}(X^{H}_{k}) = \operatorname{grad}\psi_{k}^{H}(X^{H}_{k}) = \operatorname{grad}f^{H}(X^{H}_{k})+V_{k}^{H}.
\]

\subsection{First-order coherence in the case of the orthographic retraction}\label{app:first-order_coherence}

In this section, we show the first-order coherence property at the base point $X^{H}_{k}$ for every $S^{H}\in\mathrm{T}_{X^{H}_{k}}\mathcal{M}^{H}$, with $S^{h}=\widetilde{\cI}_{H}^{h}S^{H}$, under the adjointness condition for the tangent transfer operators $\widetilde{\cI}_{h}^{H}=(\widetilde{\cI}_{H}^{h})^{\ast}$, namely,
\begin{equation}\label{app:eq:first_order_coherence}
    \innerproduct{\operatorname{grad} \psi^{H}_{k}(X^{H}_{k})}{S^{H}}_{\mathrm{F}} = \innerproduct{\operatorname{grad} f^{h}(X^{h}_{k})}{S^{h}}_{\mathrm{F}},
\end{equation}
We first develop the left-hand side of~\eqref{app:eq:first_order_coherence}:
\begin{align*}
    \innerproduct{\operatorname{grad} \psi^{H}_{k}(X^{H}_{k})}{S^{H}}_{\mathrm{F}} &\stackrel{\eqref{app:eq:rgrad_psiHk_orth_retr}}{=} \innerproduct{\operatorname{grad} f^{H}(X^{H}_{k}) + V^{H}_{k}}{S^{H}}_{\mathrm{F}} \\
    &\stackrel{\eqref{app:eq:vHk}}{=} \innerproduct{\cancel{\operatorname{grad} f^{H}(X^{H}_{k})} + \widetilde{\cI}_{h}^{H} \bigl( \operatorname{grad} f^{h}(X^{h}_{k}) \bigr) \, \cancel{- \operatorname{grad} f^{H}(X^{H}_{k})}}{S^{H}}_{\mathrm{F}} \\
    &= \innerproduct{\widetilde{\cI}_{h}^{H} \bigl( \operatorname{grad} f^{h}(X^{h}_{k}) \bigr)}{S^{H}}_{\mathrm{F}}.
\end{align*}
For the right-hand side of~\eqref{app:eq:first_order_coherence}, using $ S^{h} = \widetilde{\cI}_{H}^{h} S^{H} $, and the adjointness of $ \widetilde{\cI}_{H}^{h} $ with respect to the inner product $ \langle \cdot, \, \cdot \rangle_{\mathrm{F}} $, we obtain
\begin{align*}
    \innerproduct{\operatorname{grad} f^{h}(X^{h}_{k})}{S^{h}}_{\mathrm{F}} &= \innerproduct{\operatorname{grad} f^{h}(X^{h}_{k})}{\widetilde{\cI}_{H}^{h} S^{H}}_{\mathrm{F}} \\
    &= \innerproduct{(\widetilde{\cI}_{H}^{h})^{\ast} \bigl( \operatorname{grad} f^{h}(X^{h}_{k}) \bigr)}{S^{H}}_{\mathrm{F}}.
\end{align*}
To recap, we got, for the LHS and the RHS of \eqref{app:eq:first_order_coherence}, respectively,
\begin{align*}
	\innerproduct{\operatorname{grad} \psi^{H}_{k}(X^{H}_{k})}{S^{H}}_{\mathrm{F}} &= \innerproduct{\widetilde{\cI}_{h}^{H} \bigl( \operatorname{grad} f^{h}(X^{h}_{k}) \bigr)}{S^{H}}_{\mathrm{F}}, \\ &\qquad \text{\underline{and}} \qquad \innerproduct{\operatorname{grad} f^{h}(X^{h}_{k})}{S^{h}}_{\mathrm{F}} = \innerproduct{(\widetilde{\cI}_{H}^{h})^{\ast} \bigl( \operatorname{grad} f^{h}(X^{h}_{k}) \bigr)}{S^{H}}_{\mathrm{F}}.
\end{align*}
Clearly, the left- and right-hand sides of~\eqref{app:eq:first_order_coherence} are equal if $ \widetilde{\cI}_{h}^{H} = (\widetilde{\cI}_{H}^{h})^{\ast} $, that is, if the tangent restriction operator is the adjoint of the tangent prolongation operator with respect to the inner products on the fine and coarse tangent spaces. In the real case, when the tangent spaces are represented in compatible orthonormal coordinates, and the corresponding inner products are the standard Euclidean ones, this condition reduces to $ \widetilde{\cI}_{h}^{H} = (\widetilde{\cI}_{H}^{h})\tr $, or, equivalently, $ \widetilde{\cI}_{H}^{h} = (\widetilde{\cI}_{h}^{H})\tr $.

